%% file: ms.tex
\documentclass[final, 5p, twocolumn, 10pt, times]{elsarticle}

\input{preamble}
\begin{document}
\begin{frontmatter}
\title{Industrial scale large eddy simulations (LES) with adaptive octree meshes using immersogeometric analysis}

\author[isuMechEAddress]{Kumar Saurabh\fnref{equalContrib}}
\ead{maksbh@iastate.edu}

\author[isuMechEAddress]{Boshun Gao\fnref{equalContrib}}
\ead{boshun@iastate.edu}

\author[utahAddress]{Milinda Fernando}
\ead{milinda@cs.utah.edu}

\author[utahAddress]{Songzhe Xu}
\ead{songzhex@sci.utah.edu}

\author[isuMechEAddress,isuMathAddress]{Makrand A. Khanwale}
\ead{khanwale@iastate.edu}

\author[isuMechEAddress]{Biswajit Khara}
\ead{bkhara@iastate.edu}

\author[isuMechEAddress]{Ming-Chen Hsu}
\ead{jmchsu@iastate.edu}

\author[isuMechEAddress]{Adarsh Krishnamurthy}
\ead{adarsh@iastate.edu}

\author[utahAddress]{Hari Sundar}
\ead{hari@cs.utah.edu}

\author[isuMechEAddress]{Baskar Ganapathysubramanian \corref{correspondingAuthor}}
\ead{baskarg@iastate.edu}

\address[isuMechEAddress]{Department of Mechanical Engineering, Iowa State University, Ames, Iowa 50010, USA}
\address[isuMathAddress]{Department of Mathematics, Iowa State University, Ames, Iowa 50010, USA}
\address[utahAddress]{School of Computing, The University of Utah, Salt Lake City, Utah 84112, USA}

\cortext[correspondingAuthor]{Corresponding author}
\fntext[equalContrib]{These authors contributed equally}

\begin{abstract}
We present a variant of the immersed boundary method integrated with octree meshes for highly efficient and accurate Large-Eddy Simulations (LES) of flows around complex geometries. We demonstrate the scalability of the proposed method up to $\mathcal{O}(32K)$ processors. This is achieved by (a) rapid in-out tests; (b) adaptive quadrature for an accurate evaluation of forces; (c) tensorized evaluation during matrix assembly. We showcase this method on two non-trivial applications: accurately computing the drag coefficient of a sphere across Reynolds numbers $1-10^6$ encompassing the \emph{drag crisis} regime; simulating flow features across a semi-truck for investigating the effect of platooning on efficiency.
\end{abstract}


\end{frontmatter}

\section*{Highlights}
\begin{itemize}
    \item Deployed VMS with weak BCs on massively parallel octree-based adaptive meshes.     
    \item Developed octree parallelization, efficient matrix assembly, and rapid in-out tests. 
    \item Demonstrated the ability to capture drag crisis without any wall treatment.
    \item Demonstrated scalability of the framework up to $\mathcal{O} (32 K)$  processors.
    \item Deployed for industrial scale study of the platooning effect of semi-trucks.
\end{itemize}
\section{Introduction}
\input{introduction}


\section{Mathematical preliminaries: Variational treatment of IMGA}
\label{Sec:Formulation}

The fundamentals of immersogeometric fluid-flow analysis consist of three main components. The flow physics is formulated using a variational multiscale (VMS) method for incompressible flows ~\citep{Hughes00a,Bazilevs07b}. To capture the flow domain geometry accurately, adaptively refined quadrature rules are used in the \Intercepted elements, without modifying the background mesh~\citep{Duester:08.1,xu2016tetrahedral}. Finally, the Dirichlet boundary conditions on the surface of the immersed objects are enforced weakly in the sense of Nitsche's method~\citep{Nitsche:70.1,Bazilevs07c}. We briefly discuss each of these components next.

\subsection{Variational Multiscale Formulation (VMS)}
The VMS approach has been successfully utilized to model flow physics across a range of applications~\citep{gravemeier2007variational,xu2019immersogeometric,bazilevs2010large,koobus2004variational,loewe2012projection}. It is considered to be an LES type approach, and uses variational projections in place of the traditional filtered equations in LES. The method is derived completely from the incompressible Navier--Stokes and does not employ any eddy viscosity. See \cite{Bazilevs07b} for a detailed derivation of the VMS method.


The VMS discretization of the Navier--Stokes equations is stated as: Find fluid velocity $\mathbf{u}^h$ and pressure $p^h$ such that for all test functions $\mathbf{w}^h$ and $q^h$ (defined in appropriate function spaces):\\
\begin{align}\label{ALE-VMS-form}
B^{\text{VMS}}\left(\{\mathbf{w}^h,q^h\},\{\mathbf{u}^h,p^h\}\right)-F^{\text{VMS}}\left(\{\mathbf{w}^h,q^h\}\right)=0\text{ ,}
\end{align}

\noindent where,
\begin{align}\label{vms}
\begin{split}
&B^{\text{VMS}}\left(\{\mathbf{w}^h,q^h\},\{\mathbf{u}^h,p^h\}\right) = \\
&\quad \int_{\Omega}\mathbf{w}^h\cdot\rho\left(\frac{\partial\mathbf{u}^h}{\partial t}
+ \mathbf{u}^h\cdot\pmb{\nabla} \mathbf{u}^h\right) d\Omega \\
& \quad + \int_{\Omega}\pmb{\varepsilon}(\mathbf{w}^h):\pmb{\sigma} \left(\mathbf{u}^h, p^h \right) d\Omega + \int_{\Omega}q^h\,\pmb{\nabla} \cdot\mathbf{u}^h d\Omega \\
&\quad - \sum_{e}\int_{\Omega^e\cap\Omega}\left(\mathbf{u}^h\cdot\pmb{\nabla} \mathbf{w}^h+\frac{\pmb{\nabla}q^h}{\rho}\right)\cdot\mathbf{u}' d\Omega\\
&\quad - \sum_{e}\int_{\Omega^e\cap\Omega} p'\,\pmb{\nabla} \cdot\mathbf{w}^h d\Omega\\
&\quad + \sum_{e}\int_{\Omega^e\cap\Omega}\mathbf{w}^h\cdot(\mathbf{u}'\cdot\pmb{\nabla} \mathbf{u}^h) d\Omega\\
&\quad - \sum_{e}\int_{\Omega^e\cap\Omega}\frac{{\mathbf{\pmb{\nabla}w}}^h}{\rho}:\left(\mathbf{u}'\otimes\mathbf{u}'\right) d\Omega\\
&\quad + \sum_{e}\int_{\Omega^e\cap\Omega}\left(\mathbf{u}'\cdot\pmb{\nabla}\mathbf{w}^h\right)\overline{\tau}\cdot\left(\mathbf{u}'\cdot\pmb{\nabla} \mathbf{u}^h\right) d\Omega\text{ ,}
\end{split}
\end{align}
\noindent and $F^{\text{VMS}}$ is the force term. The fine scale velocity, $\mathbf{u}'$, and pressure, $p'$, are defined as\\
\begin{align}
\label{uprime}
\mathbf{u}' = -\tau_\text{M}\left(\rho\left(\frac{\partial\mathbf{u}^h}{\partial t}+\mathbf{u}^h\cdot\pmb{\nabla} \mathbf{u}^h-\mathbf{f}\right) - \pmb{\nabla} \cdot \pmb{\sigma}\left(\mathbf{u}^h, p^h \right)\right),
\end{align}
and 
\begin{align}\label{pprime}
p' = -\rho\,\tau_\text{C}\,\pmb{\nabla} \cdot\mathbf{u}^h \text{ .}
\end{align}
%
In the above equations,  $\Omega^e$ represents the disjoint elements, such that $\Omega \subset \cup_e  \Omega^e$, $\rho$ is the density of the fluid, and  $\pmb{\sigma}$ and $\pmb{\varepsilon}$ are the stress and strain-rate tensors, respectively.
The terms integrated over element interiors may be interpreted both as stabilization and as a turbulence model~\citep{Brooks82a,Tezduyar91c,Tezduyar00a,Hughes01a,Bazilevs07b,Hsu09a}. $\tau_\text{M}, \tau_\text{C}$ and $\overline{\tau}$ are the stabilization parameters. Their detailed expression used in this work can be found in~\citet{xu2016tetrahedral}.

\subsection{Adaptive Quadrature}
The geometric boundary of immersed object creates complex, discontinuous integration in the \Intercepted elements. In order to ensure the geometrically accurate integration of the volume integrals, we use a sub-cell based quadrature scheme~\citep{Duester:08.1,kamensky2015immersogeometric}. The basic idea is to increase the number of quadrature points in the \Intercepted elements and perform accurate evaluation by considering the quadrature points that are only outside the body. We discuss an efficient approach to do this later in \secref{sec: VGP}.

\subsection{Weak Enforcement of Boundary Conditions}
The standard way of imposing Dirichlet boundary conditions is to enforce them strongly by ensuring that they are satisfied by all trial solution functions, which is not feasible in immersed methods. Instead, the strong enforcement is replaced by weakly enforced Dirichlet boundary conditions~\citep{Bazilevs07c,BaMiCaHu07,Bazilevs09d}. The semi-discrete problem can be stated as follows: 
Find fluid velocity $\mathbf{u}^h$ and pressure $p^h$ such that for all test functions $\mathbf{w}^h$ and $q^h$:
%
\begin{align}
\label{eq:weakbc}
\begin{split}
&B^{\text{VMS}}\left(\{\mathbf{w}^h,q^h\},\{\mathbf{u}^h,p^h\}\right)-F^{\text{VMS}}\left(\{\mathbf{w}^h,q^h\}\right) \\
&\quad \quad 
-\int_{\Gamma^{D}}\mathbf{w}^h\cdot\left(-p^h\,\mathbf{n}+2\mu\,\pmb{\varepsilon}(\mathbf{u}^h)\,\mathbf{n}\right) d\Gamma 
\\&\quad \quad 
-\int_{\Gamma^{D}}\left(q^h\,\mathbf{n}+2\mu\,\pmb{\varepsilon}(\mathbf{w}^h)\,\mathbf{n}\right)\cdot\left(\mathbf{u}^h - \mathbf{g}\right) d\Gamma\\
&\quad -\int_{\Gamma^{D,-}}\mathbf{w}^h\cdot\rho\left(\mathbf{u}^h\cdot\mathbf{n}\right)\left(\mathbf{u}^h-\mathbf{g}\right) d\Gamma
\\&\quad \quad 
+\int_{\Gamma^{D}}\tau^B_{\text{TAN}}\left(\mathbf{w}^h-\left(\mathbf{w}^h\cdot\mathbf{n}\right)\mathbf{n}\right)
\\&\quad \quad \quad \quad \quad
\cdot\left(\left(\mathbf{u}^h - \mathbf{g}\right)-\left(\left(\mathbf{u}^h - \mathbf{g}\right)\cdot\mathbf{n}\right)\mathbf{n}\right) d\Gamma
\\&\quad \quad 
+\int_{\Gamma^{D}}\tau^B_{\text{NOR}}\left(\mathbf{w}^h\cdot\mathbf{n}\right)\left(\left(\mathbf{u}^h - \mathbf{g}\right)\cdot\mathbf{n}\right) d\Gamma~=~0\text{ ,}
\end{split}
\end{align}
%
\noindent where
$\mathbf{n}$ is the outward unit normal, $\mu$ is the dynamic viscosity, $\Gamma^{D}$ is the Dirichlet boundary that may cut through element interiors, $\Gamma^{D,-}$ is the \emph{inflow} part of $\Gamma^{D}$, on which $\mathbf{u}^h\cdot\mathbf{n} < 0$, $\mathbf{g}$ is the prescribed velocity on  $\Gamma^{D}$, and $\tau_{\text{TAN}}^B$ and $\tau_{\text{NOR}}^B$ are stabilization parameters.

\begin{remark}
\added{The stabilization parameter for the weak imposition of the boundary condition is chosen such that $\tau_\text{TAN}^{B} = \tau_\text{NOR}^{B} = \frac{C_b}{Re h_b} > 4$, where $h_b$ is computed as the maximum projected distance from the outside node(s) (i.e. fluid node) of the background \Intercepted elements to the triangulated surface. Complex surface geometries can sometimes result in 'sliver cut' elements, which exhibit $h_b \rightarrow 0$. This can lead to arbitrarily large values of the stabilization parameter, resulting in accuracy and convergence issues ~\cite{de2018note}. To circumvent this issue, the lower limit on $h_b$ is set to  $0.01 h$, where $h$ is the grid spacing of background element.}
\end{remark}

The weak enforcement of the boundary condition in IMGA is particularly attractive as this alleviates the need for a body-fitted mesh. The additional Nitsche terms (the third to last terms on the left-hand side of \eqnref{eq:weakbc}) are formulated independently of the mesh. The only additional overhead is a separate discretization of the domain boundary whose position in the \Intercepted elements can be determined. 

\begin{remark}
 \added{We use a fully implicit Crank Nicholson time stepping scheme. All results are reported using linear basis functions, unless explicitly stated otherwise.}
\end{remark}


\section{Algorithmic developments}\label{sec: Algorithm}
\input{Algorithms}

\section{Benchmark results}\label{sec: Results}
\input{Result}

\section{Conclusion}\label{sec:Conclusion}
\input{Conclusion}
\section*{Acknowledgment}
KS, Boshun G, BK, MAK, Baskar G were partly supported by NSF 1855902 and NSF 1935255. AK was partly supported by NSF 1644441 and NSF 1750865. HS was partly funded by NSF  OAC-1808652. Computing allocation through Directors discretionary award (\Frontera), XSEDE CTS110007 (\Stampede) and Iowa State University (\Nova) are gratefully acknowledged


\bibliographystyle{./bibliography/elsarticle-num-names}
\bibliography{./ms}




\label{sec:ae}
\input{ae.tex}


\end{document}

%% file: preamble.tex
\usepackage{xcolor} 
\usepackage[noend]{algpseudocode}
\usepackage{array}
\usepackage{amsfonts}
\usepackage{graphicx}
\usepackage{bm,mathrsfs}
\usepackage{epstopdf}
\usepackage{amssymb}
\usepackage{tikz,pgfplots}
\usetikzlibrary{snakes,arrows,shapes,trees}
\usepackage{amssymb,amsmath,amsthm}
\usepackage{amsopn}
\usepackage{listings}
\usepackage{adjustbox}
\usepackage{longtable}
\usepackage{multirow}
\usepackage{hyperref}
\usepackage{amssymb}
\usepackage{pgfplots}
\usepackage{pgfplotstable}
\usepackage{grffile}
\usetikzlibrary{arrows,shapes,plotmarks}
\usepgfplotslibrary{groupplots}
\usetikzlibrary{matrix}
\PassOptionsToPackage{normalem}{ulem}
\usepackage{ulem}

\providecolor{added}{rgb}{0,0,0}
\providecolor{deleted}{rgb}{1,0,0}
\newcommand{\added}[1]{{\color{added}{}#1}}

\newcommand{\cref}[2]{\hyperref[#2]{#1~\ref*{#2}}}
\newcommand{\figref}[1]{\hyperref[#1]{Fig.~\ref*{#1}}}
\newcommand{\secref}[1]{\hyperref[#1]{Section~\ref*{#1}}}
\newcommand{\tabref}[1]{\hyperref[#1]{Table~\ref*{#1}}}
\newcommand{\eqnref}[1]{\hyperref[#1]{Eq.~(\ref*{#1})}}
\hypersetup{
	colorlinks,
	linkcolor={blue!50!black},
	citecolor={blue!50!black},
	urlcolor={blue!80!black},
	anchorcolor = {blue!80!black},
	filecolor = {blue!80!black},
	menucolor = {blue!80!black},
	runcolor = {blue!80!black}
}

\pgfplotsset{compat=1.8}
\usepackage{soul}
\usepackage{rotating}
\usepackage{url}
\usepackage{algorithm}
\usepackage{enumitem}
\usepackage{subcaption}

\newcommand{\Tensor}[1]{\underline{\underline{\mathbf{#1}}}}

\newcommand{\Vector}[1]{\mathbf{#1}}

\newcommand{\Intercepted}{\textsc{Intercepted }}
\newcommand{\In}{\textsc{In }}
\newcommand{\Out}{\textsc{Out }}

\newcommand{\Stampede}{\href{https://www.tacc.utexas.edu/systems/stampede2}{Stampede2}}
\newcommand{\Nova}{\href{https://www.hpc.iastate.edu/guides/nova}{Nova}}
\newcommand{\Frontera}{\href{https://frontera-portal.tacc.utexas.edu/}{Frontera}}

\newcommand{\petsc}{\href{https://www.mcs.anl.gov/petsc/}{PETSc}}

\newdimen\HilbertLastX
\newdimen\HilbertLastY
\newcounter{HilbertOrder}

\def\DrawToNext#1#2{%
  \advance \HilbertLastX by #1
  \advance \HilbertLastY by #2
  \pgfpathlineto{\pgfqpoint{\HilbertLastX}{\HilbertLastY}}
}

\def\Hilbert[#1,#2,#3,#4,#5,#6,#7,#8] {
  \ifnum\value{HilbertOrder} > 0%
  \addtocounter{HilbertOrder}{-1}
  \Hilbert[#5,#6,#7,#8,#1,#2,#3,#4]
  \DrawToNext {#1} {#2}
  \Hilbert[#1,#2,#3,#4,#5,#6,#7,#8]
  \DrawToNext {#5} {#6}
  \Hilbert[#1,#2,#3,#4,#5,#6,#7,#8]
  \DrawToNext {#3} {#4}
  \Hilbert[#7,#8,#5,#6,#3,#4,#1,#2]
  \addtocounter{HilbertOrder}{1}
  \fi
}

\def\hilbert((#1,#2),#3){%
  \advance \HilbertLastX by #1
  \advance \HilbertLastY by #2
  \pgfpathmoveto{\pgfqpoint{\HilbertLastX}{\HilbertLastY}}
  \setcounter{HilbertOrder}{#3}
  \Hilbert[1mm,0mm,-1mm,0mm,0mm,1mm,0mm,-1mm]
  \pgfusepath{stroke}%
}

\ifpdf
  \DeclareGraphicsExtensions{.eps,.pdf,.png,.jpg}
\else
  \DeclareGraphicsExtensions{.eps}
\fi

\emergencystretch=\maxdimen
\definecolor{cpu3}{HTML}{F44336}
\definecolor{cpu4}{HTML}{2196F3}
\definecolor{cpu1}{HTML}{4CAF50}
\definecolor{cpu2}{HTML}{FFC107}
\definecolor{gpu3}{HTML}{EF9A9A}
\definecolor{gpu4}{HTML}{90CAF9}
\definecolor{gpu1}{HTML}{A5D6A7}
\definecolor{gpu2}{HTML}{FFE082}

\definecolor{cpu5}{HTML}{9932CC}

\definecolor{sq_b1}{RGB}{37,52,148}
\definecolor{sq_b2}{RGB}{44,127,184}
\definecolor{sq_b3}{RGB}{65,182,196}
\definecolor{sq_b4}{RGB}{127,205,187}
\definecolor{sq_b5}{RGB}{199,233,180}
\definecolor{sq_b6}{RGB}{255,255,204}

\definecolor{sq_r1}{RGB}{189,0,38}
\definecolor{sq_r2}{RGB}{240,59,32}
\definecolor{sq_r3}{RGB}{253,141,60}
\definecolor{sq_r4}{RGB}{254,178,76}
\definecolor{sq_r5}{RGB}{254,217,118}
\definecolor{sq_r6}{RGB}{255,255,178}

\definecolor{sq_g1}{RGB}{0,104,55}
\definecolor{sq_g2}{RGB}{49,163,84}
\definecolor{sq_g3}{RGB}{120,198,121}
\definecolor{sq_g4}{RGB}{173,221,142}
\definecolor{sq_g5}{RGB}{217,240,163}
\definecolor{sq_g6}{RGB}{255,255,204}

\definecolor{div_c1}{RGB}{230,171,2}
\definecolor{div_c2}{RGB}{102,166,30}
\definecolor{div_c3}{RGB}{231,41,138}
\definecolor{div_c4}{RGB}{117,112,179}
\definecolor{div_c5}{RGB}{217,95,2}
\definecolor{div_c6}{RGB}{27,158,119}
\definecolor{div_c7}{RGB}{215,48,39}

\definecolor{div_d1}{RGB}{215,25,28}
\definecolor{div_d2}{RGB}{253,174,97}
\definecolor{div_d3}{RGB}{255,255,191}
\definecolor{div_d4}{RGB}{171,217,233}
\definecolor{div_d5}{RGB}{44,123,182}

\definecolor{lineclr}{RGB}{0,0,0}
\definecolor{utorange}{RGB}{0,0,255}
\definecolor{utsecblue}{RGB}{255,255,0}
\definecolor{utsecgreen}{RGB}{255,0,0}
\definecolor{red!15}{RGB}{0,255,255}
\definecolor{fillclr5}{RGB}{0,255,0}
\definecolor{fillclr6}{RGB}{255,0,255}
\definecolor{fillclr7}{RGB}{255,255,255}
\definecolor{fillclr8}{RGB}{0,0,0}

\newcommand{\tikzcircle}[2][red,fill=red]{\tikz[baseline=-0.5ex]\draw[#1,radius=#2] (0,0) circle ;}%

\def\drawcubeI(#1,#2,#3,#4,#5){ 
\coordinate (O) at (#1,#2,#3);
\coordinate (A) at (#1,#2+#4,#3);
\coordinate (B) at (#1,#2+#4,#3+#4);
\coordinate (C) at (#1,#2,#3+#4);
\coordinate (D) at (#1+#4,#2,#3);
\coordinate (E) at (#1+#4,#2+#4,#3);
\coordinate (F) at (#1+#4,#2+#4,#3+#4);
\coordinate (G) at (#1+#4,#2,#3+#4);
\draw[#5] (O) -- (C) -- (G) -- (D) -- cycle;
\draw[#5] (O) -- (A) -- (E) -- (D) -- cycle;
\draw[#5] (O) -- (A) -- (B) -- (C) -- cycle;
\draw[#5] (D) -- (E) -- (F) -- (G) -- cycle;
\draw[#5] (C) -- (B) -- (F) -- (G) -- cycle;
\draw[#5] (A) -- (B) -- (F) -- (E) -- cycle;
}

\def\drawcubeII(#1,#2,#3,#4,#5,#6,#7){ 
\coordinate (O) at (#1,#2,#3);
\coordinate (A) at (#1,#2+#4,#3);
\coordinate (B) at (#1,#2+#4,#3+#4);
\coordinate (C) at (#1,#2,#3+#4);
\coordinate (D) at (#1+#4,#2,#3);
\coordinate (E) at (#1+#4,#2+#4,#3);
\coordinate (F) at (#1+#4,#2+#4,#3+#4);
\coordinate (G) at (#1+#4,#2,#3+#4);
\draw[#5,fill=#6,opacity=#7] (O) -- (C) -- (G) -- (D) -- cycle;
\draw[#5,fill=#6,opacity=#7] (O) -- (A) -- (E) -- (D) -- cycle;
\draw[#5,fill=#6,opacity=#7] (O) -- (A) -- (B) -- (C) -- cycle;
\draw[#5,fill=#6,opacity=#7] (D) -- (E) -- (F) -- (G) -- cycle;
\draw[#5,fill=#6,opacity=#7] (C) -- (B) -- (F) -- (G) -- cycle;
\draw[#5,fill=#6,opacity=#7] (A) -- (B) -- (F) -- (E) -- cycle;
}

\def\drawNodes(#1,#2,#3,#4,#5,#6,#7){ 
\foreach \x in {#1,#7,...,#2}{
	\foreach \y in {#3,#7,...,#4}{
		\foreach \z in {#5,#7,...,#6}{
				\draw[fill=red!60] (\x,\y,\z) circle (0.15);
				}
			}
	}				
		
}

\pgfplotsset{
  log x ticks with fixed point/.style={
      xticklabel={
        \pgfkeys{/pgf/fpu=true}
        \pgfmathparse{exp(\tick)}%
        \pgfmathprintnumber[fixed relative, precision=3]{\pgfmathresult}
        \pgfkeys{/pgf/fpu=false}
      }
  },
  log y ticks with fixed point/.style={
      yticklabel={
        \pgfkeys{/pgf/fpu=true}
        \pgfmathparse{exp(\tick)}%
        \pgfmathprintnumber[fixed relative, precision=3]{\pgfmathresult}
        \pgfkeys{/pgf/fpu=false}
      }
  }
}

\makeatletter
\newcommand\resetstackedplots{
\makeatletter
\pgfplots@stacked@isfirstplottrue
\makeatother
\addplot [forget plot,draw=none] coordinates{(48,0) (96,0) (192,0) (384,0) (768,0) (1536,0) (3072,0) (6144,0)};
}
\makeatother

\makeatletter
\newcommand\resetstackedplotsOne{
\makeatletter
\pgfplots@stacked@isfirstplottrue
\makeatother
\addplot [forget plot,draw=none] coordinates{(384,0) (768,0) (1536,0) (3072,0) (6144,0)};
}
\makeatother

\makeatletter
\newcommand\resetstackedplotsTwo{
\makeatletter
\pgfplots@stacked@isfirstplottrue
\makeatother
\addplot [forget plot,draw=none] coordinates{(16,0) (32,0) (64,0) (128,0) (256,0) (512,0) (1024,0) (2048,0) (4096,0) (8192,0) (16384,0) (32768,0)};
}
\makeatother

\makeatletter
\newcommand\resetstackedplotsThree{
\makeatletter
\pgfplots@stacked@isfirstplottrue
\makeatother
\addplot [forget plot,draw=none] coordinates{(2,0) (4,0) (8,0) (16,0) (32,0) (64,0)};
}
\makeatother

\makeatletter
\newcommand\resetstackedplotsFour{
\makeatletter
\pgfplots@stacked@isfirstplottrue
\makeatother
\addplot [forget plot,draw=none] coordinates{(4,0) (8,0) (16,0) (32,0) (64,0)};
}
\makeatother

\makeatletter
\newcommand\resetstackedplotsFive{
\makeatletter
\pgfplots@stacked@isfirstplottrue
\makeatother
\addplot [forget plot,draw=none] coordinates{(1,0) (2,0) (4,0) (8,0) (16,0) (32,0) (64,0) (128,0)};
}
\makeatother

\makeatletter
\newcommand\resetstackedplotsSix{
\makeatletter
\pgfplots@stacked@isfirstplottrue
\makeatother
\addplot [forget plot,draw=none] coordinates{(2,0) (4,0) (8,0) (16,0) (32,0) (64,0) (128,0)};
}
\makeatother

\newcommand{\oTn}{\textsc{o2n}}

\definecolor{armygreen}{rgb}{0.29, 0.33, 0.13}
\definecolor{aurometalsaurus}{rgb}{0.43, 0.5, 0.5}
\definecolor{applegreen}{rgb}{0.55, 0.71, 0.0}
\definecolor{darkgreen}{rgb}{0.0, 0.4, 0.25}

\usepackage{tgpagella}

\newtheorem*{remark}{Remark}
\usepackage{dirtytalk}

\usetikzlibrary{external}

%% file: introduction.tex


There has been a significant demand for the development of fast, scalable numerical methods that can run industrial scale Large Eddy Simulations (LES). An ambitious goal of the community is to perform LES over a complex geometry overnight~\citep{lohner2020overnight}. 
One bottleneck to this goal is creating an analysis-suitable, body-fitted mesh with appropriate refinement for complex geometry, which is time consuming and usually a labor intensive process. This becomes even more challenging for moving bodies, as deforming meshes or remeshing is required at every time step. The main motivation behind immersed boundary methods (IBMs) is to alleviate these time consuming and laborious process. The origin of IBM dates back to 1972 when \citet{peskin1972flow}  utilized it to solve a cardiac blood flow problem on a Cartesian grid. This highlights the major advantage of IBM to perform the complete simulation on a structured grid, and thus avoids the requirement for the grid to conform to the complex shape. Unfortunately, this also made the application of kinematic boundary conditions such as no-slip on the surface of immersed boundaries non-trivial. There have been several developments over the past two decades on the application of correct boundary condition for the immersed cells \citep{mittal2008versatile, fadlun2000combined, tseng2003ghost, kim2001immersed, taira2007immersed, pinelli2010immersed,yang2009smoothing}. Interested readers are referred to the review by \citet{mittal2005immersed}.

In this work, we consider a variant of IBM, known as Immersogeometric Analysis (IMGA), used in the context of Finite Element (FE) and Isogeometric~\cite{Hughes05a} simulations. 
In IMGA, the surface representation of the body in the form of boundary representation (B-rep, NURBS or .stl) is immersed into a non-body-fitted spatial discretization, thereby preserving the exact representation of the immersed geometry. Additionally, the Dirichlet boundary conditions on the immersed geometry surfaces are enforced weakly using Nitsche's method~\citep{Nitsche:70.1}. This variational weakening of the no-slip condition into a Neumann type condition provides a consistent and robust way of enforcing kinematic conditions, especially in the context of FE analysis~\citep{embar2010imposing}. 

IMGA has been deployed by several research groups for a variety of fluid-structure interaction (FSI) simulations, including complex biomedical applications with NURBS \citep{hsu2014fluid,kamensky2015immersogeometric,Xu18id,Wu19fm}, design optimization \citep{wu2017optimizing}, external aerodynamics simulations with tetrahedral meshes \citep{xu2016tetrahedral,Hsu-2016-IMGA,wang2017rapid,Xu19ct}, and moving objects \citep{xu2019immersogeometric,Kamensky20db}. The weak imposition of the no-slip condition has been demonstrated to be numerically very advantageous, especially for flow past complex geometries where steep gradients are produced \citep{Bazilevs07c, zhu2020immersogeometric}. However, challenges remain to practical deployment of IMGA, especially towards the goal of overnight LES. In this work, we identify and resolve the following technical bottlenecks to deploying IMGA for large scale simulations:
\begin{itemize}[itemsep=0pt]
    \item \textit{Simulations on massively parallel adaptive meshes}. In contrast to using unstructured background meshes, we utilize octree-based, parallel adaptive meshes resulting in improved scalability at extreme scales, as well as the ability to efficiently remesh.
    \item \textit{Matrix assembly}: The performance of linear algebra solvers has substantially improved. Several robust optimized libraries have been developed to perform fast numerical linear algebra calculations \cite{balay2019petsc,heroux2003trilinos,tomov2009magma,deconinck2017atlas}; this has made the other parts of the code, specifically matrix assembly a major bottleneck. An analysis of our IMGA solver suggests that up to 70\% of the time is spent in matrix assembly. Substantial improvements are possible by optimizing matrix assembly, which we accomplish using tensorized operations. This is a step towards our intent to transition to matrix-free methods.
    \item \textit{Load imbalance}: The enforcement of kinematic constraints on \Intercepted elements (see \figref{fig: IBMfig}) requires an additional surface and volume integration in those elements. The volume integration must be performed accurately on \textit{only that fraction} of each \Intercepted element that is outside the object. This can be a large fraction of assembly time for complex geometries. We deploy a weighted, dynamic partitioning to ensure load balancing, even with adaptive quadrature.
    \item \In - \Out \textit{test}: Classifying the location of a point in the background mesh with respect to the object (as inside or outside the object) is a quintessential IBM ingredient. We go beyond conventional ray-tracing approaches  \citep{zhu2020immersogeometric,xu2016tetrahedral} (which are computationally expensive) to a more efficient normal based test.
    \item \textit{Reliable numerical methods}: The utility of industrial CFD (Computational Fluid Dynamics) simulations is limited by the choice of appropriate modeling strategies that can accurately predict the region of separation especially for the turbulent cases. Currently, the most widely methods used like Reynolds Averaged Navier--Stokes (RANS) or Detached Eddy Simulation (DES) rely on additional wall treatment to achieve this. Such application specific strategies limit reliability of industrial scale CFD. In this work we demonstrate the ability to predict the drag crisis phenomena \textit{without any special treatment}, with the same method working across six orders of magnitude variation in Reynolds number.
\end{itemize}{}
\input{Plots/IBMFig}
This work is motivated by challenges articulated in the \textit{NASA CFD 2030} \cite{slotnick2014cfd} vision towards the goal of performing overnight LES: a) \textit{\say{Mesh generation and adaptivity continue to be significant bottlenecks in the CFD workflow..}}; b) \textit{\say{The use of CFD ... is severely limited by the inability to accurately and reliably predict turbulent flows with significant regions of separation}}. To the best of our knowledge, this is the first time that IMGA simulations in conjunction with VMS and weak BC on a massively parallel octree-based adaptive mesh has been performed. \figref{fig: IBMfig} shows the representative diagram of IMGA computation on octrees. In a nutshell, our main contributions include: 
\begin{enumerate}[itemsep=0pt]
    \item Deploy the variational formulation of Navier--Stokes with weak enforcement of Dirichlet boundary conditions on adaptive octree based mesh at high Reynolds numbers.
    \item  Deploy adaptive quadrature based schemes for accurate integration of \Intercepted elements.
    \item Demonstrate the ability of the numerical method to capture the drag crisis without any special wall treatment.
    \item Develop a fast normal based in - out test to accurately determine the points in and out of the boundary.
    \item Perform near optimal load balancing by accurately estimating the weight per element (using an empirical model) to account for imbalance in the load arising due to IMGA.
    \item Show scaling of our framework to $\mathcal{O}(32K)$ processors.
    \item Illustrate framework on a complex engineering problem with implications to autonomous vehicles.
\end{enumerate}

The rest of the paper is organized as follows: we begin by giving a brief overview of IMGA and weak imposition of boundary condition in \secref{Sec:Formulation}, the key algorithmic developments for mesh generation, matrix assembly, weighted partitioning and adaptive quadrature are outlined in \secref{sec: Algorithm}; numerical and scaling benchmark results are presented and discussed in \secref{sec: Results}; the framework is then deployed to study the platooning effect on a semi-trailer truck in \secref{sec: truckSimulation}; and concluding remarks and future outlook are made in \secref{sec:Conclusion}.

%% file: Plots/IBMFig.tex
\begin{figure}
\centering
    \begin{tikzpicture}[scale=0.85]
     \fill[cpu1] (0.0,0.0) rectangle (7.5,7.5);
    \fill[cpu2] (0.75,4.5) rectangle (5.25,5.25);
    \fill[cpu2] (1.50,5.25) rectangle (4.5,6.0);
    \fill[cpu2] (0.75,0.75) rectangle (1.5,4.5);
    \fill[cpu2] (4.5,0.75) rectangle (5.25,4.5);
    \fill[cpu2] (0.0,1.5) rectangle (0.75,4.5);
    \fill[cpu2] (5.25,1.5) rectangle (6.0,4.5);
    \fill[cpu2] (1.50,0.0) rectangle (4.5,3.0);
    
    \fill[cpu4] (4.5,1.5) rectangle (1.5,4.5);
    \fill[cpu4] (0.75,2.25) rectangle (1.5,3.75);
    \fill[cpu4] (4.5,2.25) rectangle (5.25,3.75);
    \fill[cpu4] (2.25,4.5) rectangle (3.75,5.25);
    \fill[cpu4] (2.25,0.75) rectangle (3.75,1.5);
    \draw[step=1.5cm, thick] (0,0) grid (7.5,7.5);
    \draw[-,black, thick] (0,0.75) -- (6,0.75);
    \draw[-,black, thick] (0,5.25) -- (6,5.25);
    \draw[-,black, thick] (0.75,0) -- (0.75,6.0);
    \draw[-,black, thick] (5.25,0) -- (5.25,6.0);
    \draw[-,black, thick] (2.25,0) -- (2.25,1.5);
    \draw[-,black, thick] (3.75,0) -- (3.75,1.5);
    \draw[-,black, thick] (2.25,4.5) -- (2.25,6.0);
    \draw[-,black, thick] (3.75,4.5) -- (3.75,6.0);
    
    \draw[-,black, thick] (0,2.25) -- (1.5,2.25);
    \draw[-,black, thick] (0,3.75) -- (1.5,3.75);
    
    \draw[-,black, thick] (4.5,2.25) -- (6,2.25);
    \draw[-,black, thick] (4.5,3.75) -- (6,3.75);

    \draw[red,ultra thick] (3,3 ) circle (2.5cm);
   
    \foreach \x in {0,1.5,...,7.5}{
            \draw [fill=red](\x,0.0) circle (0.1cm);
            \draw [fill=red](\x,7.5) circle (0.1cm);
            \draw [fill=red](\x,6.0) circle (0.1cm);
    }
    \foreach \y in {0,1.5,...,7.5}{
            \draw [fill=red](0.0,\y) circle (0.1cm);
            \draw [fill=red](6.0,\y) circle (0.1cm);
            \draw [fill=red](7.5,\y) circle (0.1cm);
    }
    
    \foreach \x in {1.5,3.0,4.5}
    \foreach \y in {1.5,3.0,4.5}
        \draw [fill=gray](\x,\y) circle (0.1cm);
        
    \foreach \y in {0.0,0.75,...,5.9}{
        \draw [fill=red](0.0,\y) circle (0.1cm);
        \draw [fill=red](0.75,\y) circle (0.1cm);
        \draw [fill=red](5.25,\y) circle (0.1cm);
    }
    \foreach \x in {0.0,0.75,...,5.9}{
        \draw [fill=red](\x,0.0) circle (0.1cm);
        \draw [fill=red](\x,0.75) circle (0.1cm);
        \draw [fill=red](\x,5.25) circle (0.1cm);
    }
    
    \foreach \x in {2.25,3.0,...,4.2}{
        \draw [fill=gray](\x,0.75) circle (0.1cm);
        \draw [fill=gray](\x,5.25) circle (0.1cm);
    }
    \foreach \y in {2.25,3.0,...,4.2}{
        \draw [fill=gray](5.25,\y) circle (0.1cm);
        \draw [fill=gray](0.75,\y) circle (0.1cm);
    }
     \draw [fill=none,color=black,ultra thick](3.0,3.0) circle (0.3cm);
    
    
    
    \end{tikzpicture}
    \caption{Sketch illustrating the IMGA on octree based adaptive mesh. The solid red line represents the boundary of the  immersed object. The elements are classified into 3 types :a) \Out (\textcolor{cpu1}{$\blacksquare$}) : all nodes are outside the body; b) \In (\textcolor{cpu4}{$\blacksquare$}): all nodes are inside the body; c) \Intercepted (\textcolor{cpu2}{$\blacksquare$}) : some of the nodes are inside and some are outside. CG (Continuous Galerkin) nodes are divided into  \In (\tikzcircle[fill=gray]{0.1 cm}) and \Out (\tikzcircle[fill=red]{0.1 cm}) nodes. Hanging nodes are neither classified \In or \Out as they do not correspond to independent degrees of freedom. The nodes which neither belong to the \Intercepted nor \Out elements and are also not the parent of hanging nodes (marked by \tikzcircle[fill=none,color=black,ultra thick]{0.1 cm}) are not solved for in IMGA simulations.}
    \label{fig: IBMfig}
\end{figure}
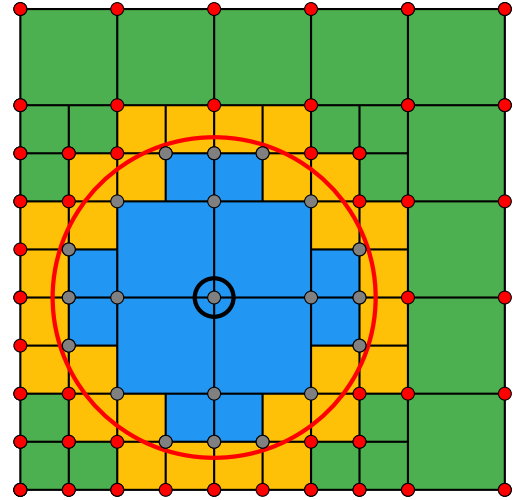{}

%% file: Algorithms.tex


In this section, we highlight some of the key algorithmic developments to accelerate IMGA computations and deploy it on massively parallel computers. Algorithm~\ref{alg:overview} gives an overview of the key steps in IMGA. We start by constructing an octree-based adaptive mesh with a constraint on the relative size of neighboring elements (2:1 balancing). The mesh is refined near the object boundary which is important to capture boundary layer effects. A brief overview of our mesh generation algorithm is given in \secref{sec: octree}. Based on the mesh boundaries, the \textit{.stl} triangles are distributed across the processors, based on its overlap with elements that are part of the background mesh. This step is important to perform fast \In - \Out tests described in \secref{sec: in-out}. The elements are labeled as \In , \Out or \Intercepted based on the position of the \textit{.stl} geometry. To avoid repeated tests these labels are stored as element markers. Elements have different routines, and computational costs, based on the element marker. For \In elements, no integration is performed and Dirichlet conditions are set on all CG nodes corresponding to \In elements. For \Out elements, integration is performed by looping through all quadrature points. For \Intercepted elements, we perform adaptive quadrature \secref{sec: VGP} and perform the integration only over \Out Gauss points. The matrix assembly is an important part here and is optimized for achieving significant speedup and is discussed in \secref{sec: MatrixAssembly}.

\subsection{Octree based Mesh Generation} \label{sec: octree}
Octrees are widely used in computational sciences~\cite{SundarSampathBiros08,BursteddeWilcoxGhattas11,Fernando2018_GR,Fernando:2017, khanwale2020simulating,xu2021octree}, due to its conceptual simplicity and ability to scale across large number of processors. 
Proceeding in a top-down fashion, an octant is refined based on a user-specified criteria. The refinement criteria is specified by a user-specified function that takes the coordinates of the octant, and returns \texttt{true} or \texttt{false}. Since the refinement happens in an element-local fashion, this step is embarrassingly parallel. In distributed memory, the initial top-down tree construction, also enables an efficient partitioning of the domain across an arbitrary number of processes. All processes start at the root node (i.e., the cubic bounding box for the entire domain). In order to avoid communication during the refinement stage, we opt to perform redundant computations on all processes. Starting from the root node, all processes refine (similar to the sequential code) until at least $\mathcal{O}(p)$ octants requiring further refinement are produced. Then using a weighted space-filling-curve (SFC) based partitioning, we partition the octants across all processes. Note that we do not communicate the octants as every process has a copy of the octants, and all that needs to be done at each process is to retain a subset of the current octants corresponding to its sub-domain. This allows us to have excellent scalability, as all processes perform (roughly) the same amount of work without requiring any communication.

Upon octree generation we enforce the 2:1 balanced constraint which ensures that neighboring octants differ by at most one level. 
Such a 2:1 balanced constraint simplifies mesh generation, and eliminates ill-conditioning introduced due to extreme scaling of neighbor elements. Following balancing, meshing is performed which generates the required data-structures to perform FE computations, intergrid transfers, spatial queries and membership tests. Additional details on our octree-based FEM framework can be found in \citet{Fernando:2017,Dendro} and \citet{Dendro5}.

\begin{algorithm}[t]
  \caption{\textsc{Immersogeometric\_algorithm}:\\ Brief Overview}\label{alg:overview}
  \footnotesize
  \begin{algorithmic}[1]
      \Require The octree and $.stl$ file
    
    \For{$elements \in leaf$}
    \State compute element Markers 
    \Comment{Algorithm~\ref{alg:marker}}
    \State compute background triangles ($\mathcal{B_T}$)
    \EndFor
    \For{time $<$ $T_\mathrm{final}$}
    \Comment{Loop over time}
    \For{$elements \in leaf$}
    \If{\textsc{\In}} 
    \State continue
    \Comment{Skip \In elements}
    \EndIf
    \If{\textsc{\Out}} 
        \State Loop over Gauss Points
        \State perform Matrix and Vector Assembly
    \EndIf
    \If{\textsc{\Intercepted}} 
        \State Fill the elements with Gauss points
        \State Check if the Gauss point $\mathcal{P}$ is \In or \Out 
        \Comment{Algorithm~\ref{alg:in-out}}
        \If{($\mathcal{P}$ is \Out)}
            \State perform Matrix and Vector Assembly
        \EndIf
    \EndIf
    \EndFor
    
    \For{$element \in leaf$}
    \If{\textsc{\Intercepted}} 
    \For{$t \in \mathcal{B_T}$}
    \State Assemble weak BC contribution to $element$
    \EndFor
    \EndIf
    \EndFor
    \State Solve system of equations
    \EndFor
  \end{algorithmic}
\end{algorithm}

\subsection{Matrix Assembly on Distributed Octrees}\label{sec: MatrixAssembly}

\begin{figure}[b!]
    \centering
    \begin{tikzpicture}
    \begin{scope}{shift={(0,0)}}
        \draw[] node at (4,3.5) {\textit{memory layout for distributed nodes/octants}};
        \draw[line width=1.5mm,-,armygreen] (2,3) -- (6,3)  node at (4,2.5) {local ($\tau_k$)};
        \draw[line width=1.5mm,-,applegreen] (2,2) -- (3,2)  node at (2,1.5) {\textcolor{black}{\small \it dependent} };
        \draw[line width=1.5mm,-,green] (3.1,2) -- (5,2)  node at (4.1,1.5) {\textcolor{black}{\small \it independent}};
        \draw[line width=1.5mm,-,applegreen] (5.1,2) -- (6,2)  node at (6.0,1.5) {\textcolor{black}{\small \it dependent}};
        \draw[line width=1.5mm,-,aurometalsaurus] (0,3) -- (2,3)  node at (1,2.5) {pre-ghost};    
        \draw[line width=1.5mm,-,aurometalsaurus] (6,3) -- (8,3)  node at (7,2.5) {post-ghost};   
    \end{scope}
    \end{tikzpicture}
    \caption{Figure depicting the memory layout for the distributed nodes/octants. The distributed octree $T$ is partitioned across $p$ processors, where each partitioned $\tau_k$ tree (local octants) has pre (from $t\tau_m m < k$) and post (from $t\tau_m m > k$) ghost octant/nodal information. The local partition $\tau_k$ is further logically decomposed in to independent $\tau_I$, and dependent $\tau_D$ disjoint octant sets such that $\tau_k = \tau_I \cup \tau_D$. This notion of independent and dependent is used in overlapping computation with communication in FEM matrix assembly.}
    \label{fig:mem_layout}
\end{figure}
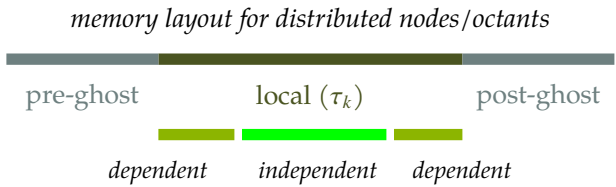{}

In FEM, the differential operator $\mathcal{L}$, after weakening and discretization becomes a matrix (stiffness matrix) \added{$K_{ij} =  (\phi _i, \mathcal{L}\phi_j )$} where $(u,v) = \int_{\Omega} u v d\Omega$, and $\phi_i$ is the basis function defined at $i^{th}$ node. In a distributed octree, each partition $\tau_k$ of $\Omega$ will compute the $K_{ij}$ restricted to $\tau_k$ denoted by $K_{ij| \tau_k}$. The octants in $\tau_k$ are further decomposed into two disjoint octant sets -- independent $\tau_I$, and dependent $\tau_D$ such that, $\tau_k = \tau_I \cup \tau_D$. An octant $e\in \tau_I$ ensures that all the nodal information related to $e$ is local while $\tau_D$ is $\tau_k\setminus \tau_I$. The notion of independent and dependent octants is used to overlap computation with communication during the matrix assembly (see \figref{fig:mem_layout}).

\begin{figure}[t!]
\resizebox{\columnwidth}{!}
{\centering
	\begin{tikzpicture}[scale=0.3,>=stealth,bluearr/.style={->,blue,shorten >= 3pt}]

		\begin{scope}[shift={(0,0)}]
			\drawcubeII(0,0,0,8,black,blue!50,0.3);
			\drawcubeII(10,0,0,4,black,green!50,0.3);
			\foreach \y in {0,2,...,4}
			{
				\foreach \z in {0,2,...,4}{
					\draw[fill=red!60] (10,\y,\z) circle (0.12);
					\draw[fill=red!80] (8,2*\y,2*\z) circle (0.15);
					\draw[bluearr] (10,\y,\z) to [bend left=45] node [midway,above,anchor=east,align=center] {} (8,2*\y,2*\z);
				}
			}	
		
			\end{scope}
				 
			\begin{scope}[shift={(18,0)}]
		
			\drawcubeII(0,0,0,8,black,blue!50,0.3);
			\drawcubeII(10,0,-4,4,black,green!50,0.3);
			\foreach \z in {0,2,...,4}{
					\draw[fill=red!60] (10,\z,0) circle (0.12);
					\draw[fill=red!80] (8,2*\z,0) circle (0.15);
					\draw[bluearr] (10,\z,0) to [bend left=45] node [midway,above,anchor=east,align=center] {} (8,2*\z,0);
			}
		
			\end{scope}
	\end{tikzpicture}
}
\caption{\label{fig:hangingElementsNodes} \small An example of a hanging face and a hanging edge where in both cases octant (\textcolor{green!50}{$\blacksquare$}) has a hanging face (left figure) and a hanging edge (right figure) with octant (\textcolor{blue!50}{$\blacksquare$}). Nodes on the hanging face/edge are mapped to the larger octant and the hanging nodal values are obtained via interpolation. 
 Note that for illustrative purposes, the two octants are drawn separately, but are contiguous.}
\vspace{0.2in}
\end{figure}
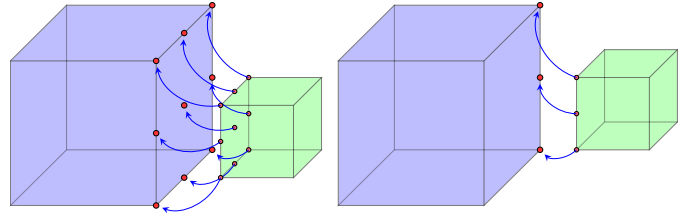

{\bf Hanging nodes:} Adaptivity in octree meshes causes non-conformity. In our framework, additional degrees of freedom on shared faces between elements of unequal sizes---the so called hanging nodes---do not represent independent degrees of freedom, and are not stored and instead are represented as a linear combination of the basis functions corresponding to the larger face (see \figref{fig:hangingElementsNodes}). 
\added{We implement this using a correction operator~\citep{Fernando2018_GR,kus2014arbitrary,sundar2008bottom} during the overall matrix assembly
(see Algorithm~\ref{alg:matassembly} and Algorithm~\ref{alg:hanging_corrections}).}

\begin{algorithm}[t!]
  \caption{\textsc{Matrix\_Assembly}}\label{alg:matassembly}
  \footnotesize
  \begin{algorithmic}[1]
    \Require Octree $T$, distributed across $p$ processors, ($\tau_k$ local partition)
    \Ensure global assembled matrix $K_{ij}$
    \State $K\leftarrow 0$ matrix
    \For{$e \in \tau_D$} \Comment{local dependent elements}
        \State $K_e \leftarrow compute\_ele\_matrix()$
        \Comment{Algorithm~\ref{alg:ele_mat_computation}}
        \State $M_e \leftarrow compute\_hanging\_corrections()$
        \Comment{Algorithm~\ref{alg:hanging_corrections}}
        \State $k_e \leftarrow M_e^T K_e M_e$
    \EndFor
    \State $start\_assembly\_comm()$ \Comment{start communication}
    \For{$e \in \tau_I$} \Comment{local independent elements}
    \State $k_e \leftarrow compute\_ele\_matrix()$ \Comment{Algorithm~\ref{alg:ele_mat_computation}}
        \State $M_e \leftarrow compute\_hanging\_corrections(\tau_k,e)$ \Comment{Algorithm~\ref{alg:hanging_corrections}}
        \State $K_e \leftarrow M_e^T k_e M_e$
        \State $K\leftarrow K + \oTn(K_e)$
    \EndFor
    \State $end\_assembly\_comm()$ \Comment{wait till communication ends}
    \For{$e \in \tau_D$}
        \State $K\leftarrow K + \oTn(K_e)$
    \EndFor
  \end{algorithmic}
\end{algorithm}

\begin{algorithm}[t!]
  \caption{\textsc{compute\_hanging\_corrections}}\label{alg:hanging_corrections}
  \footnotesize
  \begin{algorithmic}[1]
    \Require Octree $T$, distributed across $p$ processors, ($\tau_k$ local partition), local element $e\in \tau_k$
    \Ensure Hanging node correction matrix $M_e$
    \State $M_e \leftarrow I$ \Comment{$I$ is the identity matrix}
    \State $I2d_{ij} \leftarrow I1d_i \bigotimes I1d_j$ \Comment{$I1d$ is 1d interpolation matrices}
    \For{$f \in Faces(e)$ }
        \If{$is\_hanging(\tau_k,e,f)$}
            \State $ M_e[f] \leftarrow I2d_{ij} $ 
        \EndIf
    \EndFor
    \For{$edge \in Edges(e)$ }
        \If{$is\_hanging(\tau_k,e,edge)$}
            \State $ M_e[edge] \leftarrow I1d_{i} $ 
        \EndIf
    \EndFor
    \Return $M_e$
  \end{algorithmic}
\end{algorithm}
{\bf Elemental matrix assembly:}
\added{We accelerate the assembly process\footnote{Although we illustrate some promising results with a matrix-free approach \cite{ishii2019solving}, further optimization is needed---specifically in terms of tailored pre-conditioners---to make the matrix-free approach competitive.} 
by viewing each weakened differential FEM operator (FEM kernel) of the form $(\textsc{Mat\_op1},\textsc{Gp}_{i} \textsc{Mat\_op2})$ as an outer product and making use of optimized matrix-matrix multiplication libraries. Each FEM kernel comprises of two operators (\textsc{Mat\_op1} and \textsc{Mat\_op2}) obtained as a result of weakening, and \textsc{Gp}$_{i}$ denotes the value of interpolated variables at the Gauss points. Algorithm~\ref{alg:ele_mat_computation} describes the procedure along with the associated complexity of each step.  In this work, we represent \textsc{Mat\_op} as a matrix of appropriate size. For better performance, the matrix corresponding to \textsc{Mat\_op} can be precomputed and cached (~\secref{sec: MatPrecompute}). Step 1-4 denotes the computation of second term in the kernel (\textsc{Gp}$_{i}$ \textsc{Mat\_op2}). $W_i$ denotes the weight of the quadrature points. Step 5 donates the action of \textsc{Mat\_op1}. The scalar factor \textsc{Scale} is the resultant of mapping from the reference frame to the physical frame and is operator dependent \footnote{This makes use of the fact that the octree elements, in general, are cubes with equal aspect ratio $(\Delta x = \Delta y = \Delta z)$}. For example, the scale factor for the stiffness matrix is $\Delta x/2$ and for the mass matrix is $\Delta x^3/8$ for a reference element $\in [-1,1]^3$ in 3D. Clearly, this last step is computationally the most expensive. For a detailed analysis, interested readers are referred to \citet{deville2002high}. 

Our approach is different from the general approach in open source libraries (such as deal.ii \cite{dealII90}, for instance) that loop over individual Gauss points in each direction. Although the overall theoretical complexity remains the same ($\mathcal{O}(\mathrm{bf + 1})^{3d}$, where bf is the basis function order and $d$ is the number of spatial dimensions) for both implementations, implementing the operator as a matrix-matrix multiplication allows a natural way to leverage the power of \textsc{GEMM} kernel~\citep{anderson1999lapack} that is available as vendor optimized libraries (such as \textsc{Intel MKL}). These libraries are fine tuned to exploit assembly-code-level parallelism and are extremely fast on modern cache-based and pipelined architectures. As evident from the complexity analysis, optimized matrix assembly becomes especially important when using higher order basis functions. This is a careful middle ground that ensures portability across various computing platforms, while not as efficient as explicit vectorization.

It has been shown that the theoretical lower bound on the computational complexity for elemental matrix assembly is $\mathcal{O}(\mathrm{bf + 1})^{2d}$, since the elemental matrix has many non-zero entries~\citep{deville2002high}. There has been substantial effort to achieve this lower bound~\cite{bressan2019sum,melenk2001fully,deville2002high}. However, most of the proposed analysis and implementation have been limited to stiffness and mass matrix. Deploying these algorithms for complicated non-linear kernels arising as a result of weakening the Navier--Stokes equation seems non-trivial and needs further analysis. This is left as future work.}

\begin{algorithm}[t!]
  \caption{\textsc{compute\_ele\_matrix}}\label{alg:ele_mat_computation}
  \footnotesize
  \begin{algorithmic}[1]
    \Require FEM operators (\textsc{Mat\_op1} and  \textsc{Mat\_op2}), Quadrature Values (\textsc{Gp}), \textsc{Scale}, Size of matrix(n $\times$ n): n =  $(\mathrm{bf}+1)^d$
    \Ensure The computed elemental matrix \textsc{Emat}
    \State \textsc{Gp} $\leftarrow$ \textsc{Gp}$_i \times W_i$
    \Comment{$\mathcal{O}(\mathrm{n})$}
    \For{i $\leftarrow$ 1 to n}
    \For{j $\leftarrow$ 1 to n}
    \State \textsc{M}$_{ij}$ $\leftarrow$  \textsc{Mat\_op2}$_{ij} \times$ \textsc{Gp}$_j$ \Comment{$\mathcal{O}(\mathrm{n})^{\mathrm{2}}$}
    \EndFor
    \EndFor
    \State \textsc{GEMM:}\textsc{Emat} $\leftarrow$ \textsc{Scale}(\textsc{Mat\_op1}$^T$ \textsc{M}) \Comment{$\mathcal{O}(\mathrm{n})^{\mathrm{3}}$}
    
    \Return \textsc{Emat}
  \end{algorithmic}
\end{algorithm}


\subsection{\In - \Out Test}\label{sec: in-out}
\begin{algorithm}[t]
  \caption{\textsc{element\_markers}~: \In / \Out or \Intercepted}\label{alg:marker}
  \footnotesize
  \begin{algorithmic}[1]
      \Require The element $\mathcal{E}$ and  list of background triangles $\mathcal{B_T}$ associated with it.
    \Ensure The marker for $\mathcal{E}$
    \Comment{\In , \Out or \Intercepted}
    \State $counts \leftarrow 0$
    \For{$n \in n_{nodes}$}
    \State compute $\mathcal{P}$ 
    \Comment Global position of the node
    \If{\textsc{isIn($\mathcal{P}$,$\mathcal{B_T}$)}} \Comment{Algorithm \ref{alg:in-out}}
        \State increment count
    \EndIf
    \EndFor
    \If{count == $n_{nodes}$} \Comment{All nodes are \In}
     \State \Return \In 
     \EndIf
    \If{count == 0}         \Comment{All nodes are \Out}
     \State \Return \Out
     \EndIf
     \Return \Intercepted
  \end{algorithmic}
\end{algorithm}

\begin{algorithm}[t]
  \caption{\textsc{is\In}~: Normal based \In - \Out test}\label{alg:in-out}
  \footnotesize
  \begin{algorithmic}[1]
      \Require The global position $\mathcal{P}$ and  list of background triangles $\mathcal{B_T}$ associated with each octant.
    \Ensure Point $\mathcal{P}$ is \In or \Out of the given geometry
    \State $count \leftarrow 0$
    \For{$t \in \mathcal{B_T}$}
    \State Compute $\Vector{d}$ = $\mathcal{P}$ - t \Comment{distance vector from $\mathcal{P}$ to t}
    \If{$\Vector{d} \cdot \Vector{n}_t \leq 0$} 
    \State increment count
    \EndIf
    \EndFor
    \If{$count == |\mathcal{B_T}|$}
    \State \Return \Out
    \Comment{All background triangles indicate that point is \Out}
    \EndIf
    \If{$count == 0$}
    \State \Return \In
    \Comment{All background triangles indicate that point is \In}
    \EndIf
    \State perform \textsc{Ray - Tracing }
    \Comment{Normal based test fails}
  \end{algorithmic}
\end{algorithm}

In this work, we propose a new normal based identification of \In - \Out test, as opposed to ray-tracing that is conventionally use in immersed boundary simulations~\cite{khalighi2009validation,borazjani2008curvilinear,de2006recent,iaccarino2003immersed}. Algorithm~\ref{alg:in-out} describes the procedure to identify whether a given quadrature point $\mathcal{P}$ in an \Intercepted element is \In or \Out. Given a point $\mathcal{P}$ in an element $\mathcal{E}$, we identify the background triangles (that is stored in the data structure during partition of triangles) associated with the element. The dot product of the position vector (i.e.~vector from $\mathcal{P}$ to the triangle centroid) with the normal of all the background triangles $\mathcal{B_T}$ is evaluated. If the result of the dot product with all the triangles is greater than 0, then the point is outside the geometry and vice versa. The normal based test comes at no additional cost in terms of memory requirements and involves series of simple dot product evaluation between the point and background triangles. This makes it very cost efficient as compared to ray-tracing. In case of conflict, as in case of sharp corners (see \figref{fig: InOut}), we revert back to ray-tracing. 

\pagebreak

\input{Plots/InOut}

\subsection{Adaptive Quadrature} \label{sec: VGP}
\added{There has been recent advances in the development of accurate and efficient quadrature rules to account for the cut cells that arise during IMGA computations \citep{thiagarajan2014adaptively,thiagarajan2016adaptively,duczek2015efficient,schillinger2013review,barendrecht2018efficient,stavrev2016geometrically,schillinger2015finite,kudela2015efficient,xu2016tetrahedral,divi2020error}. In this work, to accurately account for the effect of fractional content of fluids in the \Intercepted cells, we use adaptive quadrature by recursive cell subdivision.} \Intercepted elements are further subdivided to include additional Gauss points. This ensures the accurate integration of the fluid domain for intersected elements. The splitting of the element is done \textit{only} to add Gauss points and does not introduce any new degrees of freedom. The \Intercepted elements are recursively sub divided to fill in Gauss points till the smallest cells are completely \In (or \Out); or a max split level is reached. \figref{fig:GP} demonstrates the case with maximum split level set to 2.


\input{Plots/GaussQaudrature}

\subsection{Dynamic Load Balancing}  

We use a Space Filling Curve (SFC) -- specifically the Hilbert curve -- to partition our octree mesh and the associated data across all processes. 
In the case of IMGA, different elements can have different computational loads depending on the \In, \Out, and \Intercepted status. \Intercepted elements perform the following additional computations compared to \In and \Out elements: a) identify whether a particular Gauss point is in or out , b) loop over additional Gauss points (due to recursive quadrature), and c) perform surface integral over the surface element for accumulation of the contribution from the weak boundary condition. We ensure load-balance across all processes by using a weighted SFC partition, wherein \Intercepted elements are assigned a higher weight. The weight associated with \Intercepted elements is proportional to the ratio of computational effort of an \Intercepted elements to an \Out element. 

\subsubsection{Estimation of Weight for the \Intercepted Elements}
Let $T_v$ be the computational cost per volumetric quadrature point and $T_s$ be the cost per surface quadrature point. The relative weight for each element, $e$, can be estimated as:

\begin{equation}
        W(e) = \frac{n_\mathrm{Vgp}(e)}{n_\mathrm{Tgp}} + \frac{T_s}{T_v}\frac{n_\mathrm{Sgp}(e)}{n_\mathrm{Tgp}} 
    \label{eq: Wpart}
\end{equation}

\noindent where: $n_\mathrm{Vgp}$ is the number of Gauss point that are outside the geometry, $n_\mathrm{Sgp}$ is the total number of surface Gauss points belonging to the \Intercepted elements and $n_\mathrm{Tgp}$ represents the total number of Gauss point in the volumetric elements which scales as $(\mathrm{bf}+1)^{\mathrm{nsd}}$, where bf is the order of basis function. For a completely \Out element and a completely \In element, \eqnref{eq: Wpart} reduces to 1 and 0, respectively.

%% file: Plots/InOut.tex
\begin{figure}[t!]
        \centering
        \begin{subfigure}[b]{0.2\textwidth}
                \centering
                \begin{tikzpicture}[scale=0.8]
                 \draw[step=4.0cm, thick] (0,0) grid (4.0,4.0);
                \draw[red, thick] (-0.1,0.1) -- (2.5,1.8);
                \draw[red, thick] (2.5,1.8) -- (4.1,4.1);
                \draw [fill=black](0.5,2.7) circle (0.1cm);
                \node at (0.5,3.1) {$\mathcal{P}$};
                
                \draw [->] (1.2,0.95) -- (1.6925,0.1967);
                \draw [->,dashed] (0.5,2.7) -- (1.2,0.95);
                \node at (1.8,3.2) {$\Vector{d}_1$};
                \node at (3.5,2.5) {$\Vector{n}_1$};
                 
                \draw [->] (3.3,2.95) -- ( 4,2.4);
                \draw [->,dashed] (0.5,2.7) -- (3.3,2.95);
                \node at (1.3,1.5) {$\Vector{d}_2$};
                \node at (1.8,0.6) {$\Vector{n}_2$};
                \end{tikzpicture}
                \caption{No sharp corners}
                \label{fig:inOutSuc}
        \end{subfigure}%
        \quad  \quad
        \begin{subfigure}[b]{0.2\textwidth}
                \centering
                \begin{tikzpicture}[scale=0.8]
                 \draw[step=4.0cm, thick] (0,0) grid (4.0,4.0);
                \draw[red, thick] (2,1) -- (4.2,1);
                \draw[red, thick] (2,1) -- (4.1,4.1);
                \draw [fill=black](0.5,2.7) circle (0.1cm);
                \node at (0.5,2.1) {$\mathcal{P}$};
                
                \draw [->] (3.0161,2.5) -- (2.2710,3.0048);
                \draw [->,dashed] (0.5,2.7) -- (3.0161,2.5);
                \node at (1.8,2.8) {$\Vector{d}_1$};
                \node at (2.7,3.1) {$\Vector{n}_1$};
                 
                \draw [->] (3,1) -- (3,0.2);
                \draw [->,dashed] (0.5,2.7) -- (3,1);
                \node at (1.5,1.4) {$\Vector{d}_2$};
                \node at (2.5,0.6) {$\Vector{n}_2$};
                \end{tikzpicture}
              \caption{Sharp corners}
              \label{fig:inOutFail}
        \end{subfigure}
        \vspace{3 mm}
        \caption{Figure illustrating the normal based \In - \Out test. \figref{fig:inOutSuc} shows the case where there is no sharp corners. In this case, the dot product of $\Vector{d}_1 \cdot \Vector{n}_1$ and $\Vector{d}_2 \cdot \Vector{n}_2$ both are  greater than  0 stating that point is \In. But in case of sharp corners, as shown in  \figref{fig:inOutFail} the background triangles can give conflicting decisions. $\Vector{d}_1\cdot \Vector{n}_1 < 0$ and $\Vector{d}_2 \cdot \Vector{n}_2 > 0$ In these cases, we resort to the ray-tracing algorithm. }
        \label{fig: InOut}
       
\end{figure}
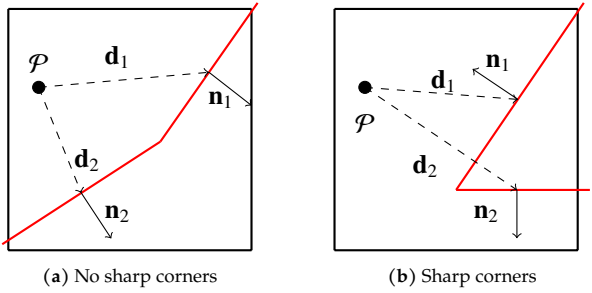

%% file: Plots/GaussQaudrature.tex
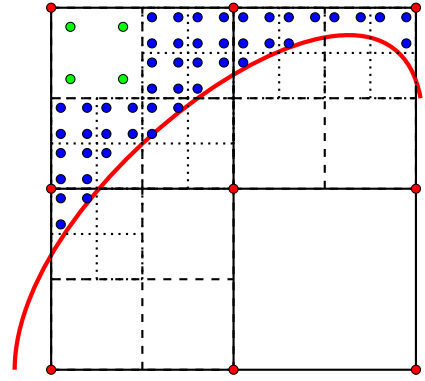
\begin{figure}
\centering
\vspace{-0.4in}
    \begin{tikzpicture}[scale=0.6]
     \draw[step=4.0cm, thick] (0,0) grid (8.0,8.0);
     \draw [ultra thick,red] (-0.8,0) to[out=90,in=100] (8.1,6);
     \draw[step=2.0cm, thick,dashed] (0,0) grid (4.0,4.0);
     \draw[step=2.0cm, thick,dashed] (0,4) grid (4.0,8.0);
     \draw[step=2.0cm, thick,dashed] (4,4) grid (8.0,8.0);
     \draw[step=1.0cm, thick,dotted] (0,2) grid (2.0,4.0);
     \draw[step=1.0cm, thick,dotted] (0,4) grid (4.0,6.0);
     \draw[step=1.0cm, thick,dotted] (2,6) grid (8.0,8.0);
    \foreach \x in {0.0,4.0,8.0}{
    \foreach \y in {0.0,4.0,8.0}{
        \draw[fill=red](\x,\y) circle (0.1cm);
        }
    }
    \draw[fill=green](2*0.2113248654,6.0 + 2*0.2113248654) circle (0.1cm);
    \draw[fill=green](2*0.2113248654,6.0 + 2*0.78867513459) circle (0.1cm);
    \draw[fill=green](2*0.78867513459,6.0 + 2*0.2113248654) circle (0.1cm);
    \draw[fill=green](2*0.78867513459,6.0 + 2*0.78867513459) circle (0.1cm);

    \draw[fill=blue](0.2113248654,4.0 + 0.2113248654) circle (0.1cm);
    \draw[fill=blue](0.2113248654,4.0 + 0.78867513459) circle (0.1cm);
    \draw[fill=blue](0.78867513459,4.0 +0.2113248654) circle (0.1cm);
    \draw[fill=blue](0.78867513459,4.0 +0.78867513459) circle (0.1cm);

     \draw[fill=blue](0.2113248654,5.0 + 0.2113248654) circle (0.1cm);
    \draw[fill=blue](0.2113248654, 5.0 + 0.78867513459) circle (0.1cm);
    \draw[fill=blue](0.78867513459,5.0 +0.2113248654) circle (0.1cm);
    \draw[fill=blue](0.78867513459,5.0 +0.78867513459) circle (0.1cm);
    
    \draw[fill=blue](1.0+0.2113248654, 5.0 + 0.2113248654) circle (0.1cm);
    \draw[fill=blue](1.0+0.2113248654, 5.0 + 0.78867513459) circle (0.1cm);
    \draw[fill=blue](1.0+0.78867513459,5.0 +0.2113248654) circle (0.1cm);
    \draw[fill=blue](1.0+0.78867513459,5.0 +0.78867513459) circle (0.1cm);
    
    \draw[fill=blue](1.0+0.2113248654, 4.0 + 0.78867513459) circle (0.1cm);
    
    \draw[fill=blue](2.0+0.2113248654, 6.0 + 0.2113248654) circle (0.1cm);
    \draw[fill=blue](2.0+0.2113248654, 6.0 + 0.78867513459) circle (0.1cm);
    \draw[fill=blue](2.0+0.78867513459,6.0 +0.2113248654) circle (0.1cm);
    \draw[fill=blue](2.0+0.78867513459,6.0 +0.78867513459) circle (0.1cm);
    
    \draw[fill=blue](2.0+0.2113248654, 7.0 + 0.2113248654) circle (0.1cm);
    \draw[fill=blue](2.0+0.2113248654, 7.0 + 0.78867513459) circle (0.1cm);
    \draw[fill=blue](2.0+0.78867513459,7.0 +0.2113248654) circle (0.1cm);
    \draw[fill=blue](2.0+0.78867513459,7.0 +0.78867513459) circle (0.1cm);
    
    \draw[fill=blue](3.0+0.2113248654, 7.0 + 0.2113248654) circle (0.1cm);
    \draw[fill=blue](3.0+0.2113248654, 7.0 + 0.78867513459) circle (0.1cm);
    \draw[fill=blue](3.0+0.78867513459,7.0 +0.2113248654) circle (0.1cm);
    \draw[fill=blue](3.0+0.78867513459,7.0 +0.78867513459) circle (0.1cm);

    \draw[fill=blue](3.0+0.2113248654, 6.0 + 0.2113248654) circle (0.1cm);
    \draw[fill=blue](3.0+0.2113248654, 6.0 + 0.78867513459) circle (0.1cm);
    \draw[fill=blue](3.0+0.78867513459,6.0 +0.78867513459) circle (0.1cm);
    
    \draw[fill=blue](4.0+0.2113248654, 6.0 + 0.78867513459) circle (0.1cm);
    
    \draw[fill=blue](4.0+0.2113248654, 7.0 + 0.2113248654) circle (0.1cm);
    \draw[fill=blue](4.0+0.2113248654, 7.0 + 0.78867513459) circle (0.1cm);
    \draw[fill=blue](4.0+0.78867513459,7.0 +0.2113248654) circle (0.1cm);
    \draw[fill=blue](4.0+0.78867513459,7.0 +0.78867513459) circle (0.1cm);
    
    \draw[fill=blue](5.0+0.2113248654, 7.0 + 0.2113248654) circle (0.1cm);
    \draw[fill=blue](5.0+0.2113248654, 7.0 + 0.78867513459) circle (0.1cm);
    \draw[fill=blue](5.0+0.78867513459,7.0 +0.78867513459) circle (0.1cm);
    
    \draw[fill=blue](6.0+0.2113248654, 7.0 + 0.78867513459) circle (0.1cm);
    \draw[fill=blue](6.0+0.78867513459,7.0 +0.78867513459) circle (0.1cm);

    \draw[fill=blue](7.0+0.2113248654, 7.0 + 0.78867513459) circle (0.1cm);
    \draw[fill=blue](7.0+0.78867513459,7.0 +0.2113248654) circle (0.1cm);
    \draw[fill=blue](7.0+0.78867513459,7.0 +0.78867513459) circle (0.1cm);
    
    \draw[fill=blue](2.0+0.2113248654, 5.0 + 0.2113248654) circle (0.1cm);
    \draw[fill=blue](2.0+0.2113248654, 5.0 + 0.78867513459) circle (0.1cm);
    \draw[fill=blue](2.0+0.78867513459,5.0 +0.78867513459) circle (0.1cm);

    \draw[fill=blue](0.2113248654, 3.0 + 0.2113248654) circle (0.1cm);
    \draw[fill=blue](0.2113248654, 3.0 + 0.78867513459) circle (0.1cm);
    \draw[fill=blue](0.78867513459,3.0 +0.78867513459) circle (0.1cm);
    \end{tikzpicture}
    \caption{Sketch demonstrating adaptive quadrature. Boundary is represented by the red line. Octant shared nodes are represented by (\tikzcircle[fill=red]{0.1 cm}). (\tikzcircle[fill=green]{0.1 cm}) represents the new Gauss points generated as the result of $1^{st}$ level of splitting and  (\tikzcircle[fill=blue]{0.1 cm}) represents the Gauss point of $2^{nd}$ level of splitting. It is to be noted that the splitting happens only for the \Intercepted  elements. Once all the Gauss points in the splitted elements are either in or out the elements is not further sub-divided.}
    \label{fig:GP}
\end{figure}

%% file: Result.tex
In this section, we present the benchmark results for our solver in terms of accuracy and speed by considering a suite of  appropriate test cases. We begin by validating our framework using a canonical problem of flow past a sphere across a wide range of Reynolds numbers, encompassing the laminar-transition-turbulent regimes (\secref{sec:sphere}). Next, we show the impact of tensorization on the matrix assembly and solve time using a benchmark lid driven cavity problem (\secref{results:mat-assembly}). We next quantify the advantage of our proposed \In - \Out test for significantly complex standard geometries in \secref{Sec:InOutTest}. In order to articulate the advantage of using adaptive quadrature and weighted partitioning on IMGA simulations, we revisit flow past a sphere in \secref{Sec:WeightedPartioning}. Finally, we present the scaling results of our solver for different mesh resolutions (\secref{Sec:Scaling}).



\subsection{Validation: Flow Past a Sphere}\label{sec:sphere}
\begin{figure}[b!]
    \centering
    \includegraphics[width=0.95\linewidth,clip,trim={2.0in 0.85in 2.0in 1.0in}]{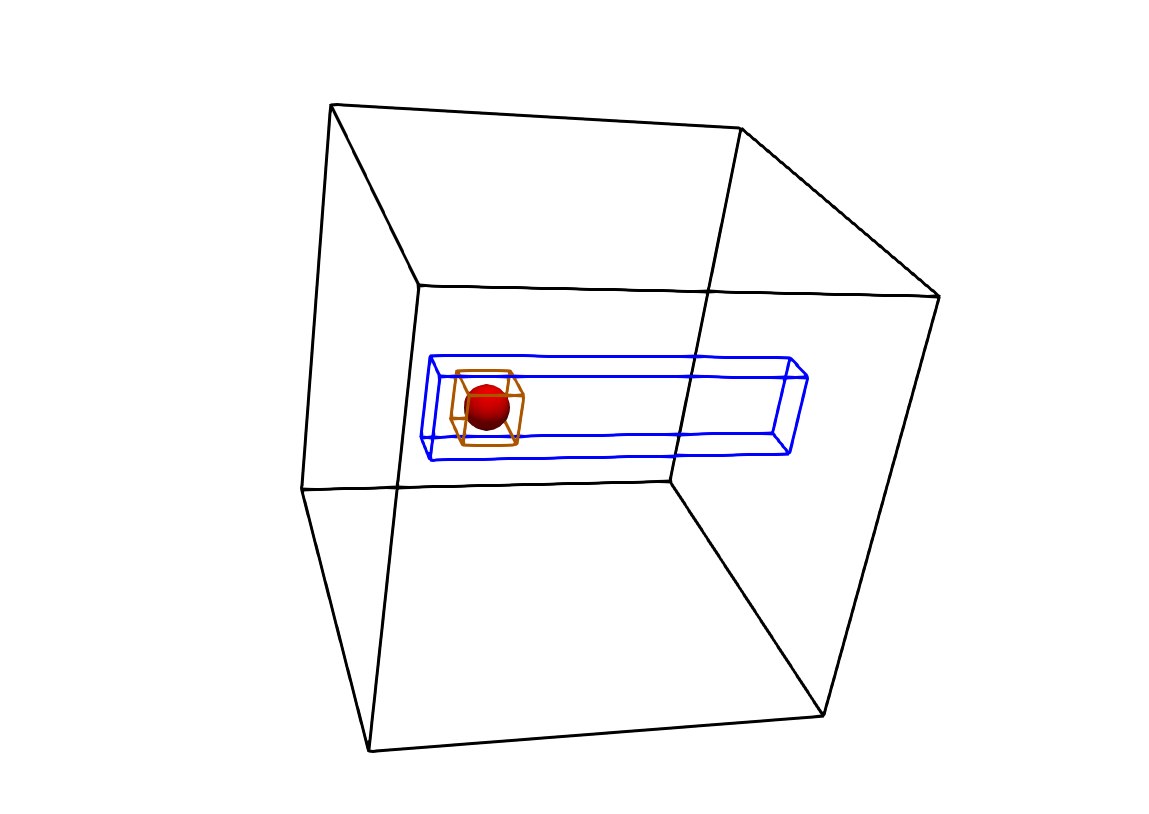}
    \caption{{Computational domain for sphere simulations. The computation domain consists of a cubic domain with dimension $10d \times 10d \times 10d$. Overall there are three different levels of refinement: $R_{bkg}$ for the background mesh, $R_{wake}$ within the blue box to capture the wake, and $R_{bdy}$ near the body within the brown box. The refinement level inside the body is equal to $R_{bkg}$. The sphere of diameter $d=1$ is positioned at $(3d,\,5d,\,5d)$.} A refinement level of $i$ would corresponding to the resolution of $\frac{10d}{2^i}$ in each direction.}
    \label{fig:scaling-sphere}
\end{figure}
\input{Plots/DragCompar}
\input{Plots/SpherePlot}
As the first step, we validate our numerical setup by computing the non-dimensional drag coefficient $C_d$ on sphere in a range of Reynolds number from 1 to 1,000,000. \figref{fig:scaling-sphere} shows the schematic of the computational domain used to perform the simulation for flow past the sphere. The computational domain consists of length $10d$ with the center of the sphere of diameter $d = 1$ placed at distance $3d$ from inlet (at $(3d, \,5d,\, 5d)$). All the walls of the cubic domain, except the outlet have constant non-dimensional freestream velocity of $(1, \, 0, \, 0)$ and zero pressure gradient. At the outlet, the pressure is set to 0 and zero gradient velocity boundary condition is applied at the wall. The no-slip boundary condition (zero Dirichlet) for velocity is weakly imposed on the surface of sphere. The octree and surface mesh resolution was varied by successively refining the mesh depending on the Reynolds number according to the Taylor micro length scale ($R_{bdy}$ was varied from 5 to 11 with increasing $Re$). At the finest resolution for the case of $Re = 1$ million, we have the equivalent of around $500$ elements across the diameter of sphere. The spatial domain was discretized using linear basis function. The Crank Nicholson scheme was used for time discretization. 
\begin{remark}
The finest resolution for $Re = 1$ million was chosen based on the work by \citet{geier2017parametrization}. This resolution has been shown to be sufficient to capture the drag crisis. Further refinement was not considered in this work due to the high computational cost. 
\end{remark}

We evaluated the non-dimensional drag coefficient across this range of Reynolds numbers, which is plotted in \figref{fig: Drag}. We can see that the IMGA results are in excellent agreement with experimental data. Note that we are able to accurately capture the drag crisis phenomena, where a sudden drop in drag coefficient is observed. The prediction of drag crisis is of significant interest to the CFD community, with important engineering implications in aerodynamics and vehicular dynamics. We emphasize that no special numerical/modelling treatment is needed in this framework to capture the drag crisis; this is in contrast to existing work where typically a wall treatment 
\citep{hoffman2006simulation,constantinescu2004numerical} is required to predict the drag crisis.


\figref{fig: wake} show the flow structures in the wake behind the sphere for increasing $Re$.
At $Re = 100$, an attached stable ring eddy is formed behind the surface of the sphere (\figref{fig: Re100}). At $Re$ in the range of 1000, the ring eddy becomes unstable and starts to shed. 
Further increase in $Re$ results in a fully turbulent wake (\figref{fig: Re160K} and \figref{fig: Re1M}). The drag crisis is evident by noticing difference in the wake structure between \figref{fig: Re160K} and \figref{fig: Re1M}. 
In the subcritical regime, the wake behind the sphere tends to diverge as it moves away from the sphere. Above the critical $Re$, the wake tends to converge and become narrower as it moves behind the sphere. Subsequently, we can observe a high pressure zone being developed behind the sphere. The main reason for the drag crisis can be attributed to the development of this high pressure zone behind the sphere that tends to push the sphere in the inflow direction. These results are in agreement with the simulations by \citet{geier2017parametrization} using Lattice Boltzmann method to capture the drag crisis.


\begin{remark}
This shows the ability of our framework to capature the phenomenological description of the drag crisis observed as $Re$ is increased from $160,000$ to $10^6$. A more detailed analysis of drag crisis is required to answer some of the key questions from the flow physics perspective. This includes identification of features that trigger drag crisis, and characterization of transition from high drag (subcritical) to low drag (supercritical) state (into sudden vs continuous transition). These questions are beyond the scope of the current methods paper, but we anticipate addressing these intriguing questions in subsequent work.
\end{remark}

\subsection{Matrix Assembly}\label{results:mat-assembly}
We report on how tensorization can speed up matrix assembly as described in \secref{sec: MatrixAssembly} using linear and quadratic basis function. In order to compare the performance, we considered a canonical lid driven cavity problem in a unit cubic domain at $Re = 100$ on a uniform mesh of $16 \times 16 \times 16$. Table \ref{tab:matAssembly} shows the comparison of the time for matrix assembly and subsequent solve time on TACC \Stampede~SKX and KNL node and \Frontera~CLX processors. We observe a significant speedup in the matrix assembly process across all these platforms, with more substantial speedups for higher order basis functions.
\input{Data/MatrixAssemblyData}

\subsection{\In - \Out Test}\label{Sec:InOutTest}
Here, we show the performance gain by using the normal based \In-~\Out test as compared to the ray-tracing for three different complex geometries: Stanford bunny, Stanford dragon (available at \citep{levoy2005stanford}) and truck (our target problem, discussed later in \secref{sec: truckSimulation}). These examples represent fairly complex geometries creating a diverse set of \Intercepted elements. In order to mimic our simulation scenario, we generated $5^3$ Gauss quadrature points per \Intercepted elements. \tabref{tab:in-out} shows the comparison between the normal based \In-~\Out and ray-tracing. Since both of these algorithms are embarrassingly parallel, the total time reported is the cumulative sum of the time spent by each processor. We make several observations: (a) on average, normal based \In-~\Out test is $1000\,\times$ faster than ray-tracing; (b) as the mesh resolution improves, the fraction of octants needing ray-tracing decreases. This is because any sharp corners (where our normal \In-~\Out test fails) now reside within fewer elements; and (c) by combining ray-tracing with the normal based \In-~\Out we see an overall speedup by at least $6\,\times$ (for the finest resolution) without compromising on accuracy.
\input{Data/in-outTable}
\subsection{Adaptive Quadrature}
In order to test the impact of adaptive quadrature on the flow physics, we consider a benchmark problem of flow past a sphere at $Re=100$. The simulation is performed on a fixed mesh ($R_{bkg} = 5, R_{wake} = 5, R_{bdy}=6$), but with increasing levels of adaptive quadrature of \Intercepted elements. \figref{fig:AdapQuadDragandTime} shows the convergence of drag coefficient and associated percentage error to its reference value with increasing number of quadrature points in the \Intercepted elements. This indicates that adaptive quadrature is essential to accurately model \Intercepted elements. Our numerical experiments suggest that two levels of quadrature splitting (percentage error $\sim$ 1\%) produces converged results, and further increase in quadrature level does not give any significant improvement in quantitative prediction of aerodynamic coefficient, which was also reported by \citet{xu2016tetrahedral} on unstructured meshes. With increase in the quadrature levels, the associated cost with increased Gauss points in \Intercepted elements leads to a load imbalance, significantly increasing the time-to-solve. This suggests the need for weighted partitioning.
\input{Plots/AdaptiveQuadratureDragandTime}
\subsection{Weighted Partitioning}\label{Sec:WeightedPartioning}
Here, we demonstrate the advantage of using weighted partitioning using flow past a sphere problem (at $Re = 100$) as a benchmark. \eqnref{eq: Wpart} provides a theoretical estimate of the weights. This requires identification of the ratio, $\frac{T_V}{T_S}$, which can in principle be evaluated in an architecture specific way. Alternatively, simple Monte-Carlo sampling of volumetric and surface quadrature points provides a good estimate of this fraction. With this fraction, we accurately identify the weight associated with each octant using \eqnref{eq: Wpart}. Table \ref{tab:wpart} shows the estimated cost by running the simulation on \Stampede~ SKX and KNL nodes. The experiments reveal the relative cost for matrix assembly to be $\approx 3.3$.
\input{Data/wpart}
\input{Plots/WpartPlot}

\figref{fig:WPartPrediction} compares the predicted weight (using \eqnref{eq: Wpart}) with the experimental weight, \added{averaged over 10 runs}, for two different meshes and surface discretization. We can see that the predicted weight is in good agreement with the actual weight. The estimated ratio of $\frac{Tv}{Ts}$, on a specific architecture only depends on the type of PDE being solved (which governs the number of FLOPS and data movement across memory hierarchy required per volume and surface quadrature points) and is independent of the number of mesh elements, the nature of surface discretization, as well as the parameters (like Reynolds numbers). Further, we utilize this model to show the impact of correct load balancing on the actual solve time. 

\input{Plots/WpartComparNew}
\input{Data/Wpart/FullRun/dataWpartRun}
\input{Plots/scaling} 
\figref{fig:WpartRun} shows the element and weight distribution across different processors. In case of equal partition, the number of elements are equally distributed across different processors. This results in load imbalance which is circumvented by making use of weighted partitioning. The weighted partitioning ensures that the elements are distributed such that each processor receives equal amount of work. \tabref{tab:wpart-comparison} shows the overall speedup achieved in the assembly time and total solve time. We observe that the correct distribution of work can help to achieve a substantial speedup. \eqnref{eq: Wpart} generalizes to moving meshes, where after each remeshing step, re-partitioning is performed. The enumeration of the viable Gauss points and surface elements in each \Intercepted element is performed in linear time, involving a single pass over the elements.



\subsection{Scaling Studies}\label{Sec:Scaling}
Finally, in this section, we report on the scaling behaviour of our solver by simulating the flow across a sphere case discussed in detail for validating the numerical method in \secref{sec:sphere}.
\subsubsection{\textbf{Strong Scaling}}
\input{Data/scalingDomain}
 For strong scaling, we consider four different adaptive meshes: \textsc{M1} consisting of 76,868 elements, \textsc{M2} consisting of 519,800 elements, \textsc{M3} consisting of around 4 million elements and \textsc{M4} consisting of around 15 million elements. The maximum refinement region in all the meshes was near the sphere region. Table~\ref{tab: scaling-domain} shows the refinement level in different regions for the scaling studies for different mesh. The sphere surface discretization for carrying out the integration of weak boundary conditions was kept constant and comprises of around 0.3 million triangle elements. Surface discretization of the sphere was chosen to ensure that we have at least 3-4 Gauss quadrature points per \Intercepted elements at the finest octree mesh resolution (Mesh \textsc{M4}). We used linear basis function for spatial discretization. 

\figref{fig:Scaling} shows the strong scaling result for the different meshes on TACC's \Frontera.\footnote{We have shown the scaling to the maximum allocation that we had on \Frontera~ supercomputer.} We use a bi-conjugate gradient solver with additive Schwartz preconditioner. The simulation was carried out for 10 time steps with $\Delta t $ of 0.01. In total there were 21 rounds of non-linear solve, which resulted in 21 rounds of matrix assembly and 31 rounds of vector assembly. The average per iteration matrix and vector assembly time is reported. 
The same results are also presented in \figref{fig: RelativeSpeedup} as the relative speedup for different problem sizes. 
The results reveal near ideal scaling till the grain size (the number of octants per processor) is around 320. Since we are solving for 4 degrees of freedom per node, this translates to roughly 5000 degrees of freedom per processor. 
At smaller grain sizes, the amount of data being exchanged with neighboring processor increases compared to the amount of work being done. This makes it difficult to effectively overlap communication and computation, leading to a loss in scalability. 
We also observed that the number of linear solve iterations increases with the number of processors, primarily due to block preconditioning, which additionally contributes to the deviation from the ideal scaling. Finally, at small grain sizes, load balancing for the intercepted elements becomes challenging, and contributes to the loss of scalability. Note that for most practical problems, we will be operating with much larger grain sizes avoiding these challenges. These extreme strong scaling results are presented to motivate additional research into efficient preconditioners and load-balancing techniques for IMGA.

\subsubsection{\textbf{Weak Scaling}}
For weak scaling, we considered 3 different adaptive meshes with number of elements per processor varying from $\sim 500$ to $\sim 2000$. Similar to the strong scaling case, the discretization for the object is kept constant. \figref{fig: WeakScaling} plots the results of our weak scaling experiments. We compare the weak scaling result for the matrix assembly, vector assembly, the time taken for each \textsc{ksp} iteration\footnote{\textsc{ksp} is the Krylov subspace context used in \petsc~for solving linear systems.} iteration and total solve time. We achieve excellent weak scaling performance, with similar trends as the strong scaling. Overall, scaling is better at larger grain sizes, especially for matrix and vector assembly. We achieve a weak scaling efficiency $\sim 0.5$ for total solve time while increasing the number of processor and problem size by a factor of 16. The overall time to solve suffers because of degradation of the preconditioner and poor conditioning of the matrix, which results in increase in number of KSP iteration with problem size.

\begin{remark}
We hypothesize that a multilevel preconditioner can improve the scalability (both strong and weak) as the number of iteration count remains approximately constant with increase in the number of processor and the problem size. However, in our experiments using Algebraic Multigrid (AMG) via \petsc~interface, while we observed a constant number of iteration with increase in the number of processor as well as mesh independence (number of iterations remaining approximately constant with increase problem size), the setup costs for AMG are high and exhibit poor scalability for large number of processor counts (especially $>$ 128 ~\Frontera~ nodes $\sim$ 7000 processors) (\secref{sec:AMG}). Currently, for our target problems, we achieved a better time to solve using a single level Additive Schwarz preconditioner. The use of Geometric multigrid (GMG) should address the high setup costs associated with AMG, as well as its scalability. This is left as future work.
\end{remark}

\section{Flow past a complex geometry: Semi-trailer truck}\label{sec: truckSimulation}
In this section, we illustrate the utility of our framework on a practical application problem. We explore the flow physics across a realistic semi-truck geometry travelling at 65 MPH (a Reynolds number of $30 \times 10^6$). Then, we quantify the advantage of platooning multiple semi-trucks which is one of the compelling fuel saving features of next-generation autonomous vehicles. Compared with modern motor vehicles, commercial vehicles are aerodynamically inefficient due to their bulky design. Any reduction of drag can lower the cost of fuel consumption and thus possess potential economic and environment benefit~\citep{Takizawa15bd}. We demonstrate the capability of our framework to investigate the effect of individual parts of a semi-trailer truck on the drag as well as the platooning effect of two trucks. 

\input{Plots/FrontTruckStreamLine}

\subsection{Geometry of the Truck}
The geometry of semi-trailer consists of several parts such as tractor, trailer, tanks, tires and axis. Each part is verified to be a water-tight triangulated manifold. The non-dimensional length of the truck is $1$ (normalized by the truck length, which is 15 m). The Reynolds number of the problem, estimated based on the truck length and cruising at 65 mph, is around $30 \times10^6$. We solve in a moving reference frame, where the truck is stationary and air is moving past at 65 mph.

\subsection{Computational Domain and Boundary Conditions}
The computational domain has a dimension of $16\times2\times2$, with inlet velocity set at $1$ and the outlet pressure set at $0$. The surrounding walls are no-slip walls moving at the same speed as the incoming flow. The truck is positioned $5$ unit lengths behind the inlet. No-slip boundary condition is weakly imposed on the surface of the truck. We model the rotating wheels by enforcing a no-slip condition corresponding to the wheels rotating with an angular speed of $1/r$, producing a linear speed of 1.0. Additional details are presented in \ref{sec:AEtruck-simulation}.

\subsection{Mesh Generation}
The mesh is refined adaptively using multiple refine regions around the truck and in the wake region. The base level of refinement is set to $8$ (i.e. $2^8$ divisions along a dimension) and the smallest element around the truck has refinement level of $12$ ($2^{12}$ divisions along a dimension). This level of refinement is chosen according to the Taylor length scales at this Reynolds number, $Re = 30 \times 10^6$. This resulted in around $3.1$ million elements. A time step ($\Delta t$) of 0.00125 is chosen for the simulation, resulting in $CFL$ number ranging from $0.02$ to $0.32$ from the largest element to the smallest element.
The simulation is started with a lower $Re$, larger $\Delta t$ and a less refined mesh. As the solution converges, we ramp up the $Re$, decrease $\Delta t$ and further refine mesh near the region of interest. This is done to remove the initial transient quickly. The removal of transient is a major bottleneck in high fidelity CFD simulations; with transient removal taking several days of simulation time before any statistically reliable data can be collected. The octree based framework provides a principled approach to remove these transients by performing simulation starting from a relatively coarse mesh, and successively refining the mesh. Here, to remove transients, we used \Stampede~ SKX and KNL processors with number of processors ranging from 192 to 2176. Once the initial transient was removed, the simulation was carried out on 64 \Frontera~ nodes with 3584 processors. In approximately 4 hours, we were able to collect the statistics for about 3.5 seconds.

\input{Plots/TruckDragSingle}
\input{Plots/TwoTruckCp}
\input{Plots/TwoTruckDrag}

\subsection{Flow Quantities of Interest}

\figref{fig: TruckDragSingle} shows the (coefficient of) drag as a function of time. Note that we are plotting for time after removal of the initial transients, when statistically consistent results are produced. The truck head contributes the most to the drag force, whereas the trailer actually contributes a small negative quantity. This is mainly due to the pressure difference at the front and back of the trailer surface. The average non-dimensional drag coefficient $C_d$ comes around to be 0.695. The reference area $A$ is chosen to be the projected frontal area of the truck. This result is comparable to the previous conducted experimental studies ~\cite{Gotz1987295,chowdhury2013study,Englar2001} which had reported $C_d$ in the range of 0.6--0.9 for heavy vehicles. The previous numerical result reported the drag coefficient to be about 0.57 for RANS (Reynolds Averaged Navier--Stokes) simulation and 0.62 for DES (Detached Eddy Simulation) simulations for the truck travelling at the speed of 55 MPH. \citep{viswanathan2019platooning}. It must be noted that numerical results based on RANS and DES are susceptible to the proper choice of parameters for wall treatment; while our framework does not rely on any additional treatments. This makes our approach significantly more robust than existing state-of-art approaches.

\subsection{Flow Past Multiple Complex Objects: Platooning of Semi-trailer Trucks}

Studies have shown that by platooning large, blunt vehicles, the overall fuel saving can be significant. This is especially significant for autonomous vehicles, where a platoon of trucks can operate safely and efficiently (with nearly 25\% improved efficiency \citep{guttenberg2017evaluating,viswanathan2019platooning}). We explore this concept using the detailed LES simulations afforded by this framework. We placed two identical semi-trailer trucks one truck-length $D$ apart, with the computation domain and boundary conditions identical to previous simulation. We can see that the The full scale two-truck simulation resulted in a mesh with about 6 million elements. The complete simulation was carried out on 128 ~\Frontera~ CLX nodes (7168 processor) for about 16 hours \footnote{This time includes the time taken to remove the transients for two truck simulation.} to collect the data over 12 seconds, after removal of transients. Leveraging the ability to refine and coarsen the grid, we used a checkpoint solution from the one-truck simulation as the initial guess, which resulted in substantial reduction in overhead for transient removal during the two truck simulations. \figref{fig: TwoTruckDrag} shows the $C_d$ history for major components of two trucks. The averaged $C_d$ for the second truck is $0.475$, showing a $32\%$ decrease in total drag. This reduction in drag is comparable with the ones reported in literature where a reduction of 30\% is observed using DES simulation~\citep{viswanathan2019platooning}, as well as a full scale experimental study~\citep{torabi2018fuel}. \figref{fig:front-pageStreamlines} shows the streamlines for the flow over two trucks. We can see the streamlines of the first truck is effecting the flow around the second truck. It is clear from \figref{fig: TwoTruckPressureTractor} that the main source of drag reduction is due to the lowered drag of the tractor head. The figure plots the non--dimensional pressure coefficient $C_p$ at the centerline of the two tractors. The front section of the trailing truck has a lower pressure than the leading truck, with larger difference around the lower half of the truck. On the surface of the trailer, shown in \figref{fig: TwoTruckPressureTrailer}, the platooning truck shows similar pressure magnitude with both the front and back surface while the leading truck clearly have higher pressure at the front surface of the trailer. Additional analysis of the flow physics is reported in \ref{sec:AEtruck-simulation}


%% file: Plots/DragCompar.tex

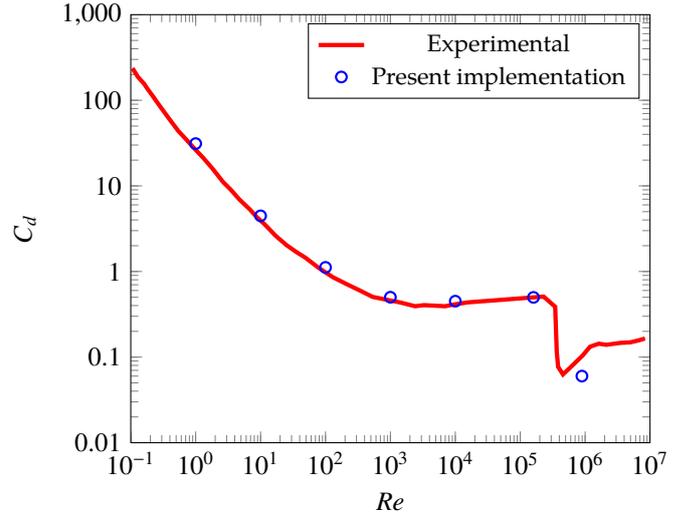
\begin{figure}[t!]
\centering
\begin{tikzpicture}
      \begin{loglogaxis}[
          width=0.95\linewidth, 
          xlabel=$Re$, 
          ylabel=$C_d$,
          log y ticks with fixed point,
        ymax = 1000,
          ymin =  1E-2,
          xmin = 1E-1,
          xmax = 1E7,
          x tick label style={rotate=0,anchor=north} 
        ]
        \addplot[mark=none,color = red,ultra thick] 
        table[x expr={\thisrow{Re}))},y expr=\thisrow{Drag},col sep=comma]{Data/SphereDrag.txt};
        
        \addplot[only marks,mark=o,color = blue,thick] 
        table[x expr={\thisrow{Re}))},y expr=\thisrow{Drag},col sep=comma]{Data/computation.dat};

        
        
         \legend{\small {Experimental},\small {Present implementation} }
        \end{loglogaxis}
        \end{tikzpicture}
        \caption{\textit{Drag crisis:} Variation of drag coefficient $C_d$ with $Re$ for flow past a sphere compared with experimental data \cite{namburi2016crystallographic,achenbach1972experiments}. Notice that the $\mathrm{x}$-axis is $\log(Re)$ representing variation for a wide range of $Re$ number spanning multiple regimes (from laminar to fully developed turbulence) each with distinct flow physics. The framework demonstrates good comparison with experimental data in all these regimes without any special adjustments to the numerical scheme.  Moving towards the right on a log scale as $Re$ is increased becomes increasingly computationally expensive, e.g adding another point on the right of the last point would increase the Reynolds number ten-folds and would require a decrease in finest element size by 3 times.} 
          \label{fig: Drag}
\end{figure}

%% file: Plots/SpherePlot.tex

\begin{figure*}[t!]
\centering
\begin{subfigure}{.33\textwidth}
    \includegraphics[width=1.0\linewidth]{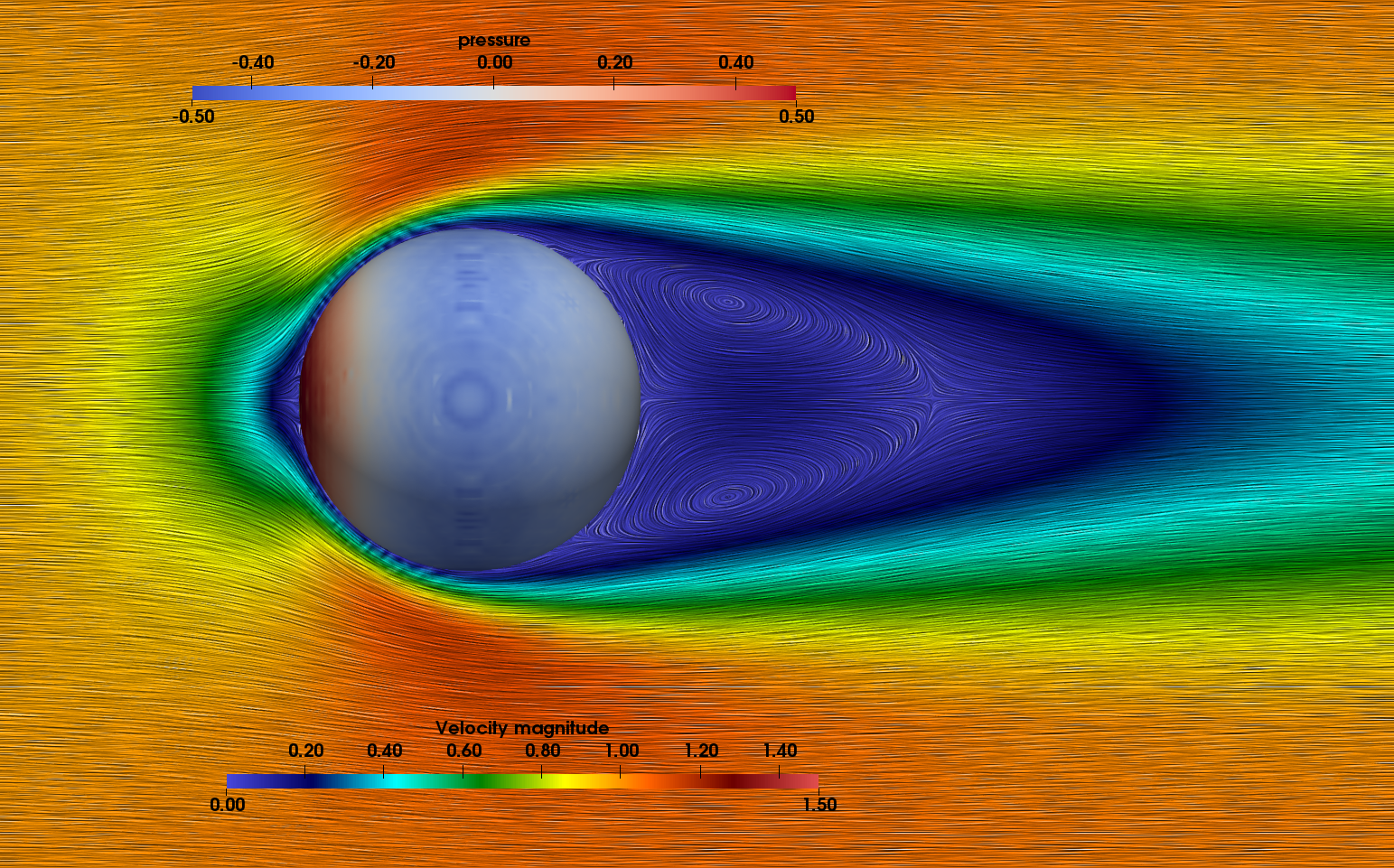}
    \caption{$Re$ = 100}
\label{fig: Re100}
\end{subfigure}
\begin{subfigure}{.33\textwidth}
    \includegraphics[width=1.0\linewidth]{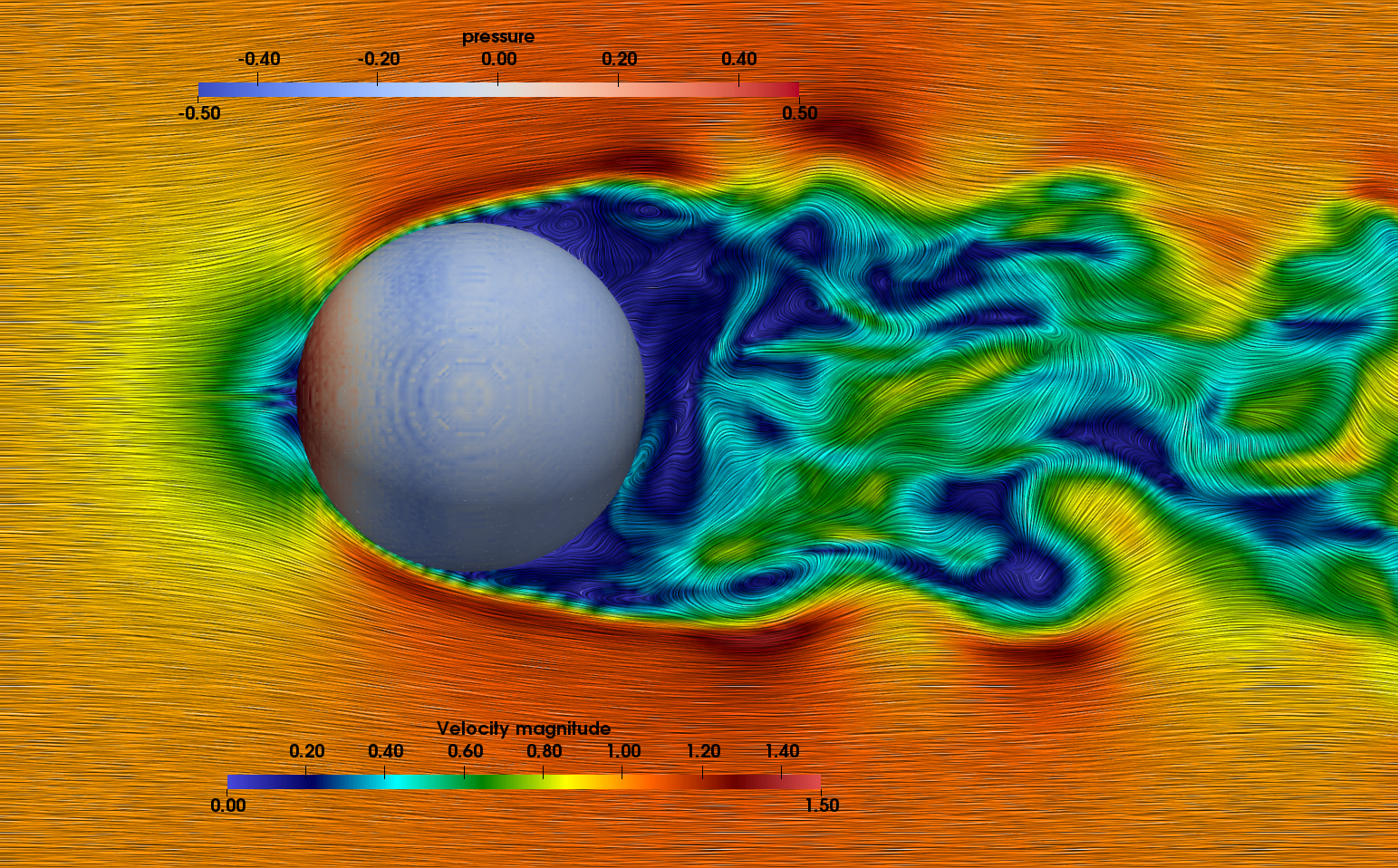}
    \caption{$Re$ = 160,000}
    \label{fig: Re160K}
\end{subfigure}
\begin{subfigure}{.33\textwidth}
    \includegraphics[width=1.0\linewidth]{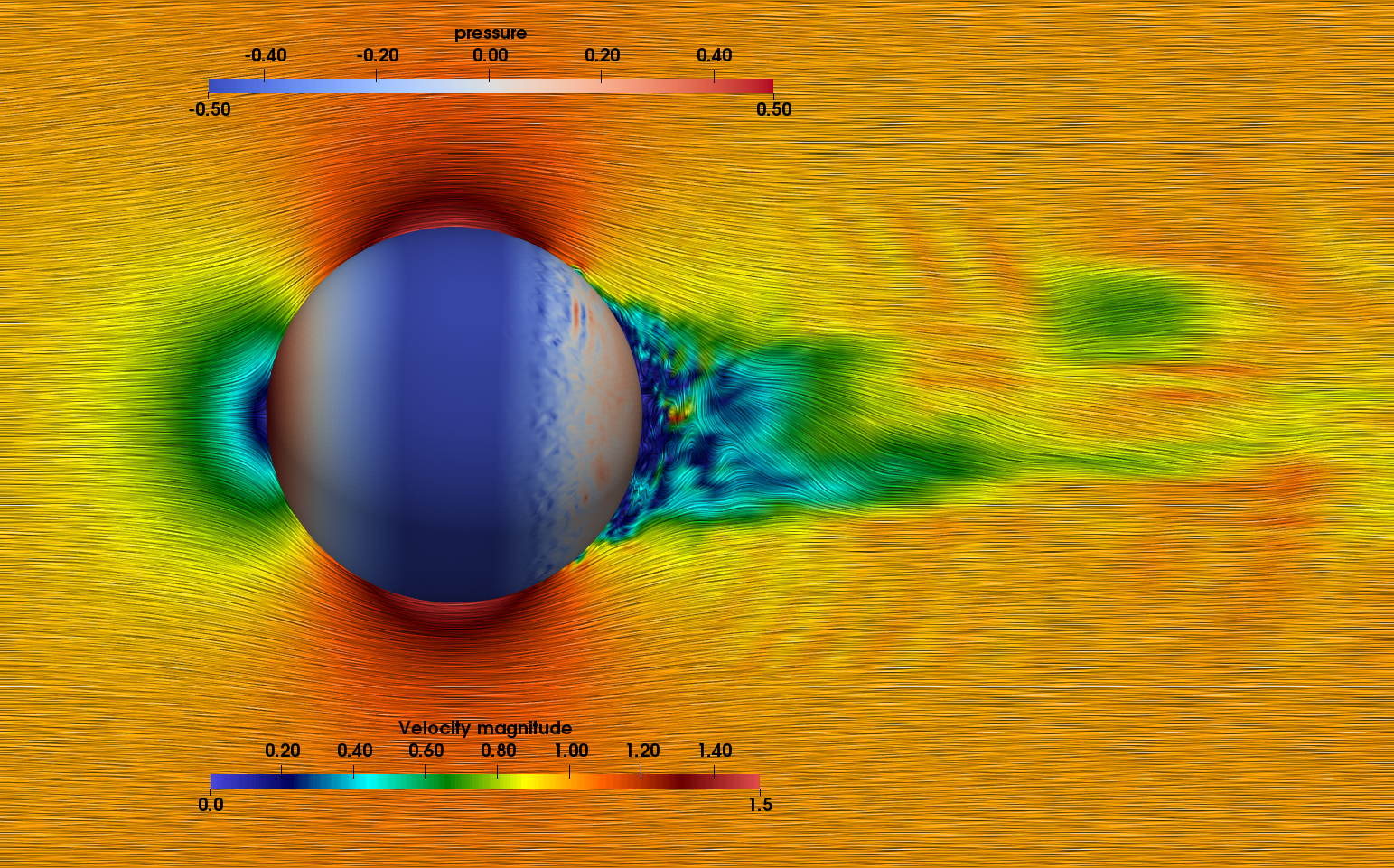}
    \caption{{$Re$ = 1,000,000}}
    \label{fig: Re1M}
\end{subfigure}
\vspace{1mm}
\caption{The wake structures and pressure distribution on sphere at different Reynolds number. At low $Re$, the wake remains axisymmetric. As $Re$ increases, it starts shedding and boundary layer becomes turbulent. The picture gives the phenomenological description of the drag crisis. The drag crisis is evident by noticing the wake structure as it changes from being divergent at $Re=160,000$ (high drag state) to being convergent at $Re = 10^6$ (low drag state). At the same time, we observe a high pressure region being developed behind the sphere. The development of this high pressure zone is attributed to the low drag state. Below the drag crisis, the flow separates at an angle smaller than $90^0$  in the hemisphere facing the flow but above the drag crisis, the separation angle is pushed backward to the hemisphere pointing away from the flow.}
\label{fig: wake}
\end{figure*}

%% file: Data/MatrixAssemblyData.tex
{\renewcommand{\arraystretch}{1.5} 
\begin{table}[h!]
\centering
\captionsetup{font=small}
\resizebox{0.95\linewidth}{!}{%
\begin{tabular}{|c|c|c|c|c|c|c|c|}
\hline
\multicolumn{2}{|c|}{\multirow{2}{*}{}} & \multicolumn{3}{c|}{Assembly Time} & \multicolumn{3}{c|}{Total Time} \\ \cline{3-8} 
\multicolumn{2}{|c|}{}                  & GP (s)         & MM (s)        & Speedup   & GP (s)        & MM (s)       & Speedup \\ \hline
{ \parbox[t]{2mm}{\multirow{3}{*}{\rotatebox[origin=c]{90}{Linear}}}}        & SKX    & 1.0617     & 0.7754    & 1.37      & 4.2598    & 3.9201    & 1.10    \\ \cline{2-8} 
                               & KNL    & 5.4996     & 3.775     & 1.46      & 7.83192   & 6.33589   & 1.24    \\ \cline{2-8} 
                               & CLX    & 0.80538    & 0.562     & 1.43      & 10.33978  & 10.17841  & 1.02    \\ \hline
{ \parbox[t]{2mm}{\multirow{3}{*}{\rotatebox[origin=c]{90}{Quadratic}}}}  & SKX    & 71.132     & 10.245    & 6.94      & 234.202   & 175.03    & 1.34    \\ \cline{2-8} 
                               & KNL    & 385.59     & 40.409    & 9.54      & 664.05    & 346.15    & 1.92    \\ \cline{2-8} 
                               & CLX    & 59.356     & 8.4647    & 7.01      & 212.77    & 164.46    & 1.3     \\ \hline
\end{tabular}%
}
\caption{Comparison of time on different computing environment for matrix assembly and subsequent solve time for linear and  quadratic basis function on \Stampede~SKX and KNL nodes, and \Frontera~CLX nodes. GP indicates the matrix assembly by looping over the individual Gauss points whereas MM indicates the assembly by imposing FEM operator as matrix--matrix multiplication.}
\vspace{0.1in}
\label{tab:matAssembly}
\end{table}
}

%% file: Data/in-outTable.tex
{\renewcommand{\arraystretch}{1.2} 
\begin{table*}[t!]
\captionsetup{font=small}
\resizebox{\linewidth}{!}{%
{
\begin{tabular}{|c|c|c|c|c|c|c|c|c|l|c|c|}
\hline
\multirow{3}{*}{}       & \multirow{3}{*}{\begin{tabular}[c]{@{}c@{}}Grid \\ Size\end{tabular}} & \multirow{3}{*}{\begin{tabular}[c]{@{}c@{}}Total  \\ points\end{tabular}} & \multicolumn{7}{c|}{Proposed scheme}                                                                                                                                        & \multirow{2}{*}{Ray--tracing} & \multirow{3}{*}{Speedup} \\ \cline{4-10}
                        &                                                                       &                                                                           & \multicolumn{3}{c|}{Normal--based} & \multicolumn{3}{c|}{Ray--tracing} & \multicolumn{1}{c|}{\multirow{2}{*}{\begin{tabular}[c]{@{}c@{}}Total\\ time (s)\end{tabular}}} &                                &                          \\ \cline{4-9} \cline{11-11}
                        &                                                                       &                                                                           & $P_s$    & Time (s)   & $T_{avg}$(s)     & $P_f$   & Time (s)   & $T_{avg}$(s)    & \multicolumn{1}{c|}{}                                                                            & Total time (s)                 &                          \\ \hline
 {\parbox[t]{2mm}{\multirow{3}{*}{\rotatebox[origin=c]{90}{Bunny}}}}  & $64^3$                                                                    & 1,193,640                                                                   & 0.79        & 0.996      & 1E-6     & 0.21        & 697.012    & 0.003   & 698.008                                                                                          & 3247.6                         & 4.66                     \\ \cline{2-12} 
                        & $128^3$                                                                   & 4,149,957                                                                   & 0.87        & 2.205      & 6E-7     & 0.13        & 1496.81    & 0.003   & 1499.01                                                                                          & 11264.5                        & 7.52                     \\ \cline{2-12} 
                        & $256^3$                                                                   & 10,317,097                                                                  & 0.91        & 4.234      & 4E-7     & 0.09        & 2555.78    & 0.003   & 2560.02                                                                                          & 28156.6                        & 11                       \\ \hline
{ \parbox[t]{2mm}{\multirow{3}{*}{\rotatebox[origin=c]{90}{Dragon}}}} &  $64^3$                                                                      &      837,263                                                                     &     0.56        &     0.558       & 5E-7         &             0.44 &     1086.05       &     0.003    &  1086.32                                                                                                &         2213.46                       &   2.03                       \\ \cline{2-12} 
                        & $128^3$                                                                    &       3,351,110                                                                    &    0.74         &  0.624          &      2.53E-7    &    0.26         &    2615.81        &    0.003     &    2616.43                                                                                              &                       8885.58           &                3.4        \\ \cline{2-12} 
                        & $256^3$                                                                   & 12,120,591                                                                  & 0.85        & 1.382      & 1.3E-7   & 0.15        & 6096.96    & 0.003   & 6098.34                                                                                          & 35392.4                        & 5.81                     \\ \hline
{\parbox[t]{2mm}{\multirow{3}{*}{\rotatebox[origin=c]{90}{Truck}}}}  &  $64^3$   & 383,474 & 0.64 & 0.99 & 4E-6 & 0.36 & 621.469 & 0.004 & 622.459 & 1638.21 & 2.63 \\ \cline{2-12} 
 & $128^3$  & 1,602,839 & 0.80 & 1.84 & 1E-7 & 0.20 & 1724.55 & 0.005 & \multicolumn{1}{c|}{1726.39} & 7468.15 & 4.33 \\ \cline{2-12} 
 & $256^3$ & 6,205,213 & 0.89 & 3.535 & 6.3E-7 & 0.11 & 3816.79 & 0.005 & 3820.33 & 34675.8 & 9.08 \\ \hline                                                    
\end{tabular}%
}
}
\caption{Comparison of the performance of the proposed \In - \Out test with the conventional ray-tracing for three different complicated geometries. $P_s$ represents the success fraction of normal based test and $P_f$ represents the failure fraction for which we resort to ray-tracing. $T_{avg}$ represents the average time per evaluation of the normal based and ray-tracing.} 
\label{tab:in-out}
\end{table*}
}

%% file: Plots/AdaptiveQuadratureDragandTime.tex
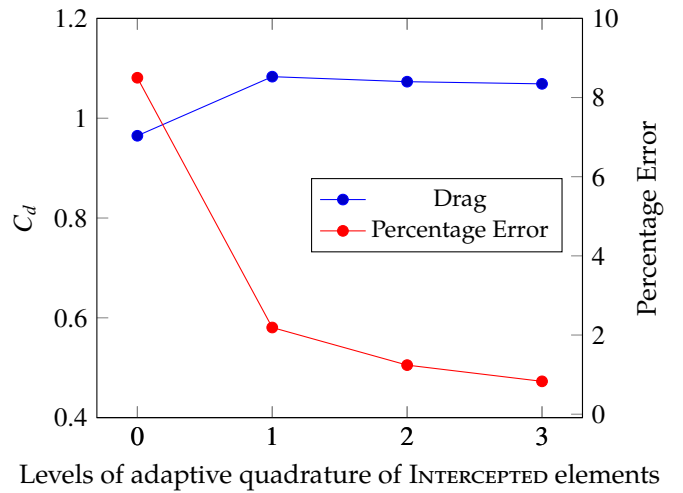
\begin{figure}[t!]
\centering
\begin{tikzpicture}
\pgfplotsset{
xtick={0,1,2,3}
}
      \begin{axis}[
          width=0.9\linewidth, 
          xlabel=Levels of adaptive quadrature of \Intercepted elements, 
          ylabel=$C_d$,
          ymin = 0.4,
          ymax = 1.2,
          /pgf/number format/precision=5,
          ytick pos=left,
          x tick label style={rotate=0,anchor=north} 
        ]
        \addplot
        table[x expr={\thisrow{Level}))},y expr=\thisrow{Drag},col sep=space]{Data/quadratureData.dat};
        \label{AdapQuadDrag}
      \end{axis}
      \begin{axis}[
          axis y line*=right,
          width=0.9\linewidth, 
          ylabel=Percentage Error,
          ymax = 10,
          legend style={at={(0.7,0.6)},anchor=north,legend columns=1}, 
      ]
      \addlegendimage{/pgfplots/refstyle=AdapQuadDrag}\addlegendentry{\small {Drag}}
        \addplot[mark=*,red]
        table[x expr={\thisrow{Level}))},y expr=\thisrow{PercentError},col sep=space]{Data/quadratureData.dat};
        \addlegendentry{\small {Percentage Error}}
      \end{axis}
\end{tikzpicture}
\caption{Drag and percentage error for different levels of adaptive quadrature. We see that with increase in the number of quadrature point for \Intercepted nodes, $C_d$ converges to the experimental observed value of 1.06-1.096 at $Re = 100$ \cite{xu2016tetrahedral}. Reference value of 1.06 ~\cite{marella2005sharp} was chosen to compute percent error.}
\label{fig:AdapQuadDragandTime}
\end{figure}

%% file: Data/wpart.tex
{\renewcommand{\arraystretch}{1.2} 
\begin{table}[h!]
\centering
\captionsetup{font=small}
\begin{tabular}{|c|c|c|c|c|c|c|}
\hline
\multirow{2}{*}{} & \multicolumn{3}{c|}{Matrix Assembly} \\ \cline{2-4} 
 & Tv & Ts & Ratio \\ \hline
SKX & 2.72E-5$\pm$7.5E-6 & 8.4E-6 $\pm$1.85E-6 & 3.23 \\ \hline
KNL & 1.6E-4 $\pm$ 2E-5 & 4.83E-5 $\pm$ 1.51E-5 & 3.32 \\ \hline
\end{tabular}%
\caption{The ratio $\frac{T_v}{T_s}$ estimated by running the simulation on \Stampede~ SKX and KNL nodes.}
\label{tab:wpart}
\end{table}
}

%% file: Plots/WpartPlot.tex
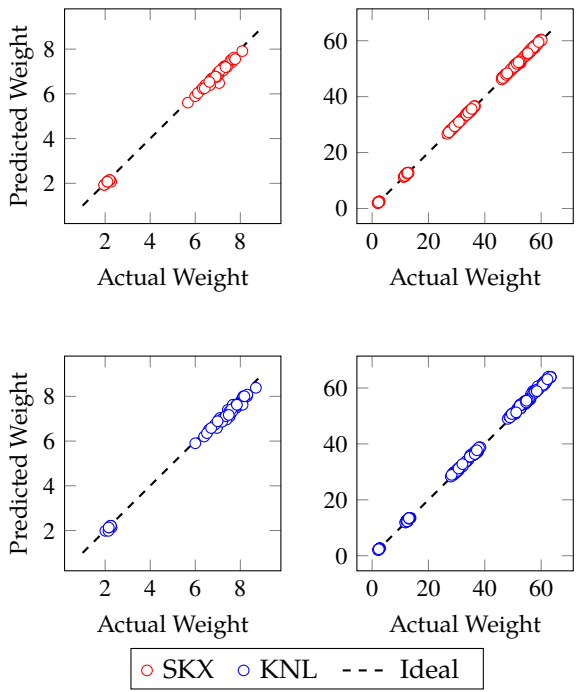
\begin{figure}[b!]
\centering
\begin{tikzpicture}
  \begin{groupplot}[
      group style={group size=2 by 2,ylabels at=edge left,vertical sep=50pt},
      ylabel style={text height=0.02\textwidth,inner ysep=0pt},
      xlabel= Actual Weight, 
      ylabel= Predicted Weight,
      height=0.5*\linewidth,width=0.5\linewidth,/tikz/font=\small,
    ]
    \nextgroupplot[]%
     \addplot[only marks,mark options={solid, red, fill=white}]
        table[x expr={\thisrow{AWeight}},y expr={\thisrow{PWeight})},col sep=space]{Data/Wpart/SKX/MatLevel4_Avg.txt};
          \addplot [no markers, color=black, thick,dashed ] table {
        1 1
        9 9
        };
         \coordinate (top) at (rel axis cs:1,0);

        \nextgroupplot[]%
    \addplot[only marks,mark options={solid, red,fill =white}]
        table[x expr={\thisrow{AWeight}},y expr={\thisrow{PWeight})},col sep=space]{Data/Wpart/SKX/MatLevel5_Avg.txt};
        \label{plot:part-SKX}
          \label{plot:part-KNL}
          \addplot [no markers, color=black, thick,dashed ] table {
        1 1
        65 65
        };
   \coordinate (bot) at (rel axis cs:1,0);
    \nextgroupplot[]%
        \addplot[only marks,mark options={solid, blue, fill=white}]
          table[x expr={\thisrow{AWeight}},y expr={\thisrow{PWeight})},col sep=space]{Data/Wpart/KNL/MatLevel4_Avg.txt};
          \addplot [no markers, color=black, thick,dashed ] table {
        1 1
        9 9
        };
         \coordinate (top) at (rel axis cs:1,0);

        \nextgroupplot[]%
        \addplot[only marks,mark options={solid, blue,fill=white}]
          table[x expr={\thisrow{AWeight}},y expr={\thisrow{PWeight})},col sep=space]{Data/Wpart/KNL/MatLevel5_Avg.txt};
          \label{plot:part-KNL}
          \addplot [no markers, color=black, thick,dashed ] table {
        1 1
        65 65
        };
        \label{plot:part-ideal}
   \coordinate (bot) at (rel axis cs:1,0);
     
  \end{groupplot}
  \path (top|-current bounding box.south)--
        coordinate(legendpos)
        (bot|-current bounding box.south);
  \matrix[
      matrix of nodes,
      anchor=north,
      draw,
      inner sep=0.2em,
    ]at([xshift=-10ex]legendpos)
    { \ref{plot:part-SKX}& SKX &[5pt]
      \ref{plot:part-KNL}& KNL &[5pt]
      \ref{plot:part-ideal}& Ideal &[5pt]\\};
\end{tikzpicture}
\caption{\added{Comparison of predicted weight and  actual weight for \Intercepted elements for two different spatial  and surface mesh on \Stampede~ SKX and KNL nodes. The left panel corresponds to the $16^3$ uniform grid with surface discretization comprising of 3000 triangles and right panel correspond to the $32^3$ uniform grid and surface discretization with 82000 triangles.}}
\label{fig:WPartPrediction}
\end{figure}{}

%% file: Plots/WpartComparNew.tex
\begin{figure}[t!]
\centering
\begin{tikzpicture}
  \begin{groupplot}[
      group style={group size=2 by 1,ylabels at=edge left,vertical sep=50pt},
      ylabel style={text height=0.02\textwidth,inner ysep=0pt},
      xlabel= Number of processor, 
      ylabel= Relative fraction,
      height=0.5*\linewidth,width=0.5\linewidth,/tikz/font=\small
    ]
    \nextgroupplot[ymode = log,ytick={1,0.1,0.01},
      yticklabels={$1$,$0.1$,$0.01$},title=Equal partition]%
     \addplot[mark=none,color=blue,ultra thick]
        table[x expr={\thisrow{P}},y expr={\thisrow{Weight})},col sep=space]{Data/Wpart/FullRun/InitialWeight.txt};
        \addplot[mark=none,color=red,ultra thick]
          table[x expr={\thisrow{P}},y expr={\thisrow{Elements})},col sep=space]{Data/Wpart/FullRun/InitialWeight.txt};
    \coordinate (bot) at (rel axis cs:1,0);
        \nextgroupplot[ymode=log,title=Weighted partition,  ytick={1,0.1,0.01},
      yticklabels={$1$,$0.1$,$0.01$}]%
    \addplot[mark=none,color=blue,ultra thick]
        table[x expr={\thisrow{P}},y expr={\thisrow{Weight})},col sep=space]{Data/Wpart/FullRun/FinalWeight.txt};
        \label{plot:WpartWeight}
        \addplot[mark=none,color=red,ultra thick]
          table[x expr={\thisrow{P}},y expr={\thisrow{Elements})},col sep=space]{Data/Wpart/FullRun/FinalWeight.txt};
          \label{plot:WpartElements}
    \coordinate (top) at (rel axis cs:1,0);
     
  \end{groupplot}
  \path (top|-current bounding box.south)--
        coordinate(legendpos)
        (bot|-current bounding box.south);
  \matrix[
      matrix of nodes,
      anchor=north,
      draw,
      inner sep=0.2em,
    ]at([xshift=-10ex]legendpos)
    { \ref{plot:WpartWeight}& Weight fraction &[5pt]
      \ref{plot:WpartElements}& Element fraction &[5pt]\\};
\end{tikzpicture}
\caption{Impact of weighted partition: Figure compares how elements and weights are distributed for the standard (left) and weighted (right) partitioning scheme. In equal partition, each processor receives the equal number of element whereas in case of weighted partitioning, the elements are partitioned in such a way that each processor receives nearly identical weight. The fraction are computed by normalizing by the maximum number of elements and weight across all processors.}
\label{fig:WpartRun}
\end{figure}
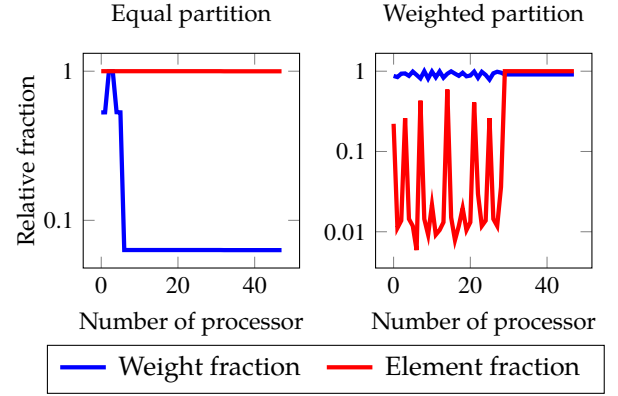{}

%% file: Data/Wpart/FullRun/dataWpartRun.tex
\begin{table}[h!]
\centering
\captionsetup{font=small}
\resizebox{\linewidth}{!}{%
{
\begin{tabular}{|c|c|c|c|}
\hline
 & Equal partition & Weighted Partition & Speedup \\ \hline
\begin{tabular}[c]{@{}c@{}}Matrix  Assembly\end{tabular} & 72.672 s & 28.01 s & 2.60 \\ \hline
\begin{tabular}[c]{@{}c@{}}Solve Time\end{tabular} & 127.4 s & 47.67 s & 2.67 \\ \hline
\end{tabular}%
}
}
\caption{Comparison of matrix assembly and solve time for equal and weighted partitioning}
\label{tab:wpart-comparison}
\end{table}

%% file: Plots/scaling.tex
\begin{figure*}[h!]
\begin{tikzpicture}
	\centering
  \begin{groupplot}[
      group style={group size=3 by 1,ylabels at=edge left,vertical sep=50pt},
      ylabel style={text height=0.02\textwidth,inner ysep=0pt},
      xlabel= Number of processors, 
      ylabel= Time (in s),
      height=0.3*\linewidth,
      width=0.35\linewidth,
      /tikz/font=\small,
      xtick={512,2000,8000,32000},
      xticklabels={$512$,$2K$,$8K$,$32K$},
    ]
    \nextgroupplot[xmode=log,log basis x={2},ymode=log,ytick={0.1,1,10},ymax=10,yticklabels={$0.1$,$1$,$10$},title={\textbf{Matrix Assembly}}]%
     \addplot
         table[x expr={56*\thisrow{Nodes}))},y expr={\thisrow{Jacobian}/\thisrow{MatTimes}},col sep=space]{Data/Scaling/8.dat};
        \addplot
         table[x expr={56*\thisrow{Nodes}))},y expr={\thisrow{Jacobian}/\thisrow{MatTimes}},col sep=space]{Data/Scaling/KSP9.dat};
       \addplot
         table[x expr={56*\thisrow{Nodes}))},y expr={\thisrow{Jacobian}/\thisrow{MatTimes}},col sep=space]{Data/Scaling/KSP10.dat};
         \addplot
         table[x expr={56*\thisrow{Nodes}))},y expr={\thisrow{Jacobian}/\thisrow{MatTimes}},col sep=space]{Data/Scaling/KSP10_new.dat};
    \coordinate (bot) at (rel axis cs:1,0);
    
    \nextgroupplot[xmode=log,log basis x={2},ymode=log,ytick={0.01,0.1,1},yticklabels={$0.01$,$0.1$,$1$},ymin=0.02,title={\textbf{Vector Assembly}}]%
   \addplot
         table[x expr={56*\thisrow{Nodes}))},y expr={\thisrow{VecAssembly}/\thisrow{VecTimes}},col sep=space]{Data/Scaling/KSP8.dat};
        \addplot
         table[x expr={56*\thisrow{Nodes}))},y expr={\thisrow{VecAssembly}/\thisrow{VecTimes}},col sep=space]{Data/Scaling/KSP9.dat};
       \addplot
         table[x expr={56*\thisrow{Nodes}))},y expr={\thisrow{VecAssembly}/\thisrow{VecTimes}},col sep=space]{Data/Scaling/KSP10.dat};
         \addplot
         table[x expr={56*\thisrow{Nodes}))},y expr={\thisrow{VecAssembly}/\thisrow{VecTimes}},col sep=space]{Data/Scaling/KSP10_new.dat};
    \coordinate (top) at (rel axis cs:0,1);
    
    \nextgroupplot[xmode=log,log basis x={2},ymode=log,
    ytick={1,10,100},yticklabels={$1$,$10$,$100$},
    title={\textbf{Solve time}}]%
   \addplot
         table[x expr={56*\thisrow{Nodes}))},y expr={\thisrow{Solve}},col sep=space]{Data/Scaling/KSP8.dat};
         \label{plot:level8}
        \addplot
         table[x expr={56*\thisrow{Nodes}))},y expr={\thisrow{Solve}},col sep=space]{Data/Scaling/KSP9.dat};
         \label{plot:level9}
       \addplot
         table[x expr={56*\thisrow{Nodes}))},y expr={\thisrow{Solve}},col sep=space]{Data/Scaling/KSP10.dat};
         \label{plot:level10}
         \addplot
         table[x expr={56*\thisrow{Nodes}))},y expr={\thisrow{Solve}},col sep=space]{Data/Scaling/KSP10_new.dat};
         \label{plot:level10N}
         \addplot [no markers, color=black, thick, dashed ] table {
         224 4.7237e+02 
         28672 3.690390625 
        };\label{plot:ideal}
        \addplot [no markers, color=black, thick, dashed ] table {
         224  6.9489e+01
         7168 2.17153125
        };
        \addplot [no markers, color=black, thick, dashed ] table {
         112  4.0678e+01
         1792 2.542375
        };
        \addplot [no markers, color=black, thick, dashed ] table {
         1792  4.3468e+02 
         28672 27.1675
        };
    \coordinate (top) at (rel axis cs:0,1);
     
  \end{groupplot}
  \path (top|-current bounding box.south)--
        coordinate(legendpos)
        (bot|-current bounding box.south);
  \matrix[
      matrix of nodes,
      anchor=north,
      draw,
      inner sep=0.2em,
    ]at([yshift=-1ex]legendpos)
    { \ref{plot:level8}& \textsc{M1} &[5pt]
      \ref{plot:level9}& \textsc{M2} &[5pt]
      \ref{plot:level10}& \textsc{M3} &[5pt]
      \ref{plot:level10N}& \textsc{M4} &[5pt]
      \ref{plot:ideal}& ideal &[5pt]\\};
\end{tikzpicture}
\caption{Strong scaling result on TACC \Frontera~on 4 different meshes. The plot shows good scaling until the grain size per processor value of $\approx$ 320 (5000 dof).}
\label{fig:Scaling}
\end{figure*}
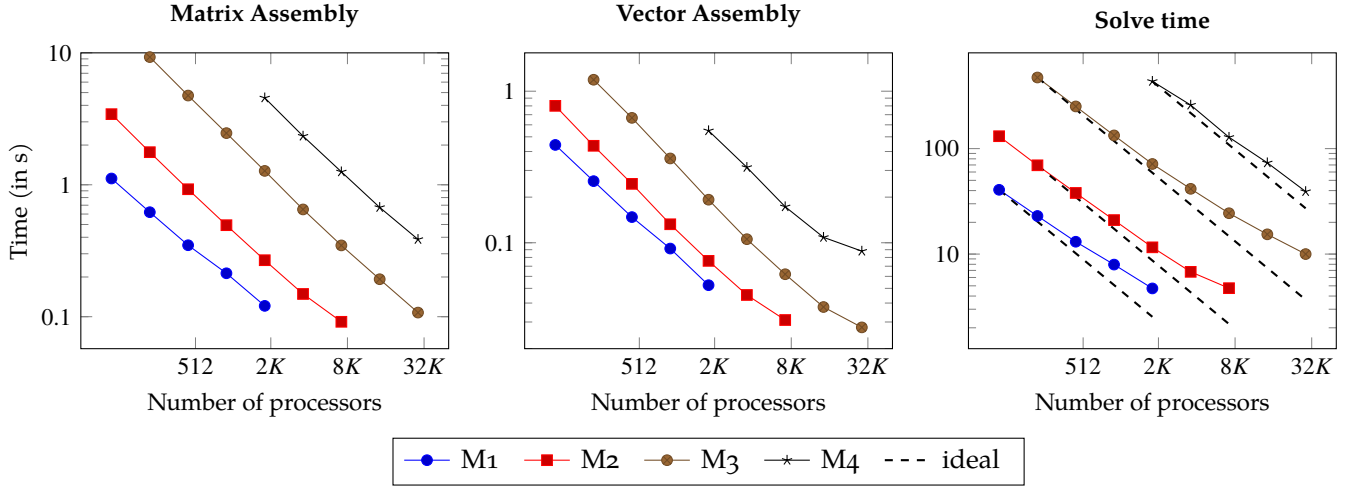{}

\begin{figure}[h!]
\centering
\begin{tikzpicture}
    \begin{loglogaxis}[
          width=0.85\linewidth, 
          height=0.9\linewidth,
          ylabel=Relative speedup, 
          xlabel=Number of Processor,
          legend style={at={(1.2,1.0)},anchor=north west,legend columns=1}, 
          x tick label style={rotate=0,anchor=north}, 
          log basis y={2},
          log basis x={2},
          legend style={at={(0.5,-0.3)},anchor= north,legend columns=5}, 
]
\addplot
     table[x expr={56*\thisrow{Nodes}},y expr={(4.0678e+01/\thisrow{Solve})},col sep=space]{Data/Scaling/KSP8.dat};
     \addplot
     table[x expr={56*\thisrow{Nodes}},y expr={(1.3108e+02/\thisrow{Solve})},col sep=space]{Data/Scaling/KSP9.dat};
 \addplot
    table[x expr={56*\thisrow{Nodes}},y expr={(4.7237e+02*2/\thisrow{Solve})},col sep=space]{Data/Scaling/KSP10.dat};
    \addplot
    table[x expr={56*\thisrow{Nodes}},y expr={(4.3468e+02*16/\thisrow{Solve})},col sep=space]{Data/Scaling/KSP10_new.dat};
    
      \addplot [no markers, color=black, thick,dashed ] table {
        112 1
        28672 256
        };
    \legend{\small \textsc{M1}, \textsc{M2},\textsc{M3}, \textsc{M4},ideal}
  \end{loglogaxis}
  \end{tikzpicture}
  \caption{Relative speedup for solve time  as a function of number of processor for different problem size.}
  \label{fig: RelativeSpeedup}
\end{figure}
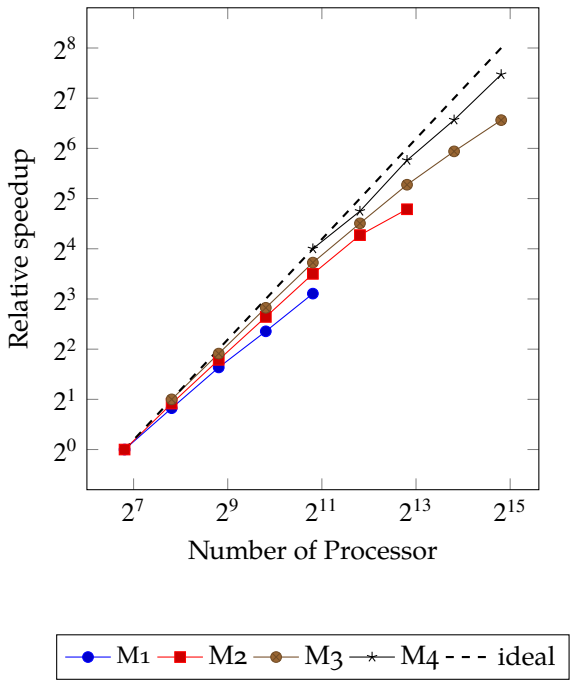{}

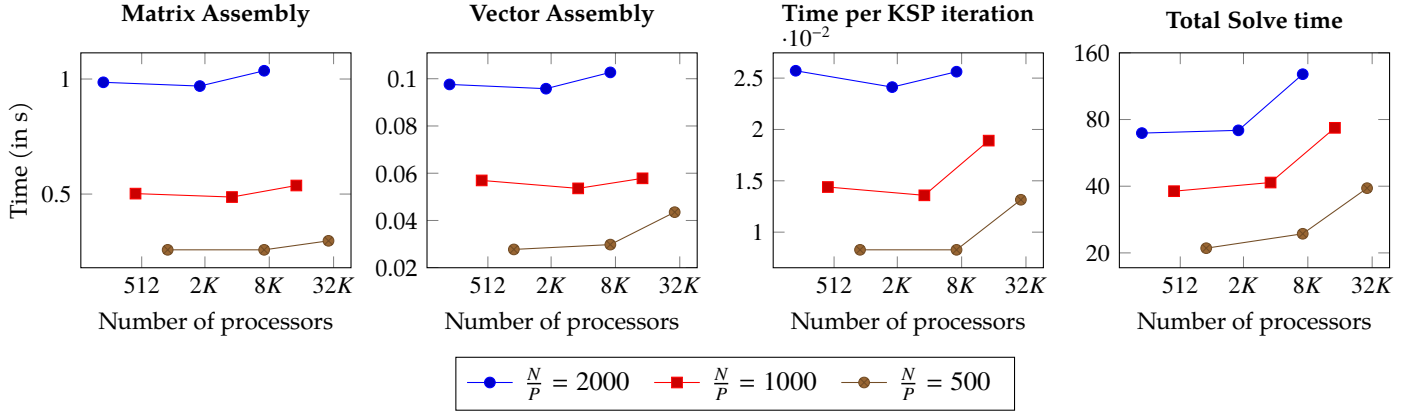
\begin{figure*}[h!]
\begin{tikzpicture}
  \begin{groupplot}[
      group style={group size=4 by 1,ylabels at=edge left,vertical sep=50pt},
      ylabel style={text height=0.02\textwidth,inner ysep=0pt},
      xlabel= Number of processors, 
      ylabel= Time (in s),
      width=0.28*\linewidth,
      /tikz/font=\small,
      xtick={512,2000,8000,32000},
      xticklabels={$512$,$2K$,$8K$,$32K$},
    ]
    \nextgroupplot[xmode=log,log basis x={2},ytick={0.1,0.5,1.0},yticklabels={$0.1$,$0.5$,$1$},title={\textbf{Matrix Assembly}}]%
     \addplot
         table[x expr={56*\thisrow{Nodes}},y expr={(\thisrow{MatVol}/\thisrow{MatTimes})},col sep=space]{Data/Scaling/WeakScaling1.dat};
        \addplot
         table[x expr={56*\thisrow{Nodes}},y expr={(\thisrow{MatVol}/\thisrow{MatTimes})},col sep=space]{Data/Scaling/WeakScaling2.dat};
       \addplot
         table[x expr={56*\thisrow{Nodes}},y expr={(\thisrow{MatVol}/\thisrow{MatTimes})},col sep=space]{Data/Scaling/WeakScaling3.dat};
    \coordinate (bot) at (rel axis cs:1,0);
    
    \nextgroupplot[xmode=log,log basis x={2},ytick={0.1,0.08,0.06,0.04,0.02},yticklabels={0.1,0.08,0.06,0.04,0.02},ymin=0.02,title={\textbf{Vector Assembly}}]%
    \addplot
         table[x expr={56*\thisrow{Nodes}},y expr={(\thisrow{VecVol}/\thisrow{VecTimes})},col sep=space]{Data/Scaling/WeakScaling1.dat};
        \addplot
         table[x expr={56*\thisrow{Nodes}},y expr={(\thisrow{VecVol}/\thisrow{VecTimes})},col sep=space]{Data/Scaling/WeakScaling2.dat};
       \addplot
         table[x expr={56*\thisrow{Nodes}},y expr={(\thisrow{VecVol}/\thisrow{VecTimes})},col sep=space]{Data/Scaling/WeakScaling3.dat};
    \coordinate (top) at (rel axis cs:0,1);
    
    \nextgroupplot[xmode=log,log basis x={2},
    title={\textbf{Time per KSP iteration}}]%
   \addplot
         table[x expr={56*\thisrow{Nodes}))},y expr={\thisrow{KSPSolve}/\thisrow{KSPIter}},col sep=space]{Data/Scaling/WeakScaling1.dat};
          \addplot
         table[x expr={56*\thisrow{Nodes}))},y expr={\thisrow{KSPSolve}/\thisrow{KSPIter}},col sep=space]{Data/Scaling/WeakScaling2.dat};
          \addplot
         table[x expr={56*\thisrow{Nodes}))},y expr={\thisrow{KSPSolve}/\thisrow{KSPIter}},col sep=space]{Data/Scaling/WeakScaling3.dat};
    \coordinate (top) at (rel axis cs:0,1);
        \nextgroupplot[xmode=log,log basis x={2},ymode=log,log basis x={10},
    ytick={20,40,80,160},yticklabels={$20$,$40$,$80$,160},
    ymax= 160,
    title={\textbf{Total Solve time}}]%
   \addplot
         table[x expr={56*\thisrow{Nodes}))},y expr={\thisrow{Solve}},col sep=space]{Data/Scaling/WeakScaling1.dat};
         \label{plot:NP2K}
          \addplot
         table[x expr={56*\thisrow{Nodes}))},y expr={\thisrow{Solve}},col sep=space]{Data/Scaling/WeakScaling2.dat};
         \label{plot:NP1K}
          \addplot
         table[x expr={56*\thisrow{Nodes}))},y expr={\thisrow{Solve}},col sep=space]{Data/Scaling/WeakScaling3.dat};
         \label{plot:NP500}
    \coordinate (top) at (rel axis cs:0,1);
  \end{groupplot}
  \path (top|-current bounding box.south)--
        coordinate(legendpos)
        (bot|-current bounding box.south);
  \matrix[
      matrix of nodes,
      anchor=north,
      draw,
      inner sep=0.2em,
    ]at([yshift=-1ex]legendpos)
    { \ref{plot:NP2K}& $\frac{N}{P} = 2000$ &[5pt]
      \ref{plot:NP1K}& $\frac{N}{P} = 1000$ &[5pt]
      \ref{plot:NP500}& $\frac{N}{P} = 500$ &[5pt]\\
      };
\end{tikzpicture}
\caption{Weak scaling result on TACC \Frontera. We considered 3 different adaptive meshes with the number of elements per processor varying from 500 to 2000.}
\label{fig: WeakScaling}
\end{figure*}{}

        


%% file: Data/scalingDomain.tex
{\renewcommand{\arraystretch}{1.0} 
\begin{table}[b!]
\captionsetup{font=small}
\centering
{\begin{tabular}{|c|c|c|c|}
\hline
 & $R_{bkg}$ & $R_{wake}$ & $R_{bdy}$ \\ \hline 
\textsc{M1} & 5 & 7 & 8 \\ \hline
\textsc{M2} & 6 & 8 & 9 \\ \hline
\textsc{M3} & 7 & 9 & 10 \\ \hline
\textsc{M4} & 8 & 10 & 11 \\ \hline
\end{tabular}%
}
\caption{{The level of refinement in the various region of the domain for scaling studies.}}
\label{tab: scaling-domain}
\end{table}
}

%% file: Plots/FrontTruckStreamLine.tex
\begin{figure*}[t!]
    \centering
    \includegraphics[width=0.95\linewidth,clip,trim={1.0in 1.0in 1.0in 2.8in}]{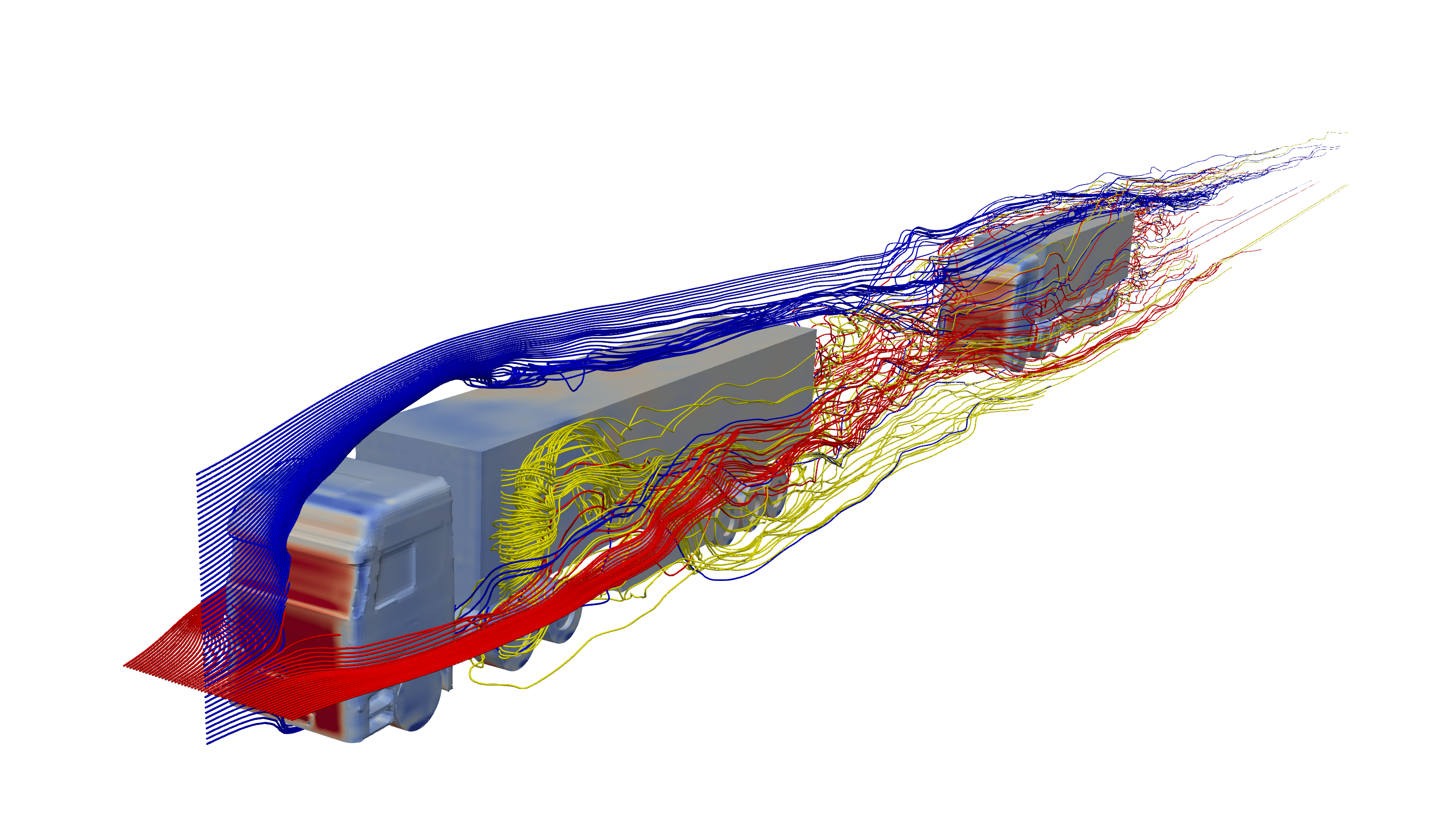}
    \captionof{figure}{Flow streamlines demonstrating the platooning effect with two trucks travelling at 65 MPH, corresponding to a Reynolds number, $Re$ of $30 \times 10^6$.}
    \label{fig:front-pageStreamlines}
\end{figure*}


%% file: Plots/TruckDragSingle.tex

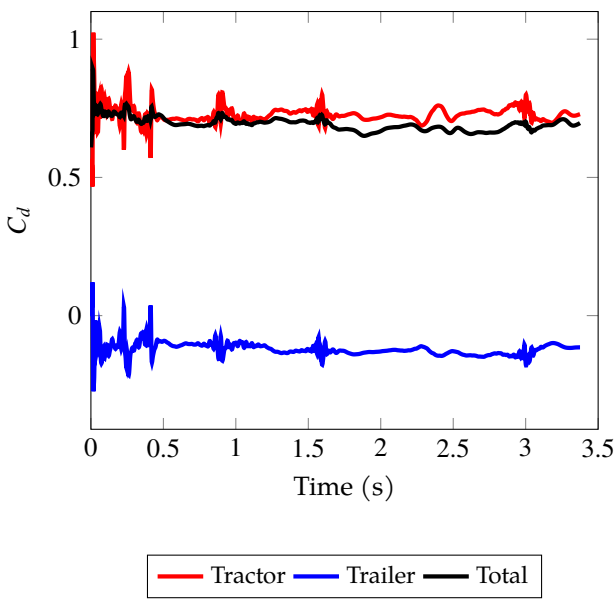
\begin{figure}[b!]
    \centering
    \begin{tikzpicture}
\begin{axis}[
          width=0.45\textwidth,
          xlabel=Time (s), 
          ylabel=$C_d$,
          xmin=0,
          xmax=3.5,
          ymax = 1.1,
          legend style={at={(1.2,0.7)},anchor=north west,legend columns=1}, 
          x tick label style={rotate=0,anchor=north}, 
          legend style={at={(0.5,-0.3)},anchor= north,legend columns=4}, 
]
    \addplot [no marks] gnuplot [raw gnuplot, color=red,ultra thick] {
      set datafile separator ",";
      plot "Data/SingleTruck/processed-Force_truck-head-new.dat" using
      ($1-400):($2*2/0.04373247949) with lines
     };
     
    \addplot [no marks] gnuplot [raw gnuplot, color=blue,ultra thick] {
      set datafile separator " ";
      plot "Data/SingleTruck/processed-Force_load.dat" using
      ($1-400):($2*2/0.04373247949) with lines
     };
     \addplot [no marks] gnuplot [raw gnuplot, color=black,ultra thick] {
      set datafile separator ",";
      plot "Data/SingleTruck/total.csv" using
      ($3-400):($2*2/0.04373247949) with lines
     };
     \legend{\small Tractor ,\small Trailer,\small Total}
  \end{axis}
  \end{tikzpicture}
  \caption{Time evolution of drag on tractor, trailer and total drag on the semi - truck.}
  \label{fig: TruckDragSingle}
\end{figure}

%% file: Plots/TwoTruckCp.tex
\begin{figure*}[t!]
\vspace{-1em}
\begin{subfigure}{0.33\textwidth}
\includegraphics[width=1.0\linewidth,clip, trim=0 0 100 200]{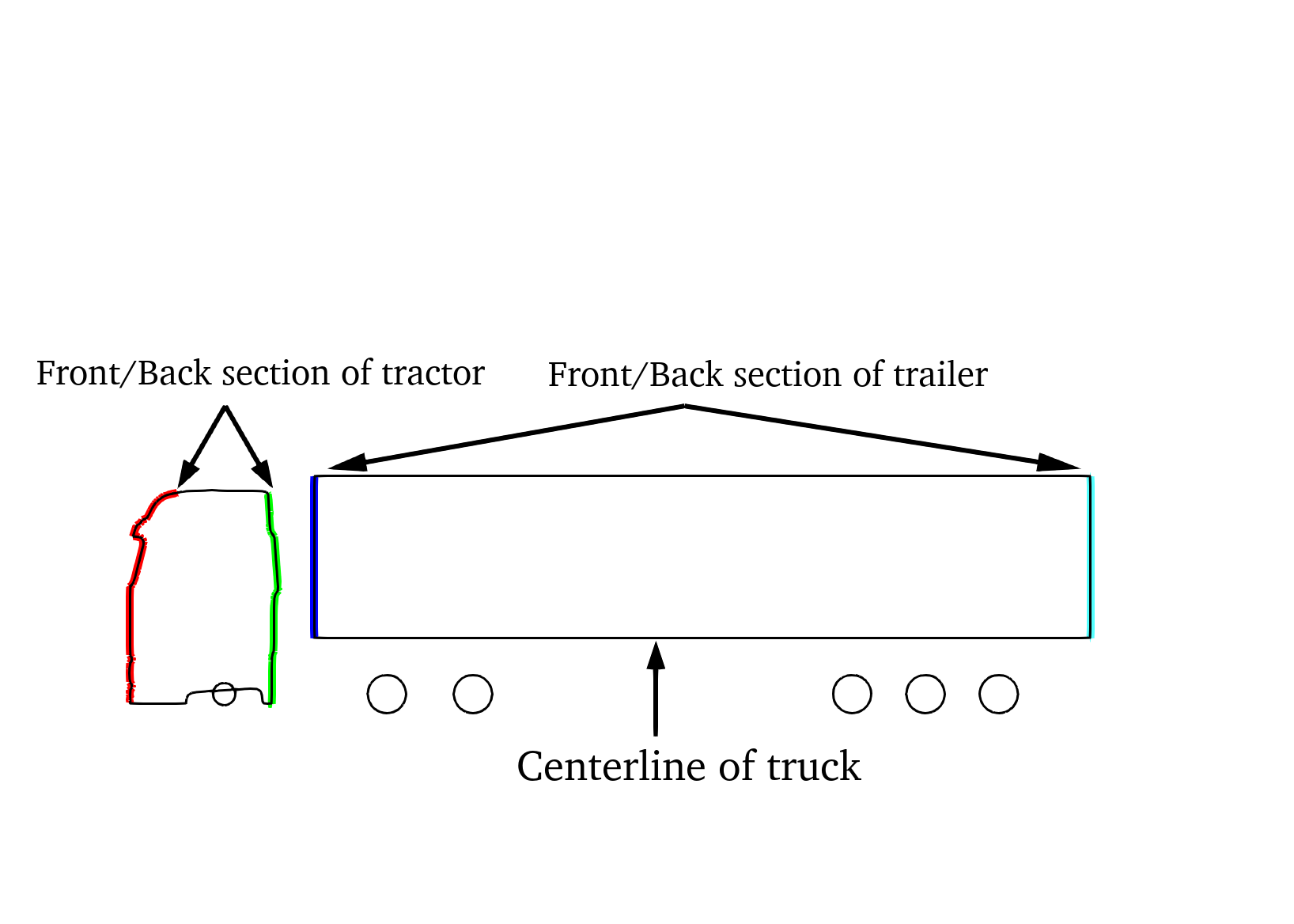}
\begin{minipage}{0.5in}\vfill\end{minipage}
\caption{Pressure probe points}
\label{fig: pprobe}
\hspace{-5mm}
\end{subfigure}
\begin{subfigure}{0.32\textwidth}
\centering
\begin{tikzpicture}
      \begin{axis}[
          width=1.0\linewidth, 
          legend columns=2, 
          xlabel=$C_p$, 
          ylabel=$y$,
          ymax = 0.3,
          ymin = 0.0,
          legend style={at={(0.5,-0.5)},anchor= north,legend columns=3}, 
          x tick label style={rotate=0,anchor=north}, 
          legend image post style={scale=3.0},
        ]

        \addplot[color = red, ultra thick] 
        table[x expr={2*\thisrow{p_ave}},y expr=\thisrow{y},col sep=comma]{Data/PlatooningTruck/CpSorted/head_front_1st.csv};
        
        \addplot[color=blue, ultra thick] 
        table[x expr={2*\thisrow{p_ave}},y expr=\thisrow{y},col sep=comma]{Data/PlatooningTruck/CpSorted/head_front_2nd.csv};
        
        \addplot[dashed,color=red,ultra thick] 
        table[x expr={2*\thisrow{p_ave}},y expr=\thisrow{y},col sep=comma]{Data/PlatooningTruck/CpSorted/head_back_1st.csv};
        
        \addplot[dashed,color=blue,ultra thick] 
        table[x expr={2*\thisrow{p_ave}},y expr=\thisrow{y},col sep=comma]{Data/PlatooningTruck/CpSorted/head_back_2nd.csv};

    \end{axis}
\end{tikzpicture}
\caption{Tractor}
\label{fig: TwoTruckPressureTractor}
\end{subfigure}
\begin{subfigure}{0.32\textwidth}
\centering
\begin{tikzpicture}
      \begin{axis}[
          width=1.0\linewidth, 
          legend columns=2, 
          xlabel=$C_p$, 
          ylabel=$y$,
          ymax = 0.3,
          ymin = 0.05,
          xmax = 0.05,
          legend style={at={(0.5,-0.5)},anchor= north,legend columns=3}, 
          x tick label style={rotate=0,anchor=north}, 
          legend image post style={scale=2.5},
        ]
    
         \addplot[color = red, ultra thick]  
        table[x expr={2*\thisrow{p_ave}},y expr=\thisrow{y},col sep=comma]{Data/PlatooningTruck/CpSorted/load_front_1st.csv};\label{fv1}
        
       \addplot[color = blue, ultra thick]  
        table[x expr={2*\thisrow{p_ave}},y expr=\thisrow{y},col sep=comma]{Data/PlatooningTruck/CpSorted/load_front_2nd.csv};\label{fv2}
        
       \addplot[dashed,color = red, ultra thick]  
        table[x expr={2*\thisrow{p_ave}},y expr=\thisrow{y},col sep=comma]{Data/PlatooningTruck/CpSorted/load_back_1st.csv};\label{bv1}
        
        \addplot[dashed,color = blue, ultra thick]
        table[x expr={2*\thisrow{p_ave}},y expr=\thisrow{y},col sep=comma]{Data/PlatooningTruck/CpSorted/load_back_2nd.csv};\label{bv2}
        
                
        \end{axis}
\end{tikzpicture}
\caption{Trailer}
\label{fig: TwoTruckPressureTrailer}
\end{subfigure}
\vspace{2 mm}
\caption{{\textit{Platooning effect:} Time-averaged non-dimensional pressure coefficient ($C_p)$ at the center-line of the truck; Front of vehicles 1 and 2 are represented by thick red and blue line and back is represented by dashed red and blue lines. \figref{fig: pprobe} shows the pressure probe points.}}
\label{fig: TwoTruckPressure}
\end{figure*}
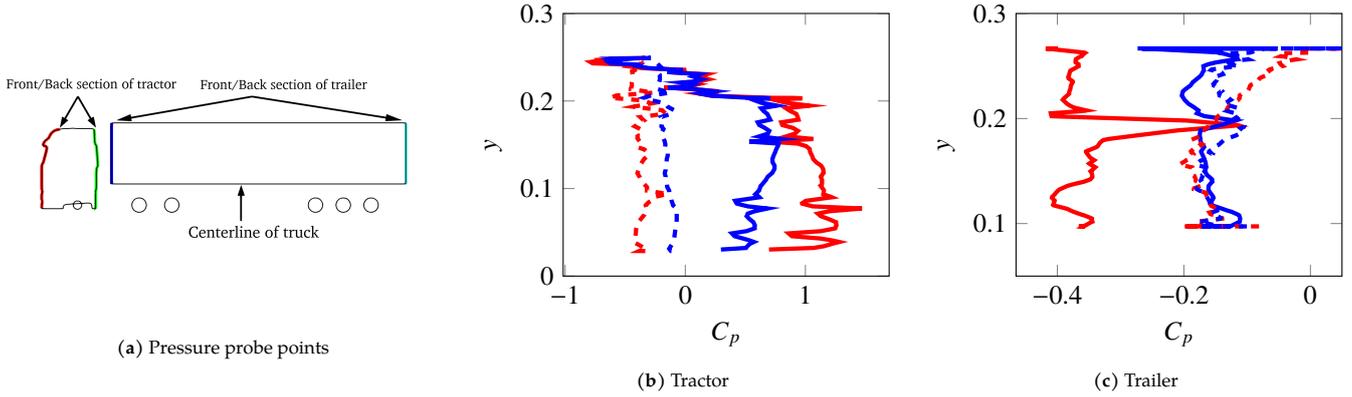

%% file: Plots/TwoTruckDrag.tex

\begin{figure}[t!]
 \centering
    \begin{tikzpicture}
\begin{axis}[
          width=0.45\textwidth,
          xlabel=Time (s), 
          ylabel=$C_d$,
          ymax = 1.0,
          xmax = 12.0,
          xmin=0,
          legend style={at={(1.2,0.7)},anchor=north west,legend columns=1}, 
          x tick label style={rotate=0,anchor=north}, 
          legend style={at={(0.5,-0.3)},anchor= north,legend columns=4}, 
			]
    \addplot [no marks] gnuplot [raw gnuplot, color=red, ultra thick] {
      set datafile separator ",";
      plot "Data/PlatooningTruck/processed-Force_truck-head-new.dat" using
      ($1-368):($2*2/0.04373247949) with lines
     };
     
    \addplot [no marks] gnuplot [raw gnuplot, color=blue, ultra thick] {
      set datafile separator ",";
      plot "Data/PlatooningTruck/processed-Force_load.dat" using
      ($1-368):($2*2/0.04373247949) with lines
     };
     
      \addplot [no marks] gnuplot [raw gnuplot, color=cpu4,ultra thick] {
      set datafile separator ",";
      plot "Data/PlatooningTruck/processed-Force_truck2-truck-head-new.dat" using ($1-368):($2*2/0.04373247949) with lines
     };
     
    \addplot [no marks] gnuplot [raw gnuplot, color=black,ultra thick] {
      set datafile separator ",";
      plot "Data/PlatooningTruck/processed-Force_truck2-load.dat" using ($1-368):($2*2/0.04373247949) with lines
     };
     \legend{\small Tractor1 ,\small Trailer1,\small Tractor2, \small Trailer2}
  \end{axis}
  \end{tikzpicture}
  \caption{Time evolution of drag on tractor and trailer demonstrating the platooning effect. The drag on second truck is significantly  lower than the first one.}
  \label{fig: TwoTruckDrag}
\end{figure}
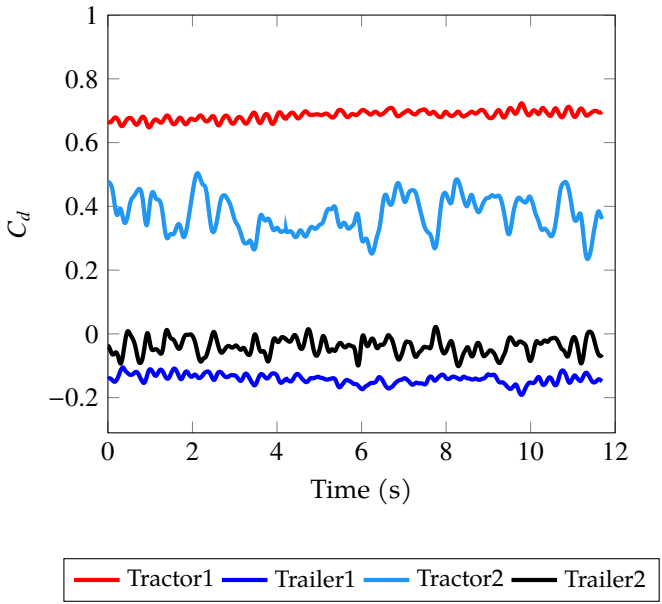{}

%% file: Conclusion.tex
We present a highly scalable, adaptive IMGA framework to solve industrial scale LES problem. We highlighted some of the key algorithmic challenges and improvements over the existing state-of-the-art IMGA methods. We have demonstrated excellent scaling results for our framework on current supercomputers, and shown preliminary application to practical problems. We believe that this will serve as a step towards achieving the goal of conducting overnight large scale LES simulations. 
\added{
Some future direction that can be explored in the context of IMGA are:
\begin{itemize}
    \item Carefully designed PDE and method specific preconditioners \cite{esmaily2015bi,de2017condition,jomo2019robust} that can exploit the matrix-free method and accelerate the convergence.
    \item Extension to higher-order finite element spaces. Higher order methods are more difficult to converge due to poor conditioning of matrices ~\cite{gahalaut2014condition}, and therefore careful design of the preconditioners is of utmost importance to achieve faster solve time.
    \item Designing fast(er) algorithms for assembly of FEM kernels (both matrix and matrix-free) that can achieve the lower bound on complexity.
    \item Efficient integration schemes for the \Intercepted elements.
    \item Scalable Multigrid methods~\citep{saberi2020parallel,jomo2020hierarchical,de2020multigrid} in the context of IMGA.
\end{itemize}
}

%% file: ae.tex
\appendix
\section{Precomputation of operators for FEM}\label{sec: MatPrecompute}
\added{A careful analysis of the FEM kernels resulting from the weakening of the Navier--Stokes equations reveals that the evaluation of these four basis function values (${\partial \phi \over \partial x};{\partial \phi \over \partial y}; {\partial \phi \over \partial z};\phi$) at the Gauss points is repeatedly performed. Thus, the values corresponding to them, over a reference element at the Gauss quadrature points can be pre-computed and cached in a matrix form. These matrices can be explicitly represented as:
\begin{equation*}
\begin{split}
    \Tensor{ [ \mathbf{\nabla \phi ]}}^k & = \bigg[{\partial \phi \over \partial k}\bigg]_{ij}  \quad  \mathrm{for} \quad i \in 1 \cdots \mathrm{nbf}; j  \in 1 \cdots \mathrm{ngp}; k \in (x,y,z) \\
    \Tensor{\mathbf{\phi}} & = [\phi]_{ij}  \quad  \mathrm{for} \quad i \in 1 \cdots \mathrm{nbf}; j  \in 1 \cdots \mathrm{ngp}; 
\end{split}
\end{equation*}
where nbf donates the number of basis functions, ngp donates the number of Gauss points, $[{\partial \phi \over \partial k}]_{ij}$ denotes the value of the  $k^{th}$ direction derivative of the $i^{th}$ basis function  at $j^{th}$ Gauss point. Similarly $\phi_{ij}$ denotes the value of the $i^{th}$ basis function at  Gauss point $j$. Using standard Gauss - Legendre quadrature rule, each of these forms are matrices of size $\mathrm{(nbf+1)}^{d} \times \mathrm{(nbf+1)}^{d}$, where $d$ is the spatial dimension ($d=3$). Note that the derivative along each direction $k$ is stored separately as a matrix.

Using these precomputed forms, the various contributing expressions can be efficiently computed. For example, the stiffness matrix, can be computed as:
\begin{equation*}
\begin{split}
    \Tensor{K} = &\textsc{compute\_ele\_matrix}(\textsc{Mat\_op1 =} \Tensor{ [ \mathbf{\nabla \phi ]}}^x, \textsc{Mat\_op2 =} \Tensor{ [ \mathbf{\nabla \phi ]}}^x) + \\
    & \textsc{compute\_ele\_matrix}(\textsc{Mat\_op1 =} \Tensor{ [ \mathbf{\nabla \phi ]}}^y, \textsc{Mat\_op2 =}\Tensor{ [ \mathbf{\nabla \phi ]}}^y) + \\
    & \textsc{compute\_ele\_matrix}(\textsc{Mat\_op1 =} \Tensor{ [ \mathbf{\nabla \phi ]}}^z, \textsc{Mat\_op2 =}\Tensor{ [ \mathbf{\nabla \phi ]}}^z)
\end{split}
\end{equation*}

Other operators for a complicated FEM kernel, for instance $(\mathbf{u} \cdot \nabla \phi_i,\mathbf{u} \cdot \nabla \phi_j)$,  can be evaluated as a two step process: first interpolating nodal values (here, $\mathbf{u}$) onto Gauss points (Complexity: $\mathcal{O}(d(\mathrm{bf+1})^{d+1})$ ~\cite{deville2002high}) and then performing the element-wise matrix scalar multiplication. (Complexity: $\mathcal{O}\mathrm{(bf+1)}^{2d}$). Thus, the overall complexity of matrix assembly remain bounded by the complexity of matrix-multiplication for all FEM kernels.
}
\section{Details of solver selection}
\petsc~ was used to solve all the linear algebra problems. In particular, bi-conjugate gradient descent (\texttt{-ksp\_type bcgs}) solver was used in conjunction with Additive - Schwartz (\texttt{-pc\_type asm}) preconditioner to solve the linear system of equations. The \textsc{NEWTONLS} class by \petsc, that implements a Newton Line Search method, was used for the nonlinear problems. Both the relative residual tolerance and the absolute residual tolerance for linear and non - linear solve are set to $10^{-6}$ in all numerical results.

The correctness of the code has been validated by solving the Navier--Stokes with a known analytical solution. The validation case for flow past a sphere performed in \ref{sec:sphere} gives us confidence in the correctness of the framework. The timings reported are measured through \petsc~ \texttt{log\_view} routine. Below we provide additional simulation details for the semi - truck simulations.

\section{Additional Simulation Details}

\subsection{Scaling studies with AMG preconditioners}\label{sec:AMG}
\figref{fig: AMG} shows the strong scaling result with AMG. We see that with increasing number of processors, the time for initial AMG setup increases. This result is consistent with what is observed previously in the simulations by ~\citet{sundar2012parallel}. This limits our ability to deploy AMG on a large number of processor. We are currently working on GMG to avoid this setup cost.

\input{Plots/scalingAMG}
\subsection{Truck Simulations}\label{sec:AEtruck-simulation}
\input{Plots/TruckDomain}
\input{Plots/TruckCAD}

\input{Plots/TruckStreamline}

\input{Plots/TruckTireComparison}
\input{Plots/TruckLICyplane}
\input{Plots/TwoTruckTopDownFlowStructure}
\figref{fig: TruckDomain} shows the computational domain of the truck and \figref{fig: TruckCAD} shows the   details  of  the  CAD model  of  the  truck.  \figref{fig: TruckStreamline} shows the streamlines starting from four different locations passing the truck. The blue streamline, starting from a vertical line in front of the tractor shows that the flow stagnates at $40\%$ height of the truck. The lower portion gets pushed towards the undercarriage of the vehicle and interacts with the rotating wheels. The green streamlines in \figref{fig: TruckSL3} shows the flow at the top of the tractor does not pass the top surface of the trailer smoothly, some portion of the flow is blocked by the extra height of the trailer and enters the gap between the tractor and the trailer, therefore creating extra turbulence downstream.\\
\figref{fig: TruckLICmidz} shows the vertical mid-plane flow structure around the truck. We observed flow re-circulation at the leading edge of the top surface of trailer followed by flow re-attachment downstream.
\figref{fig: TruckTireComparison} shows the affect of rotating wheel in the simulation. In case of stationary wheels, (\figref{fig: TruckTireStationary}), flow passes through the gap between tires and trailer without obstruction and no re-circulation near the ground between tires is observed, whereas in case of rotating wheel (\figref{fig: TruckTireRotate}), the clear vortex structures are seen near the tires.

\figref{fig: TruckLICy} shows the slices at different $y$ location at $0.05$, $0.1$, $0.15$ and $0.2$, it can be seen that the rotating tires causes flow separation at the leading edge of the side surface of the trailer, which maybe reduced by a side skirt device. Finally, in \figref{fig: TruckLICx}, we show the flow structure at different $x$ slices. This flow pattern resembles the vortices coming off the wing-tip of an airplane, which contributes to additional drag.

\figref{fig: TwoTruckTopDownFlowStructure} shows the comparison in flow structure of the leading truck and the trailing truck in top-down view. The platooning truck is in the turbulent wake of the leading truck. The asymmetric incoming flow triggers the separation at the side of the tractor. The flow structure shows reduction in the re-circulation region at both sides of the tractor and earlier reattachment onto the trailer surface.

%% file: Plots/scalingAMG.tex
\begin{figure}[h!]
	\centering
\begin{tikzpicture}
    \begin{loglogaxis}[
          width=0.95\linewidth, 
          ylabel=Time (in s), 
          xlabel=Number of Processor,
          legend style={at={(1.2,1.0)},anchor=north west,legend columns=1}, 
          x tick label style={rotate=0,anchor=north}, 
          log basis y={10},
          log basis x={10},
          legend style={at={(0.5,-0.3)},anchor= north,legend columns=5}, 
]
\addplot
     table[x expr={56*\thisrow{Nodes}},y expr={(\thisrow{AMG})},col sep=space]{Data/Scaling/ScalingAMG.txt};
     \addplot
     table[x expr={56*\thisrow{Nodes}},y expr={\thisrow{Time})},col sep=space]{Data/Scaling/ScalingAMG.txt};
    \legend{\small{AMG setup}, \small{Solve time}}
  \end{loglogaxis}
  \end{tikzpicture}
  \caption{{The time for AMG setup and complete solve time for mesh \textsc{M3} with increasing number of processor.}}
  \label{fig: AMG}
\end{figure}
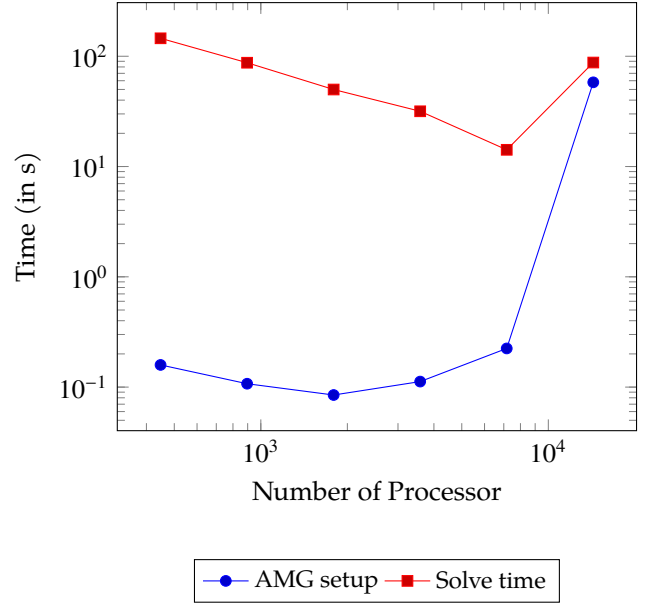{}

%% file: Plots/TruckDomain.tex
\begin{figure}[t!]
  \centering
  \includegraphics[width=0.9\linewidth]{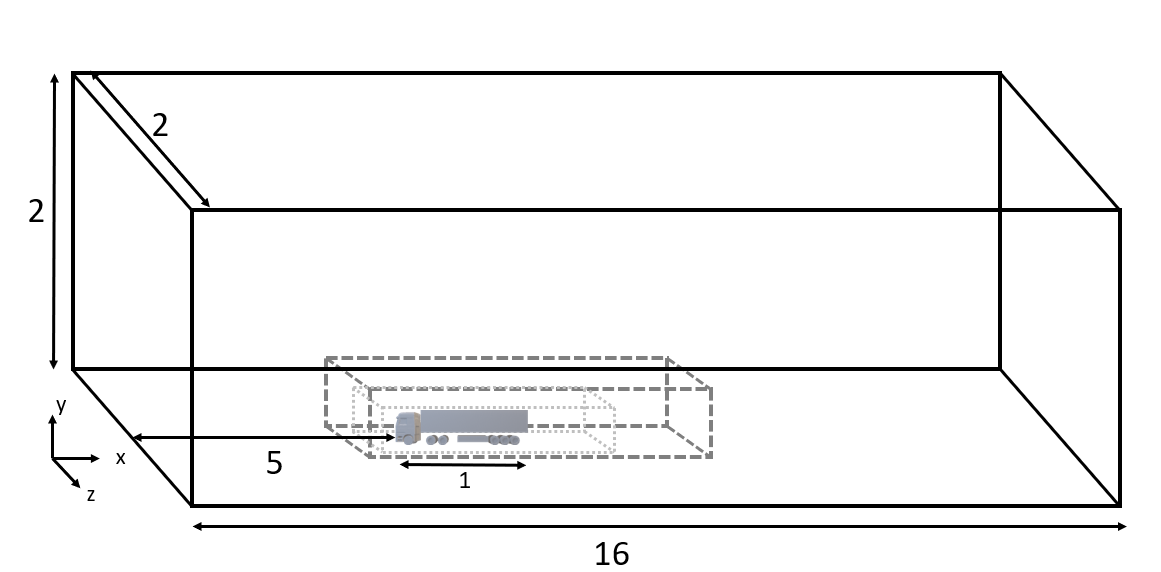}
  \caption{Computational domain used for the flow simulation over a semi-trailer truck.}
  \label{fig: TruckDomain}
\end{figure}

%% file: Plots/TruckCAD.tex
\begin{figure*}[t!]
  \includegraphics[width=1.0\linewidth]{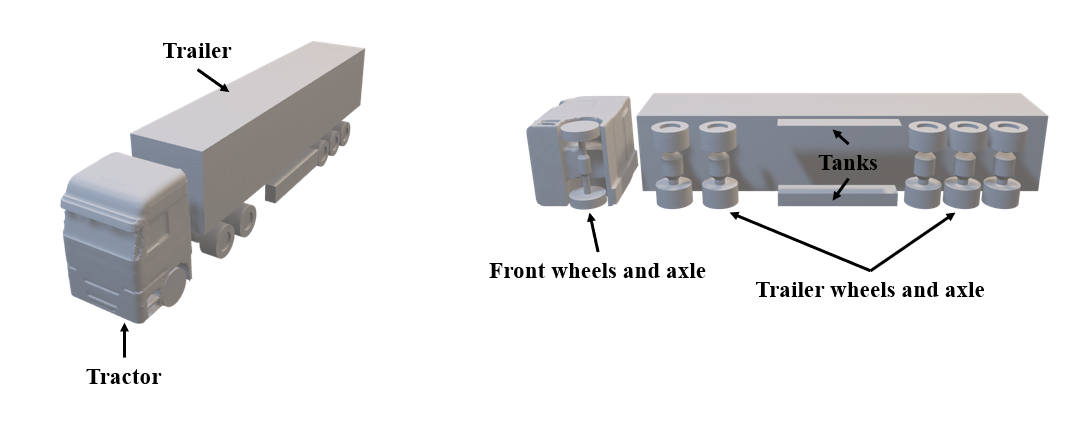}

\caption{Geometry and parts of the simulated semi-trailer truck.}
\label{fig: TruckCAD}
\end{figure*}

%% file: Plots/TruckStreamline.tex
\begin{figure*}[t!]
\centering
\begin{subfigure}{0.9\textwidth}
  \includegraphics[width=1.0\linewidth]{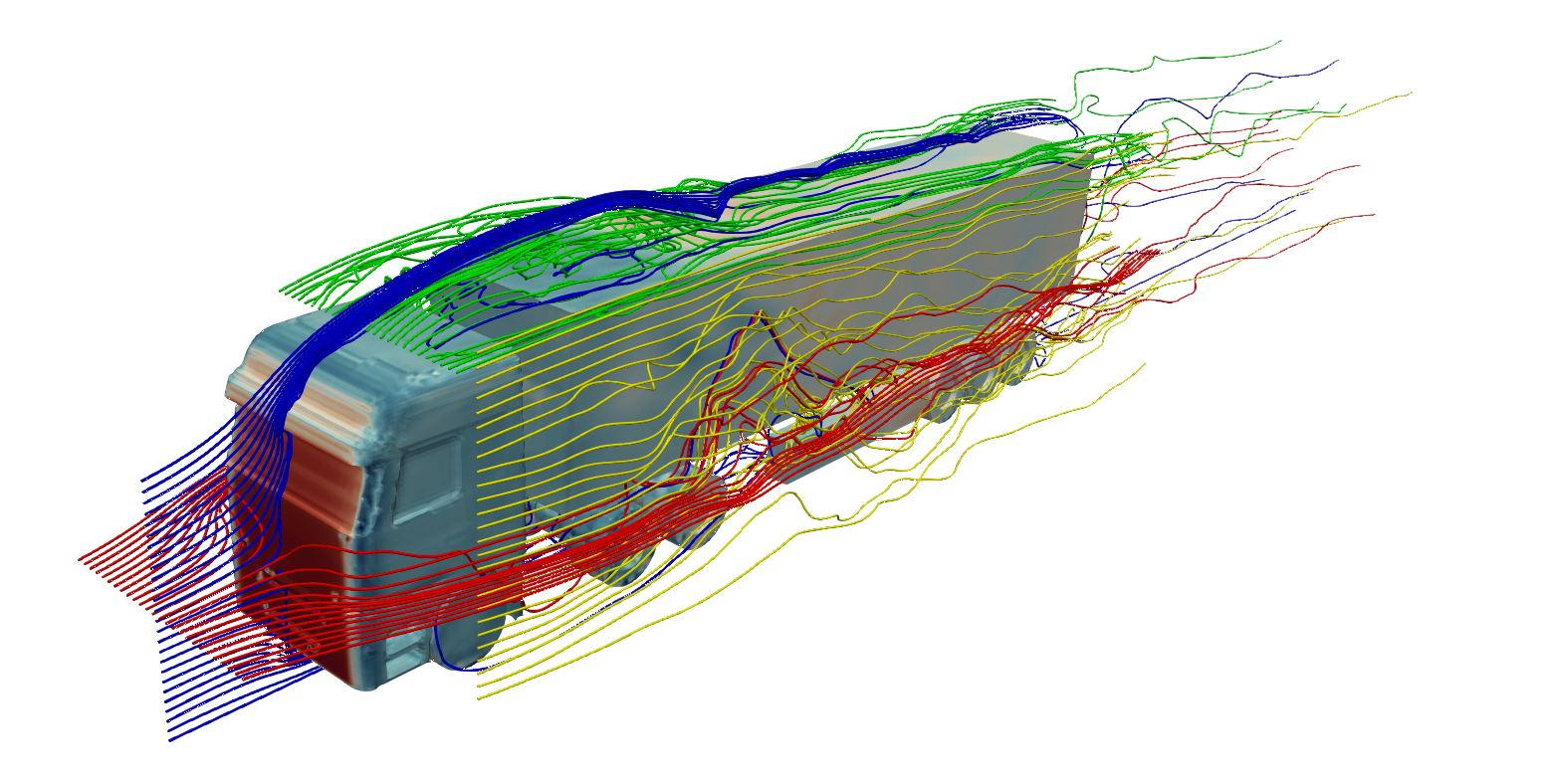}
  \caption{Overview of the flow streamlines passing the truck}
  \label{fig: TruckSL1}
\end{subfigure}
\\
\begin{subfigure}{.45\textwidth}
  \includegraphics[width=1.0\linewidth]{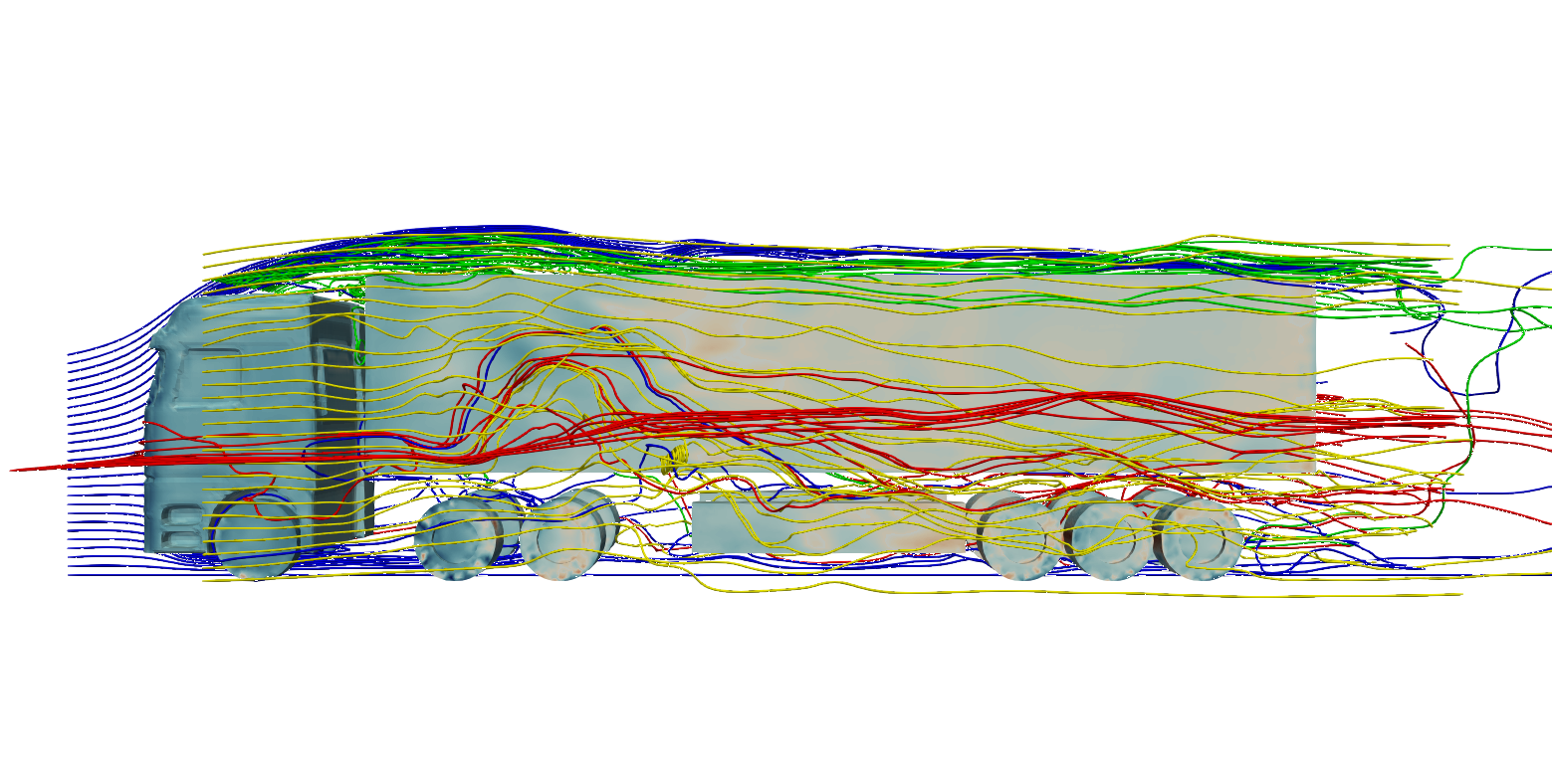}
  \caption{Side view of the flow streamlines passing the truck}
  \label{fig: TruckSL2}
\end{subfigure}
\begin{subfigure}{.45\textwidth}
  \includegraphics[width=1.0\linewidth]{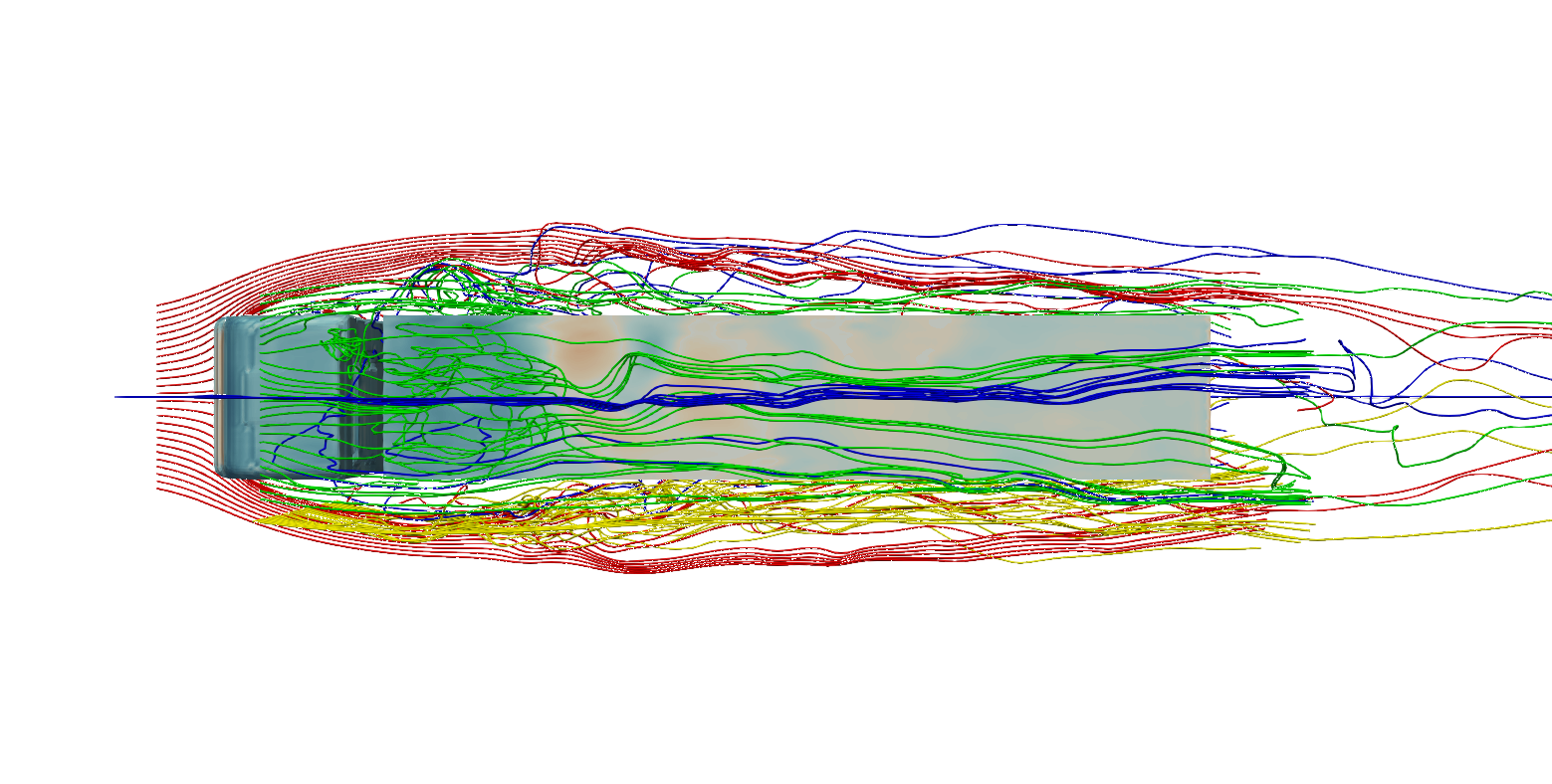}
  \caption{Top view of the flow streamlines passing the truck}
  \label{fig: TruckSL3}
\end{subfigure}
\vspace{0.1in}
\caption{Flow streamlines passing the truck at different locations.}
\label{fig: TruckStreamline}
\end{figure*}

%% file: Plots/TruckTireComparison.tex
\begin{figure*}[t!]
\centering
\begin{subfigure}{.45\textwidth}
  \includegraphics[width=1.0\linewidth]{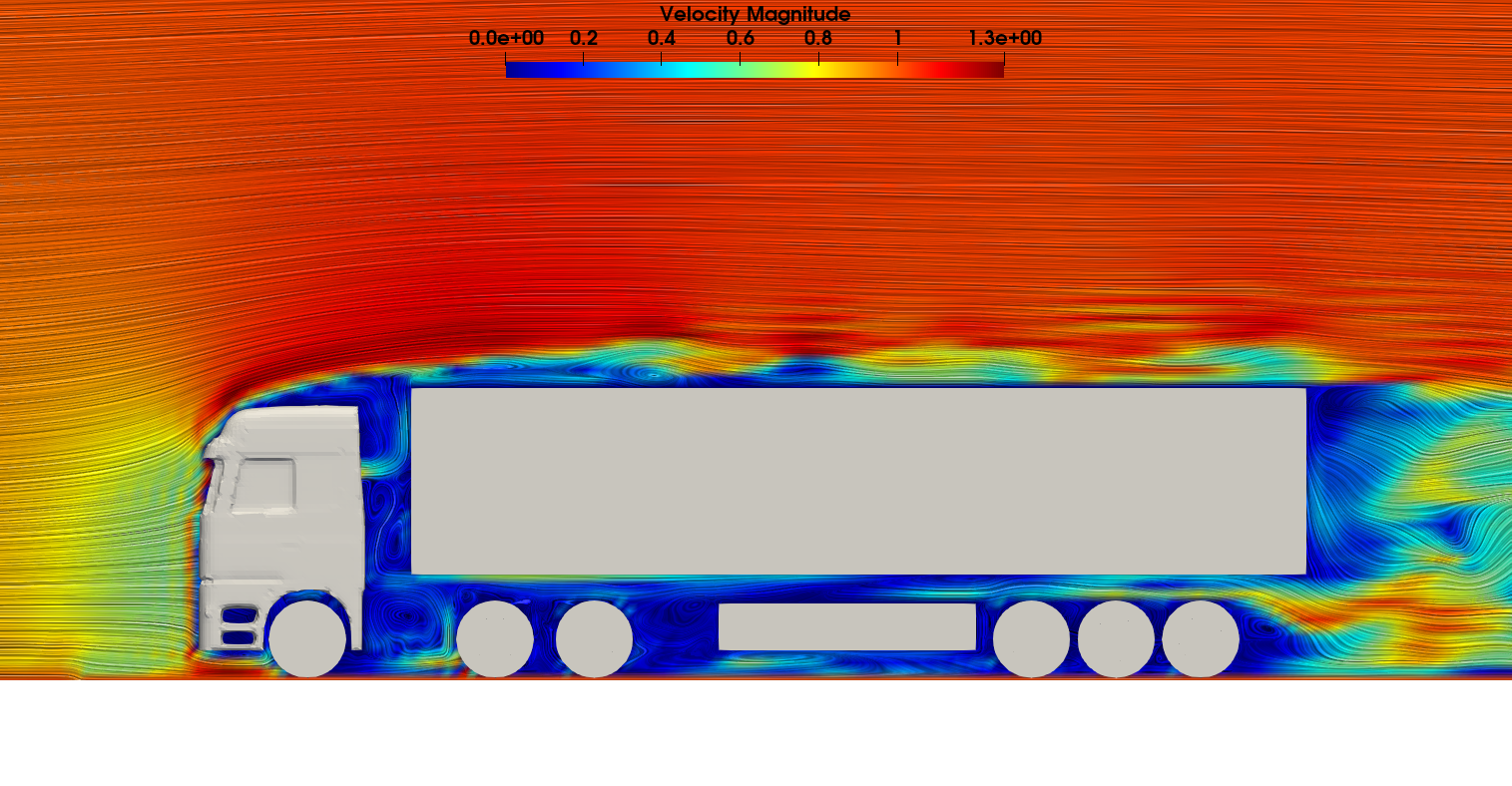}
  \caption{Snapshot of flow structures around a stationary wheels.}
  \label{fig: TruckTireStationary}
\end{subfigure}
\begin{subfigure}{.45\textwidth}
  \includegraphics[width=1.0\linewidth]{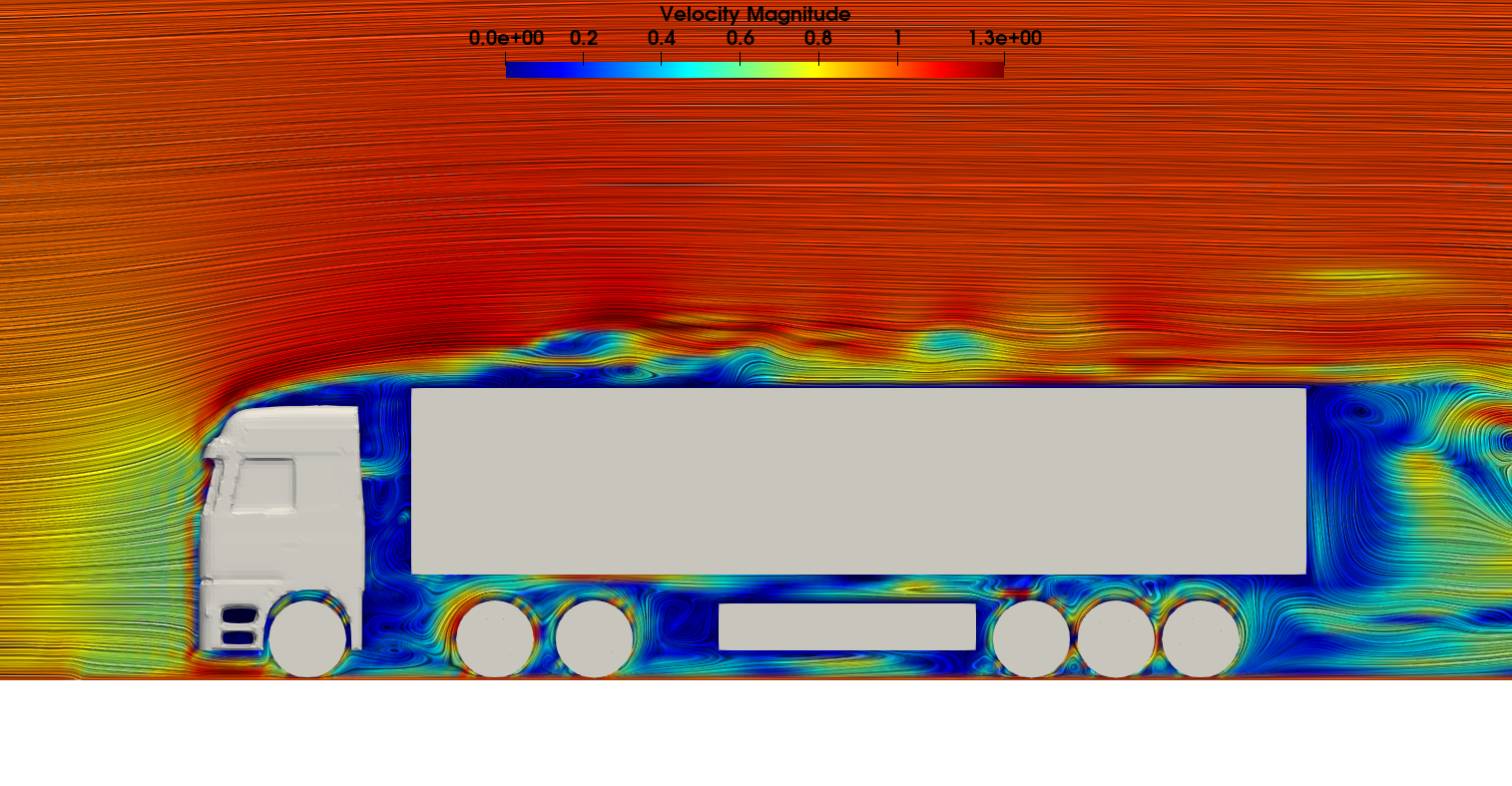}
  \caption{Snapshot of flow structures around a rotating wheels.}
  \label{fig: TruckTireRotate}
\end{subfigure}
\vspace{0.1in}
\caption{Comparison of flow structures of stationary vs. rotating wheels.}
\label{fig: TruckTireComparison}
\end{figure*}

%% file: Plots/TruckLICyplane.tex
\begin{figure*}[t!]
\centering
\begin{subfigure}{.48\textwidth}
  \includegraphics[width=1.0\linewidth]{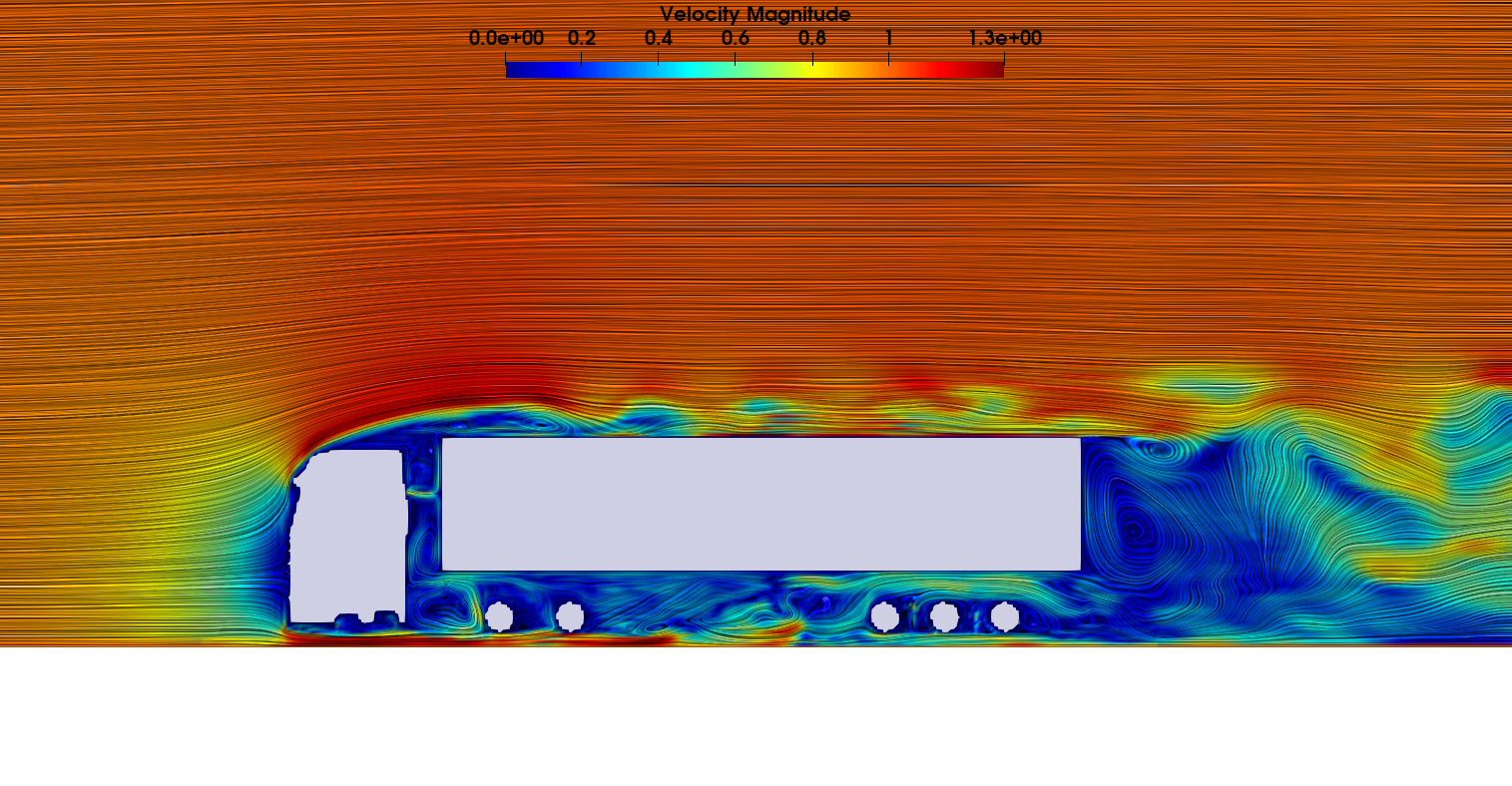}
  \caption{Flow structures at vertical middle plane.}
  \label{fig: TruckLICmidz}
\end{subfigure}
\begin{subfigure}{.48\textwidth}
  \includegraphics[width=1.0\linewidth]{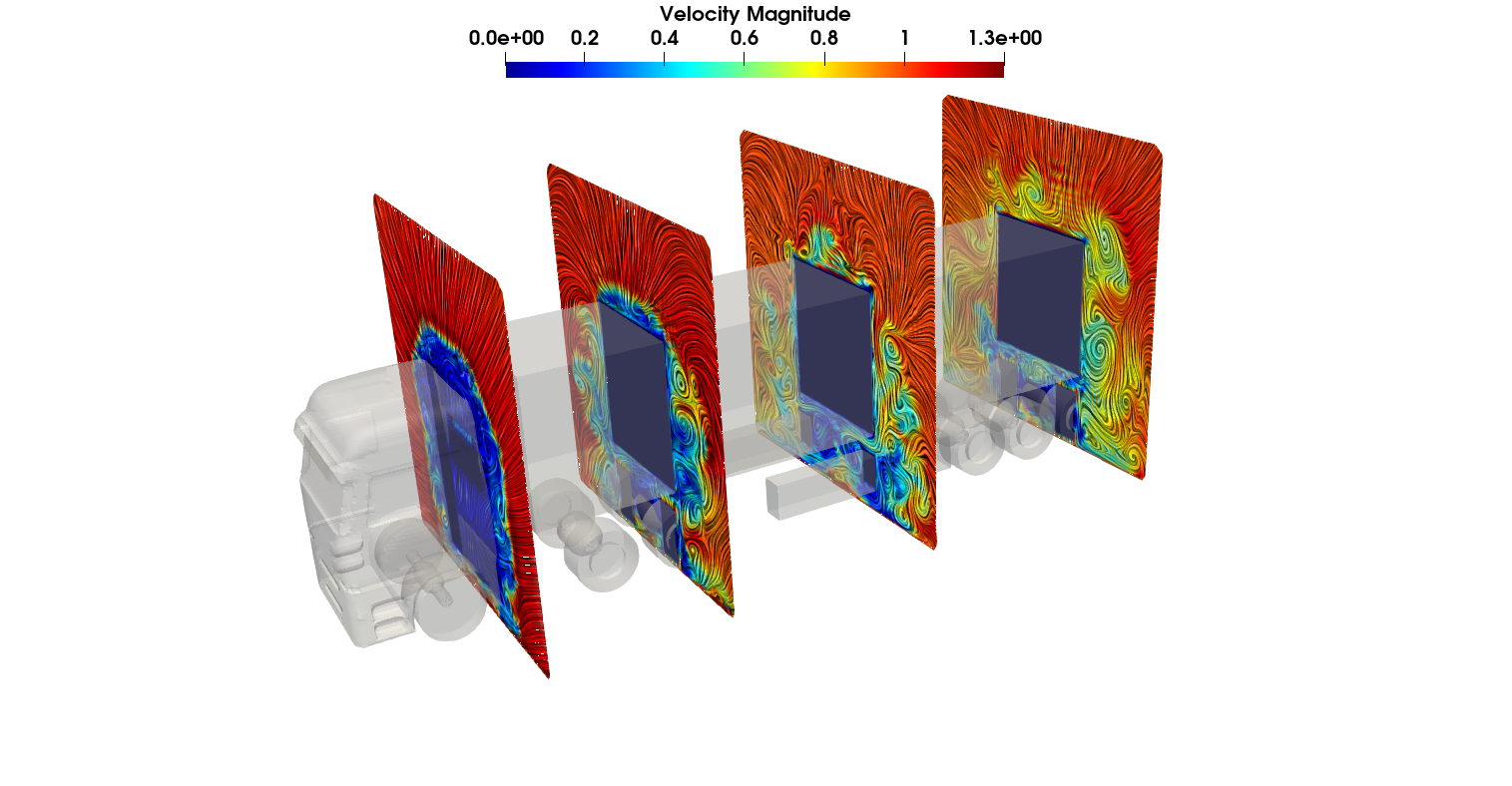}
  \caption{Flow structures at different slices in the $x$ direction.}
  \label{fig: TruckLICx}
\end{subfigure}\\
\vspace{0.1in}
\begin{subfigure}{.48\textwidth}
  \includegraphics[width=1.0\linewidth]{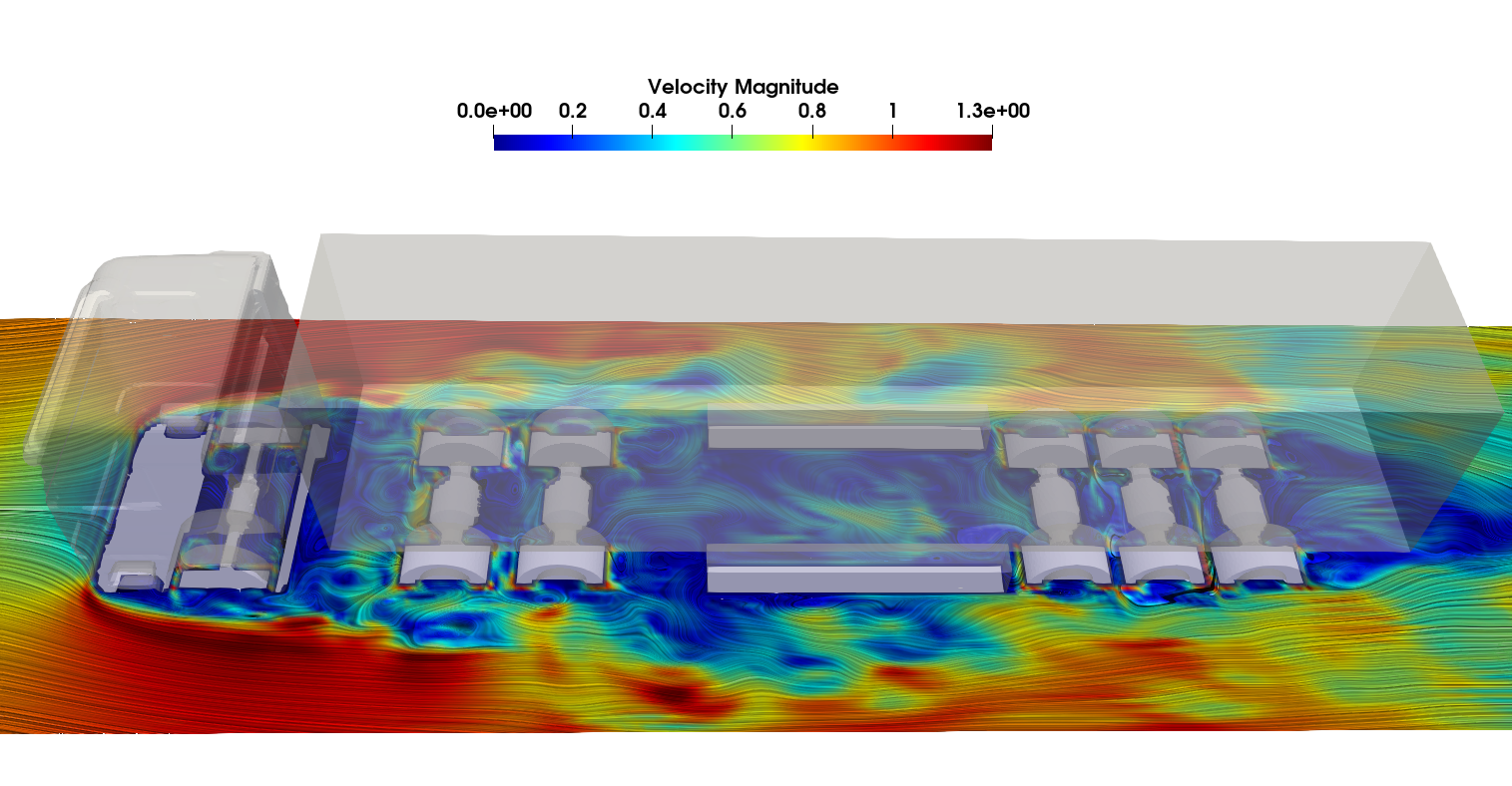}
  \caption{Flow structures near the truck at $y=0.05$}
  \label{fig: TruckLICy_1}
\end{subfigure}
\begin{subfigure}{.48\textwidth}
  \includegraphics[width=1.0\linewidth]{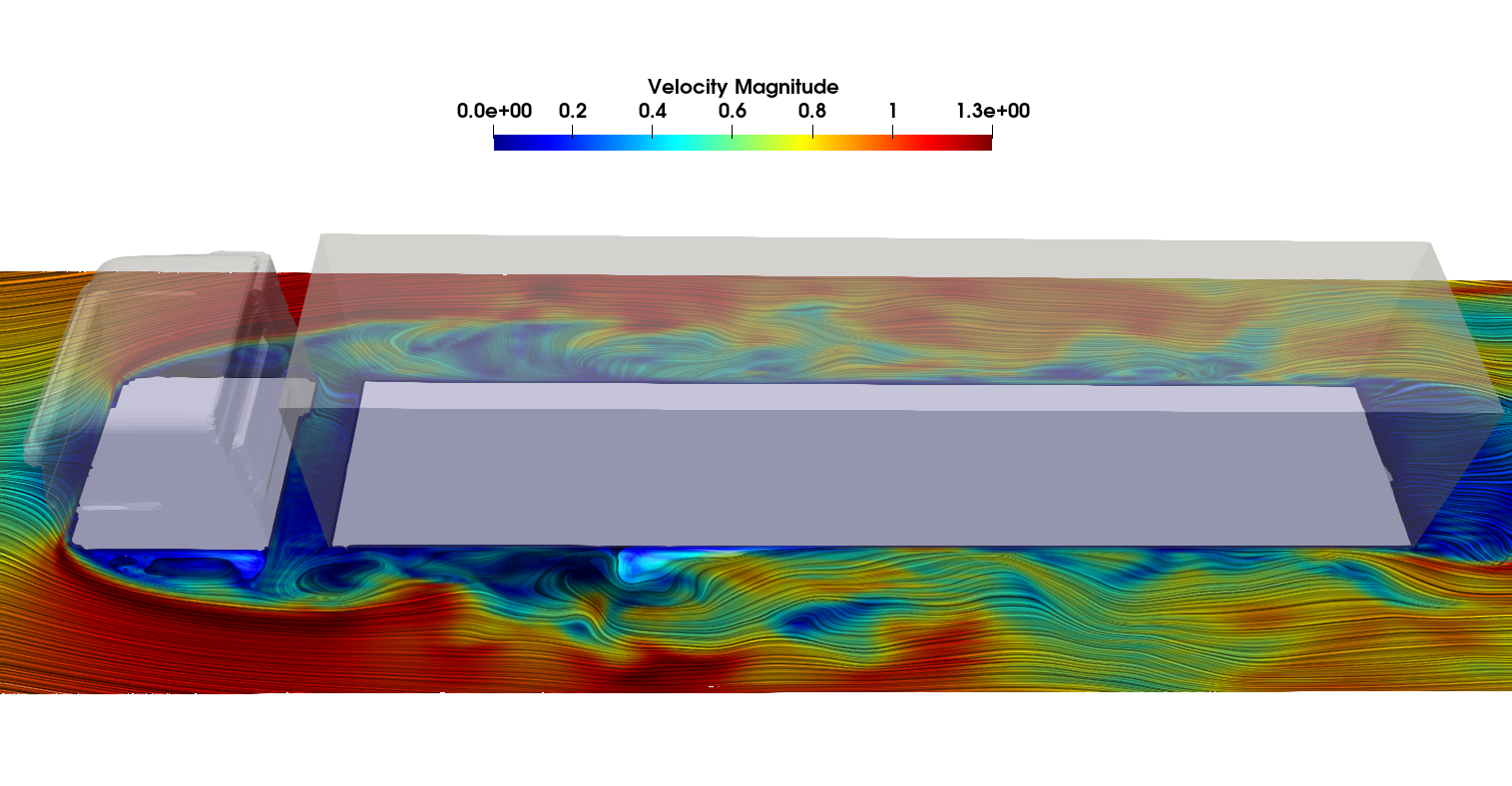}
  \caption{Flow structures near the truck at $y=0.10$}
  \label{fig: TruckLICy_2}
\end{subfigure}\\
\vspace{0.1in}
\begin{subfigure}{.48\textwidth}
  \includegraphics[width=1.0\linewidth]{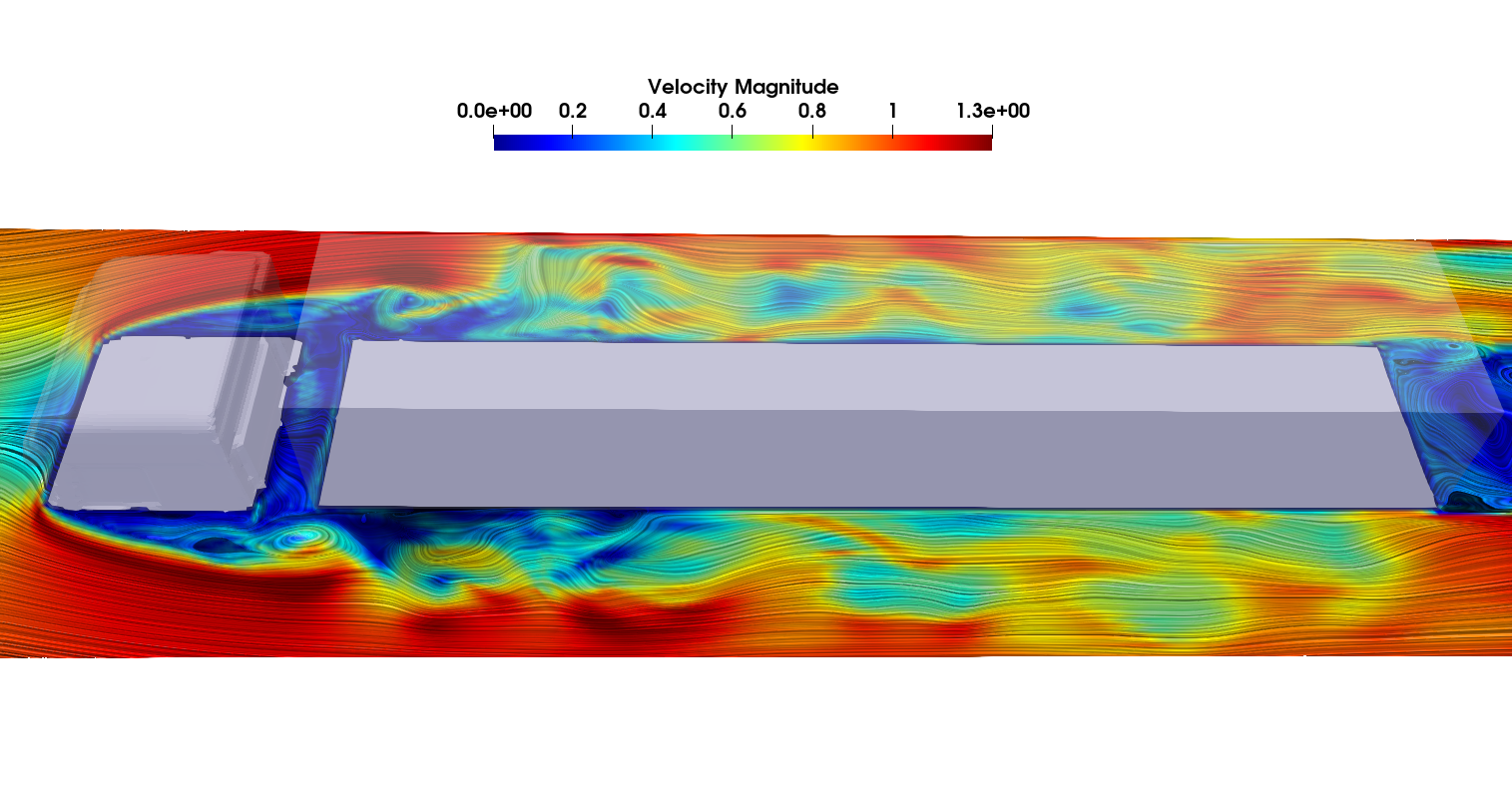}
  \caption{Flow structures near the truck at $y=0.15$}
  \label{fig: TruckLICy_3}
\end{subfigure}
\begin{subfigure}{.48\textwidth}
  \includegraphics[width=1.0\linewidth]{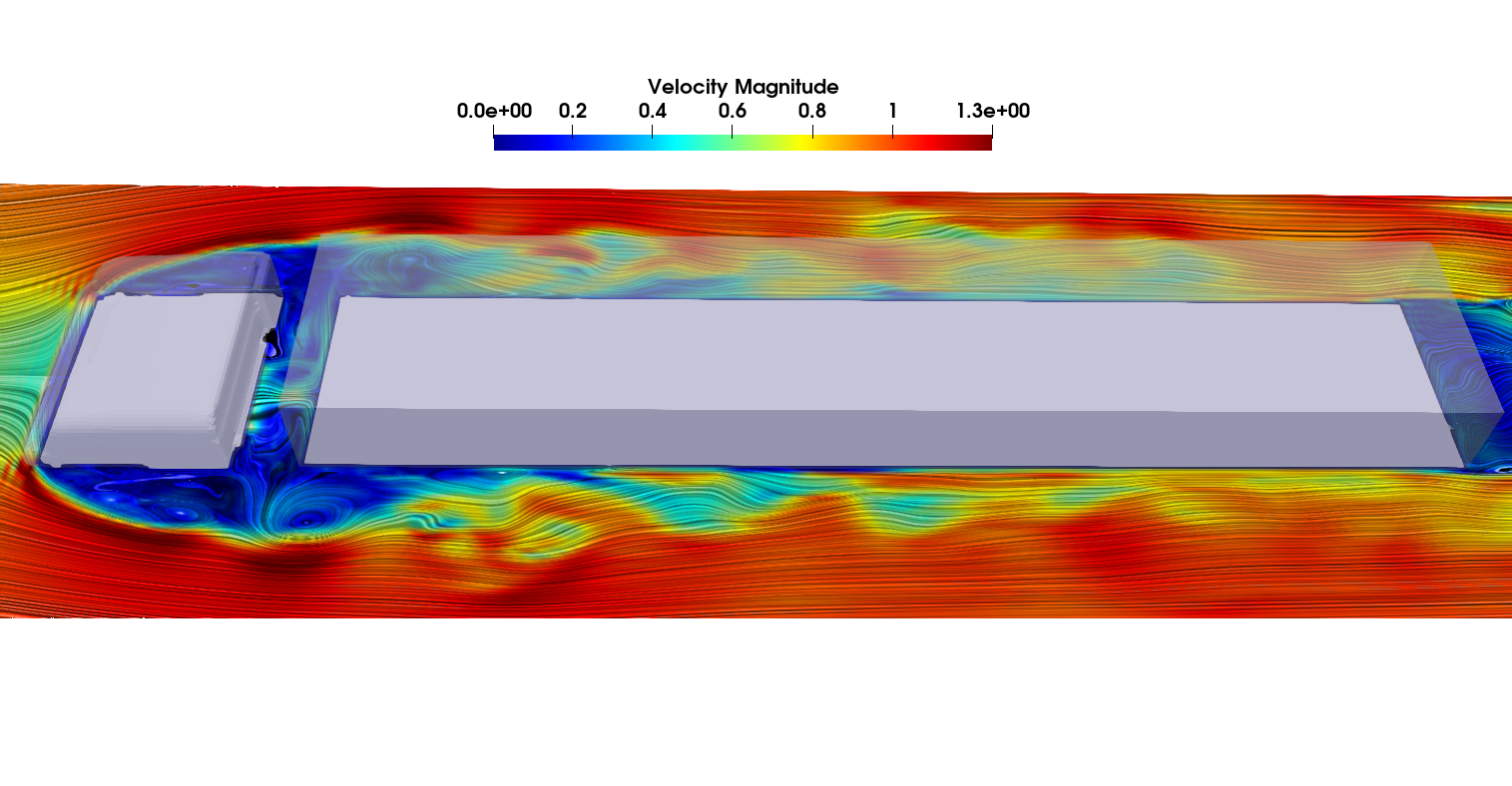}
  \caption{Flow structures near the truck at $y=0.20$}
  \label{fig: TruckLICy_4}
\end{subfigure}
\vspace{0.1in}
\caption{Flow structures around the truck at different locations.}
\label{fig: TruckLICy}
\end{figure*}

%% file: Plots/TwoTruckTopDownFlowStructure.tex
\begin{figure*}[t!]
\centering
\begin{subfigure}{.48\textwidth}
  \includegraphics[width=0.9\linewidth]{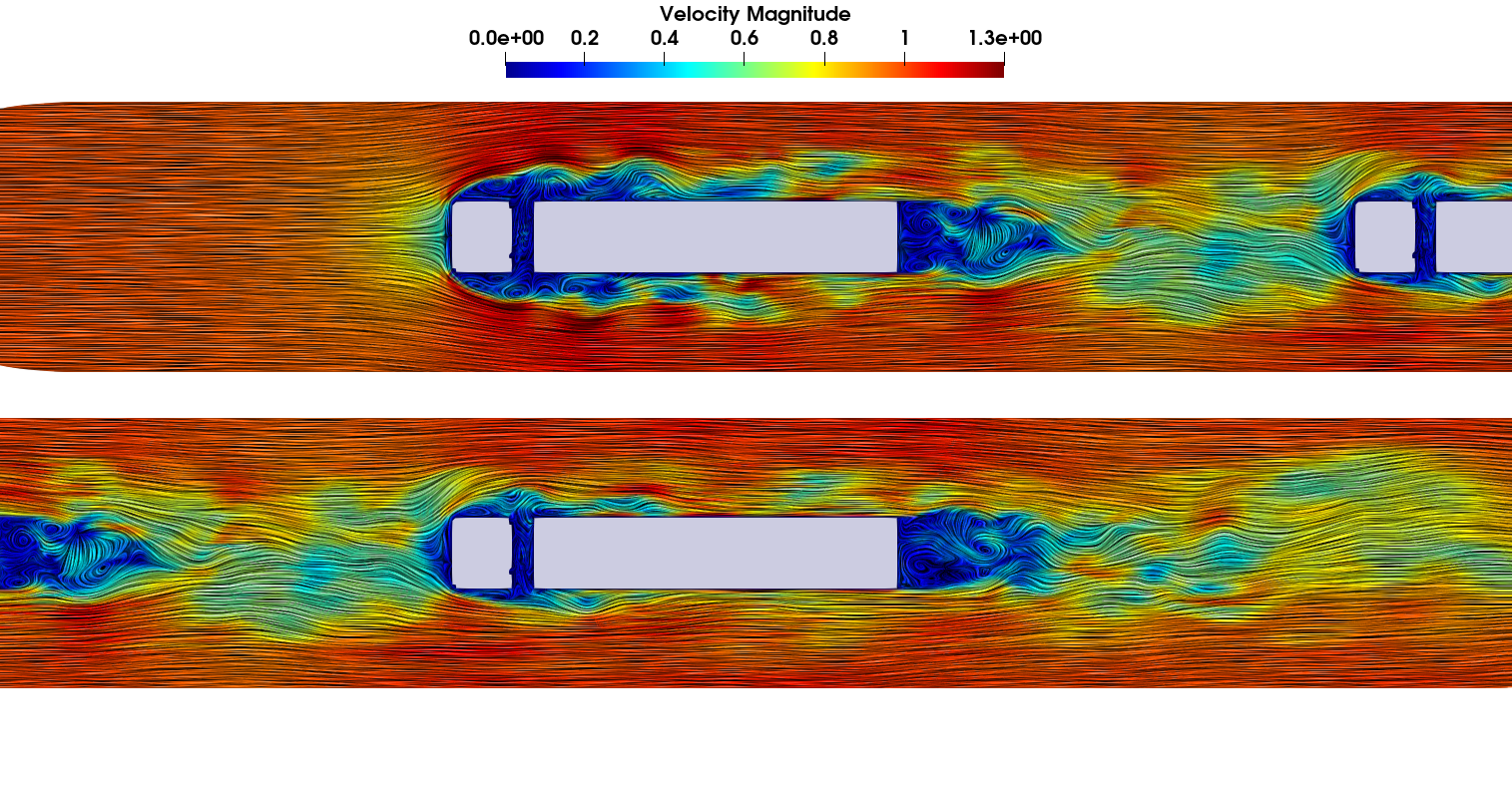}
  \vspace{-5 mm}
  \caption{Flow structures at $T = 0$}
  \label{fig: TwoTruckFlowStructure_1}
\end{subfigure}
\begin{subfigure}{.48\textwidth}
  \includegraphics[width=0.9\linewidth]{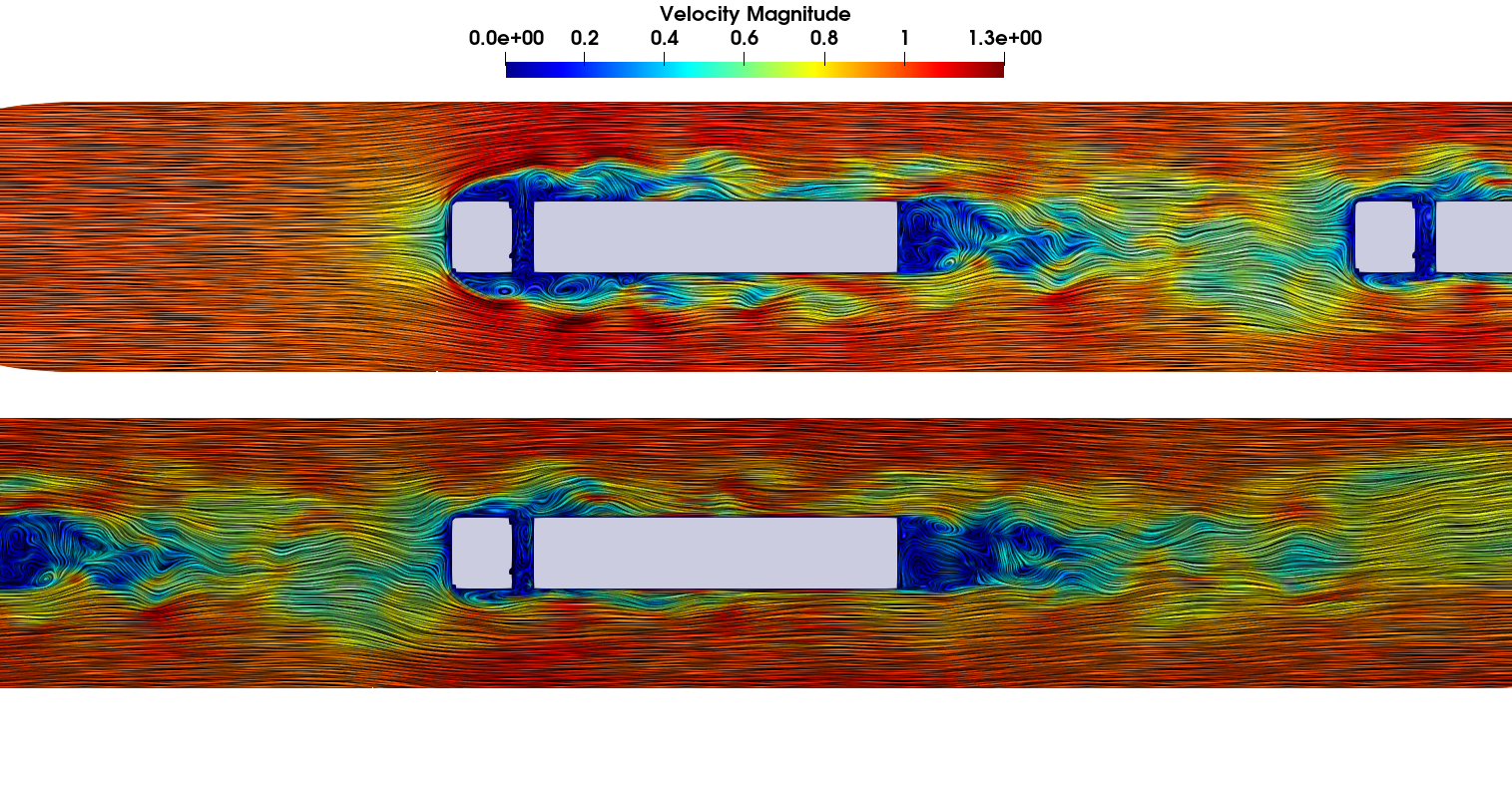}
  \vspace{-5 mm}
  \caption{Flow structures at $T = 0.25$}
  \label{fig: TwoTruckFlowStructure_2}
\end{subfigure}
\\ 
\begin{subfigure}{.48\textwidth}
\vspace{5 mm}
  \includegraphics[width=0.9\linewidth]{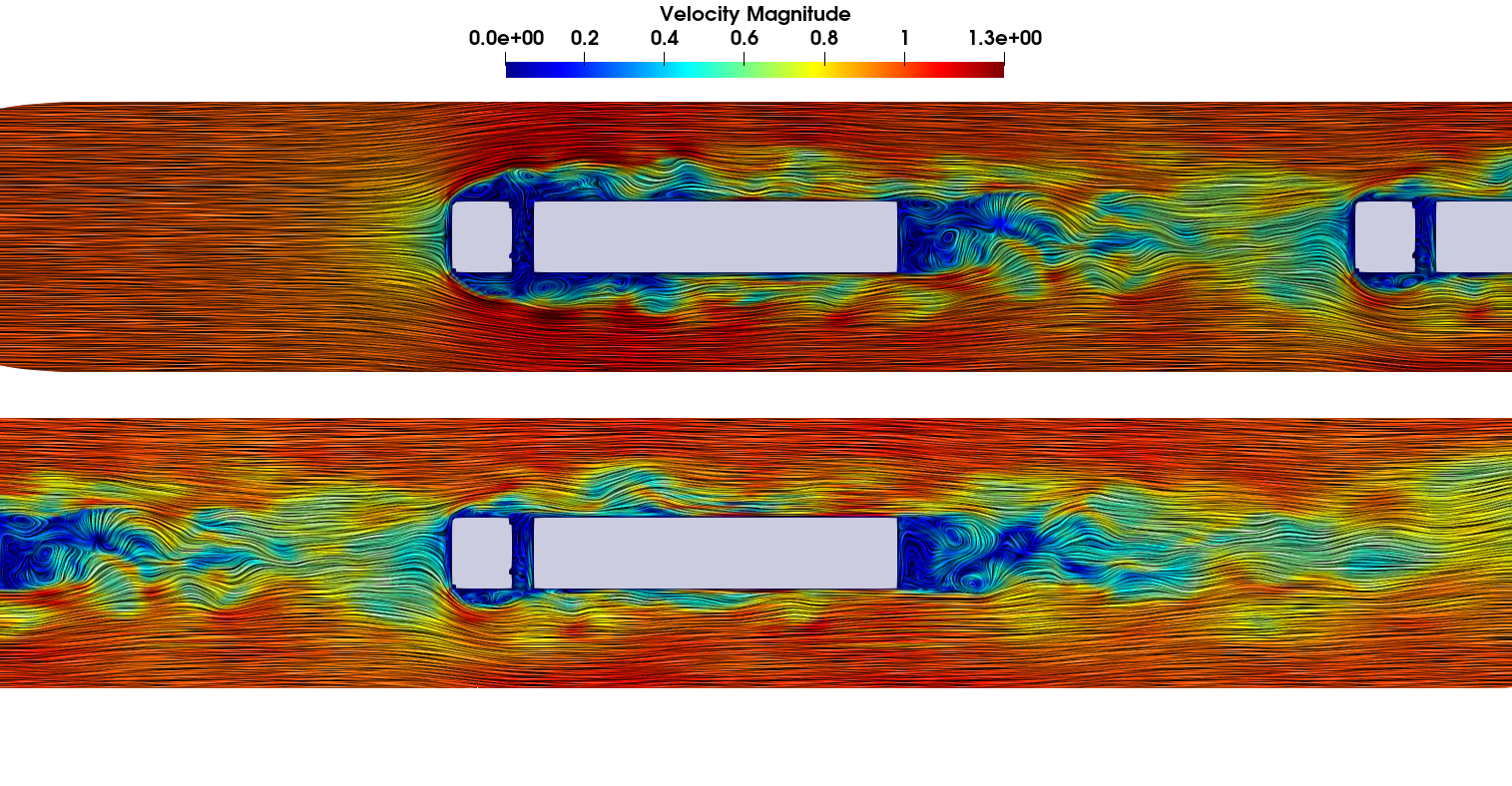}
  \vspace{-5 mm}
  \caption{Flow structures at $T = 0.5$}
  \label{fig: TwoTruckFlowStructure_3}
\end{subfigure}
\begin{subfigure}{.48\textwidth}
\vspace{5 mm}
  \includegraphics[width=0.9\linewidth]{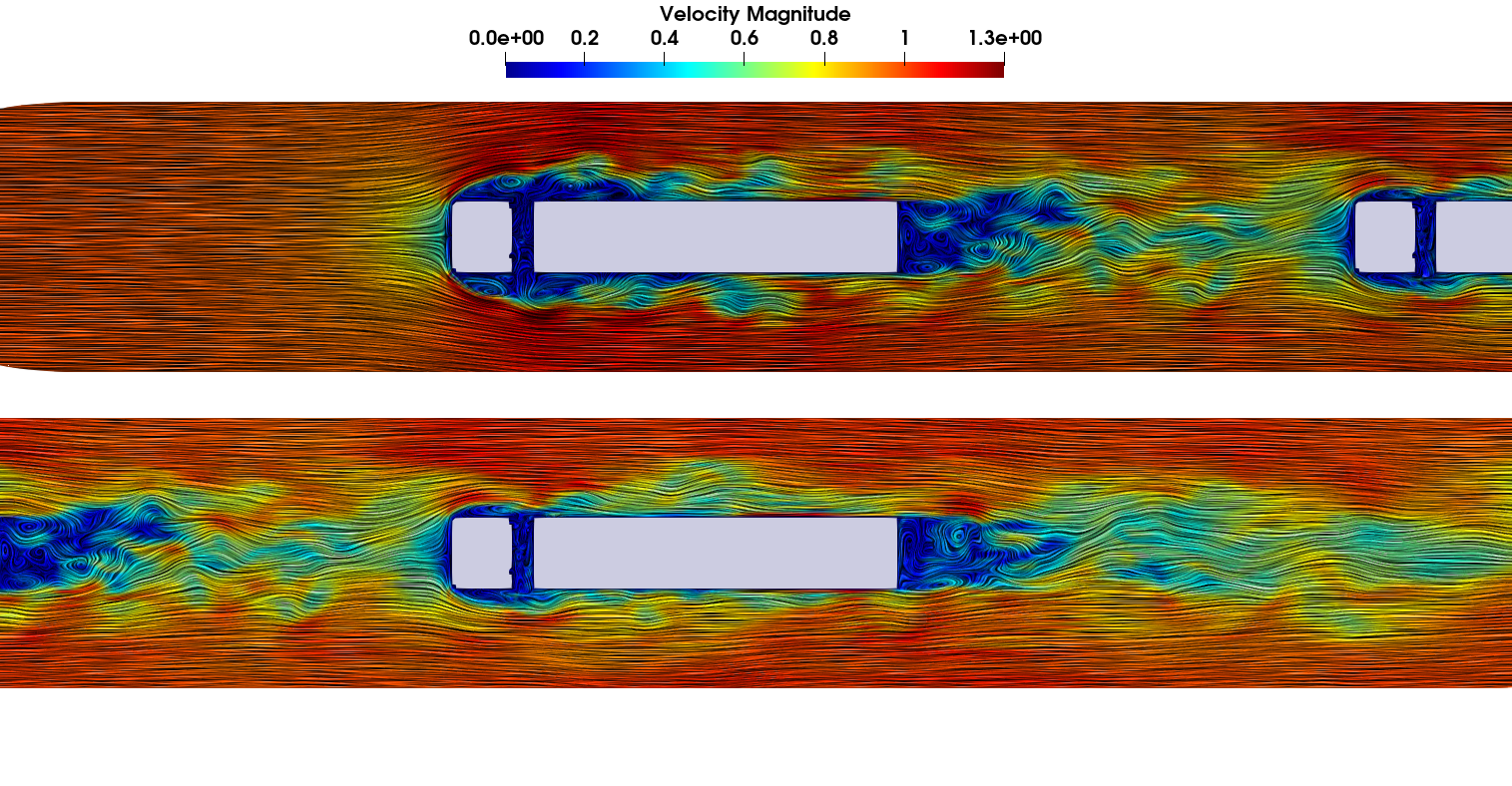}
  \vspace{-5 mm}
  \caption{Flow structures at $T = 0.75$}
  \label{fig: TwoTruckFlowStructure_4}
\end{subfigure}
\\
\begin{subfigure}{.48\textwidth}
\vspace{5 mm}
  \includegraphics[width=0.9\linewidth]{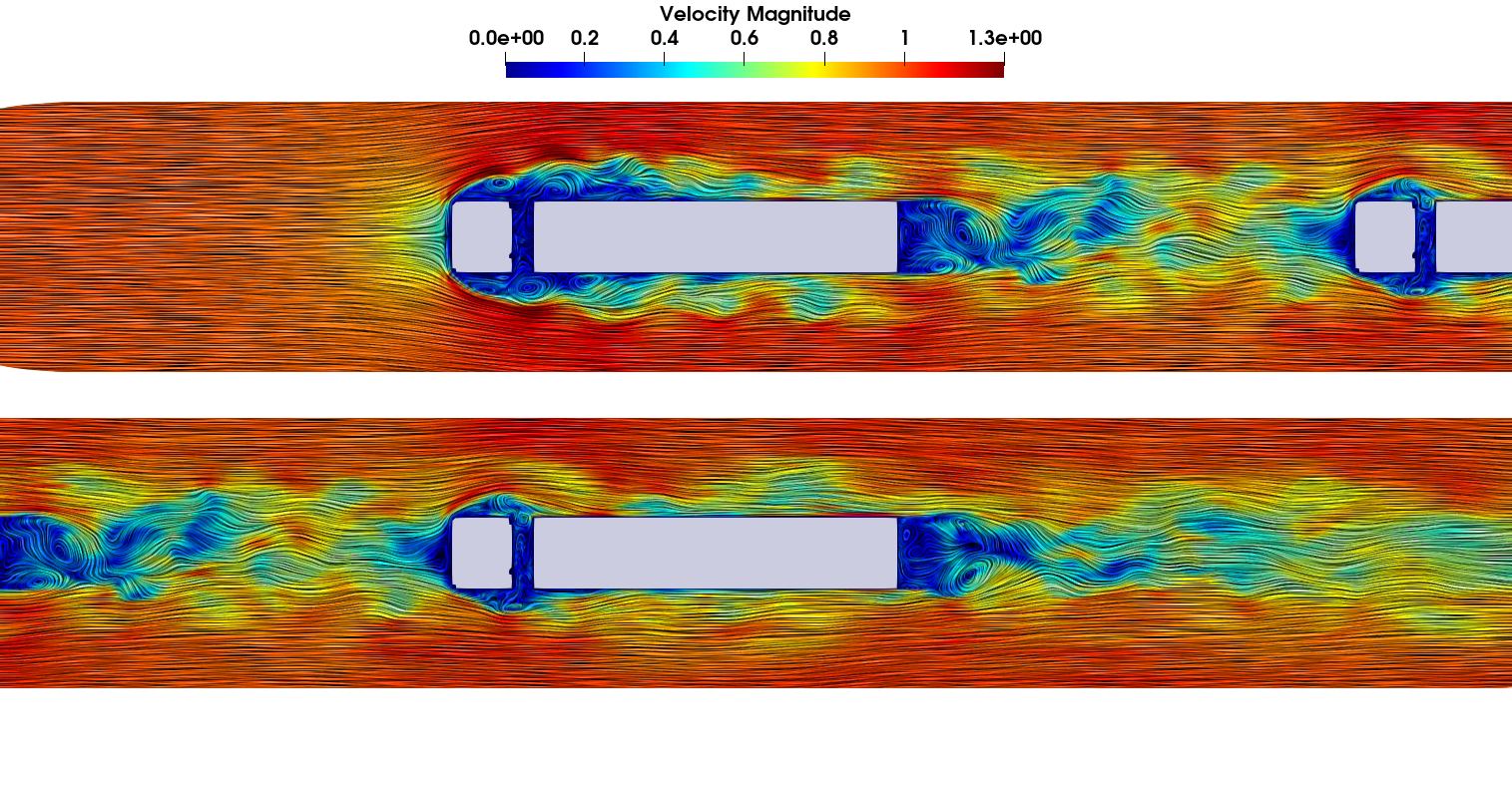}
  \vspace{-5 mm}
  \caption{Flow structures at $T = 1.0$}
  \label{fig: TwoTruckFlowStructure_5}
\end{subfigure}
\begin{subfigure}{.48\textwidth}
\vspace{5 mm}
  \includegraphics[width=0.9\linewidth]{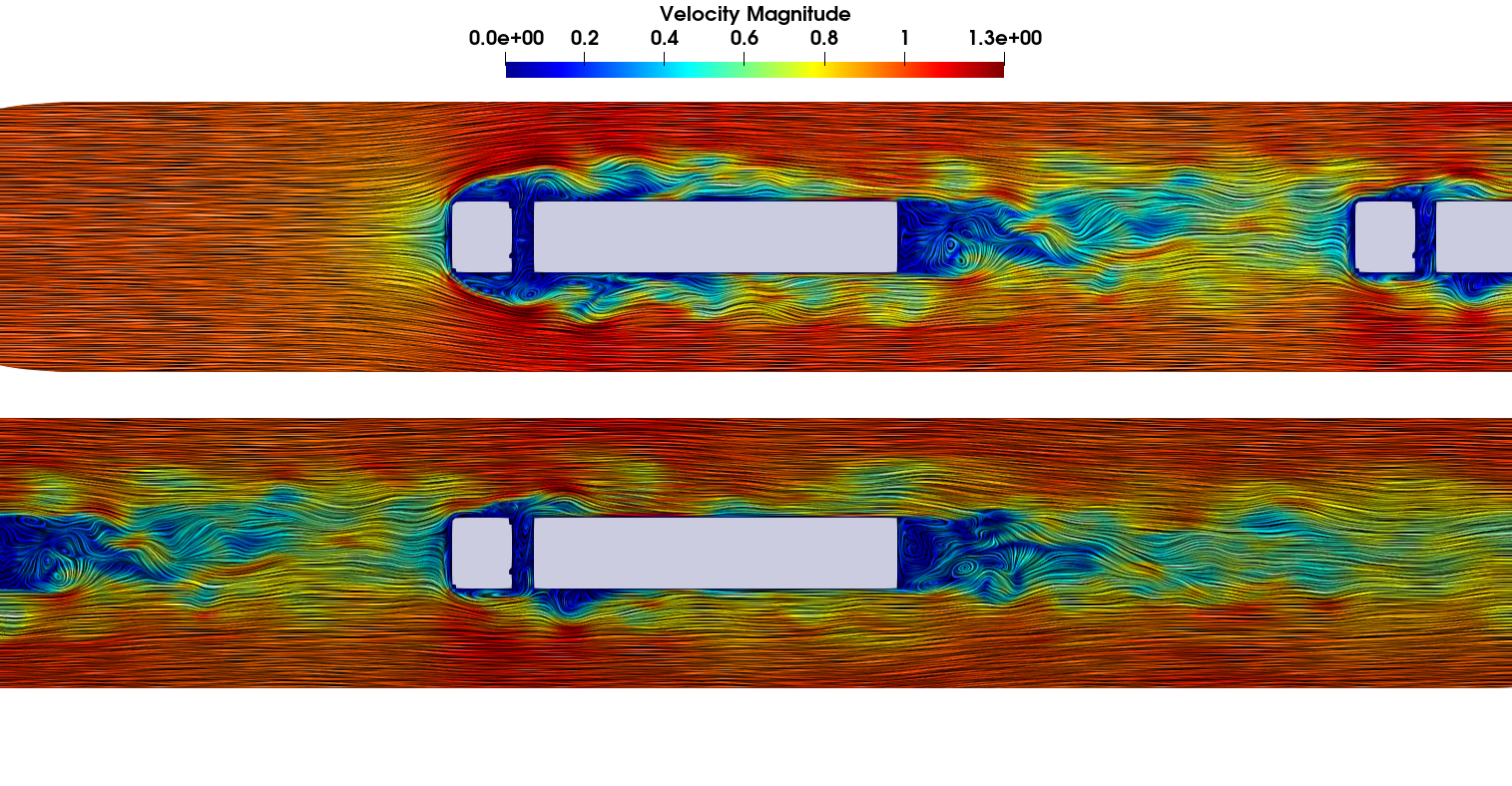}
  \vspace{-5 mm}
  \caption{Flow structures at $T = 1.25$}
  \label{fig: TwoTruckFlowStructure_6}
\end{subfigure}
\vspace{0.1in}
\caption{Flow structures for two platooning trucks in top-down view.}
\label{fig: TwoTruckTopDownFlowStructure}
\end{figure*}